\documentclass[12pt]{amsart} 

\usepackage[margin = 1in]{geometry}
\usepackage{amsfonts, amsmath,amscd, amssymb, amsbsy, cite, color, latexsym, mathrsfs, mathtools,
slashed, stmaryrd, verbatim, wasysym }
\usepackage[all]{xy}
\usepackage[pdftex]{graphicx}
\usepackage{hyperref}
\usepackage[applemac]{inputenc}
\usepackage[titletoc,title]{appendix}
\usepackage{tikz-cd}
\usepackage{oubraces}
\usepackage{cancel}
\usepackage{float}
\allowdisplaybreaks

\usepackage{import}
\usepackage{xifthen}
\usepackage{pdfpages}
\usepackage{transparent}

\usepackage{tikz}
\usetikzlibrary{decorations.pathmorphing, matrix, calc, arrows}

\newlength{\LETTERheight}
\AtBeginDocument{\settoheight{\LETTERheight}{I}}

\newtheorem{theorem}{Theorem}
\newtheorem{corollary}[theorem]{Corollary}

\newtheorem{remark}{Remark}

\numberwithin{equation}{section}
\numberwithin{theorem}{section}

\let\oldsqrt\sqrt
\def\sqrt{\mathpalette\DHLhksqrt}
\def\DHLhksqrt#1#2{%
\setbox0=\hbox{$#1\oldsqrt{#2\,}$}\dimen0=\ht0
\advance\dimen0-0.2\ht0
\setbox2=\hbox{\vrule height\ht0 depth -\dimen0}%
{\box0\lower0.4pt\box2}}

\newcommand{\bhs}[1]{\mathfrak B_{#1}}

\renewcommand{\bar}{\overline}
\renewcommand{\Re}{\operatorname{Re}}

\newcommand{\wt}[1]{\widetilde{#1}}
\newcommand{\rest}[1]{\big\rvert_{#1}} 

\DeclareMathOperator{\hook}{\lrcorner}

\newcommand{\Rp}{\mathbb{R}^+}

\newcommand\lra{\longrightarrow}
\newcommand\xlra[1]{\xrightarrow{\phantom{x} #1 \phantom{x}}}

\newcommand\pa{\partial}

\newcommand\eps\varepsilon
\renewcommand\epsilon\varepsilon

\newcommand{\sR}{\mathrm{sR}}

\newcommand\dCI{\dot{\mathcal{C}}^{\infty}}

\newcommand\CI{{\mathcal{C}}^{\infty}}

\newcommand{\lrpar}[1]{\left( #1 \right)}

\newcommand\ang[1]{\langle #1 \rangle}
\newcommand{\lrbrac}[1]{\left\lbrace #1 \right\rbrace}

\newcommand{\pmat}[1]{\begin{pmatrix} #1 \end{pmatrix}}

\newcommand\Ann{\operatorname{Ann}}

\renewcommand\det{\operatorname{det}}

\newcommand\diag{\operatorname{diag}}
\newcommand\Diff{\operatorname{Diff}}

\newcommand\End{\operatorname{End}}

\DeclareMathOperator*{\FP}{\operatorname{FP}}

\newcommand\Hom{\operatorname{Hom}}

\newcommand\Id{\operatorname{Id}}
\renewcommand\Im{\operatorname{Im}}

\newcommand\Ker{\operatorname{Ker}}

\newcommand\phg{\operatorname{phg}}

\renewcommand\Re{\operatorname{Re}}

\newcommand\Tr{\operatorname{Tr}}

\newcommand\tr{\operatorname{tr}}

\newcommand\Mand{\text{ and }}
\newcommand\Mas{\text{ as }}

\newcommand\Mfor{\text{ for }}
\newcommand\Mforall{\text{ for all }}

\newcommand\Mforsome{\text{ for some }}

\newcommand\Mif{\text{ if }}

\newcommand\Mon{\text{ on }}

\newcommand\Mst{\text{ s.t. }}

\newcommand\Mwhere{\text{ where }}
\newcommand\Mwith{\text{ with }}

\newcommand\paperintro%
        {%
         }
\newcommand\paperbody%
        {%
         }

\newcommand\bbC{\mathbb{C}}

\newcommand\bbE{\mathbb{E}}

\newcommand\bbH{\mathbb{H}}

\newcommand\bbN{\mathbb{N}}

\newcommand\bbR{\mathbb{R}}
\newcommand\bbS{\mathbb{S}}

\newcommand\cA{\mathcal{A}}

\newcommand\cD{\mathcal{D}}
\newcommand\cE{\mathcal{E}}
\newcommand\cF{\mathcal{F}}

\newcommand\cI{\mathcal{I}}
\newcommand\cJ{\mathcal{J}}
\newcommand\cK{\mathcal{K}}
\newcommand\cL{\mathcal{L}}

\newcommand\cO{\mathcal{O}}

\newcommand\cR{\mathcal{R}}

\newcommand\cU{\mathcal{U}}
\newcommand\cV{\mathcal{V}}
\newcommand\cW{\mathcal{W}}

\newcommand\sA{\mathscr{A}}

\newcommand\sE{\mathscr{E}}

\newcommand\sG{\mathscr{G}}
\newcommand\sH{\mathscr{H}}
\newcommand\sI{\mathscr{I}}

\newcommand\sK{\mathscr{K}}

\newcommand\sT{\mathscr{T}}

\newcommand\sV{\mathscr{V}}

\newcommand\tH{\operatorname{H}}

\newcommand\tN{\operatorname{N}}

\DeclareMathAlphabet{\mathpzc}{OT1}{pzc}{m}{it}

\renewcommand\mid{\text{mid}}
\renewcommand\nmid{\text{out}}
\newcommand{\density}{\Lambda} 
\newcommand{\midH}{H_{\sR}^{\mid}X}
\newcommand{\nmidH}{H_{\sR}^{\nmid}X}
\newcommand{\pmidH}{H_{\sR}^{+}X}

\begin{document}

\title[Sub-Riemmanian limit of heat kernels]{Sub-Riemannian limit of the differential form heat kernels of contact manifolds}
\author{Pierre Albin}
\address{University of Illinois, Urbana-Champaign}
\email{palbin@illinois.edu}
\author{Hadrian Quan}
\email{hquan4@illinois.edu}

\begin{abstract}
We study the behavior of the heat kernel of the Hodge Laplacian on a contact manifold endowed with a family of Riemannian metrics that blow-up the directions transverse to the contact distribution. We apply this to analyze the behavior of global spectral invariants such as the $\eta$-invariant and the determinant of the Laplacian. In particular we prove that contact versions of the relative $\eta$-invariant and the relative analytic torsion are equal to their Riemannian analogues and hence topological.
\end{abstract}

\maketitle

\paperintro
\section*{Introduction}

Sub-Riemannian geometry is a generalization of Riemannian geometry where distances are measured using curves that are tangent to a subbundle of the tangent bundle.
Perhaps the best studied setting is that of a contact manifold $M$ where the curves are required to be tangent to the contact distribution $\sH \subseteq TM.$
As $\sH$ is maximally non-integrable the Chow-Rashevskii theorem guarantees that any two points on $M$ can be connected by a curve tangent to $\sH$ \cite[Chapter 2]{Montgomery:Tour} and so a bundle metric $g_{\sH},$ known as a Carnot-Caratheodory metric, induces a distance $d_{g_{\sH}}$ on $M$ analogously to the Riemannian distance.

One natural approach to studying this sub-Riemannian geometry is through approximation: if $g_{\eps}$ is a one-parameter family of Riemannian metrics, extending $g_{\sH}$ and blowing-up in the directions transverse to $\sH$ as $\eps\to0,$ it was shown by Gromov \cite{Gromov:CC} that the metric spaces $(M, d_{g_{\eps}})$ converge in the Gromov-Hausdorff sense to $(M, d_{g_{\sH}}).$ The behavior of the Hodge Laplacians of the metrics $g_{\eps}$ was initiated by Rumin \cite{Rumin:Sub}, following work of Ge on the scalar Laplacian \cite{Ge:Collapsing}, who showed that the those parts of the spectrum that have finite limits concentrate on the spectrum of their counterparts in the Rumin contact complex from \cite{Rumin:Formes},
\begin{equation*}
	0 \to \Omega^0_{\sH}M \xlra{d_{\sH}} \cdots \xlra{d_{\sH}} \Omega^n_{\sH}M 
	\xlra{D_{\sH}} 
	\Omega^{n+1}_{\sH}M \xlra{d_{\sH}} \cdots \xlra{d_{\sH}} \Omega^m_{\sH}M \to 0.
\end{equation*}
This complex is built from the de Rham complex and the contact structure of $M;$ here $m=\dim M$ is equal to $2n+1,$ $\Omega_{\sH}^q(M)$ is a quotient of $\Omega^q(M)$ if $q\leq n$ and a subspace otherwise, $d_{\sH}$ is a first order operator and $D_{\sH}$ is a second order operator.
Interestingly this complex computes the singular cohomology of $M,$ just as the de Rham complex does.
Rumin's initial motivation was to obtain differential operators compatible with `Heisenberg dilations' in which the $\sH$ directions are scaled by $\lambda$ and the complementary directions by $\lambda^2.$
Indeed the Hodge Laplacians of the Rumin complex are not elliptic in the usual sense of invertible principal symbols; however they have symbols `Heisenberg calculus' that are invertible and as a result they are hypoelliptic operators.

Rumin points out in \cite{Rumin:Sub} that the contact complex can be derived from a spectral sequence induced, e.g., by Heisenberg dilations. We study this spectral sequence using Hodge theoretic techniques following \cite{Mazzeo-Melrose:Ad, Forman}. In particular we obtain filtrations of the de Rham complex (really of a rescaling of the de Rham complex, see \S\ref{sec:sRLim})
\begin{equation*}
\begin{gathered}
	\sE_{\infty}^p = \sE_4^p \subseteq \sE_2^p \subseteq \sE_0^p = \Omega^p(M) \Mif p \notin\{n, n+1\}, \\
	\Mand \sE_{\infty}^p \subseteq \sE_4^p \subseteq \sE_2^p \subseteq \sE_0^p = \Omega^p(M) \Mif p \in\{n, n+1\},
\end{gathered}
\end{equation*}
in which $\sE_2^p$ is isomorphic to the Rumin complex, $\sE_4^p$ corresponds to the null space of $(d_{\sH}+\delta_{\sH}),$ and $\sE_{\infty}^p$ consists of the Rumin harmonic forms and so is finite dimensional and isomorphic to the singular cohomology $\tH^p(M).$

Rumin's approach in \cite{Rumin:Sub} was to focus on the convergence of the resolvent of the Hodge Laplacians in the sub-Riemannian limit. 
They note that while their approach establishes convergence of appropriate heat kernels for large times, they do not apply to the short-time behavior.
In this article we establish the behavior of the heat kernel of the Hodge Laplacians for all time. Our approach follows Melrose's construction of the heat kernel (see \cite[Chapter 7]{Melrose:APS} and \S\ref{sec:EllHeat} below) in that we construct a manifold with corners on which the heat kernel is essentially smooth.

\begin{theorem}\label{thm:IntroThm}
 Let $M$ be a contact manifold with contact form $\theta,$ let $J$ be an almost contact structure on $\sH=\ker\theta$ such that $g_{\sH}(\cdot, \cdot)= d\theta(\cdot, J\cdot)$ is symmetric positive definite on $\sH,$ and let
\begin{equation}\label{eq:IntMet}
	g_{\eps} = g_{\sH} \oplus \frac{\theta\otimes\theta}{\eps^2}.
\end{equation}

i)
The heat kernels of the Hodge Laplacians $\Delta_{\eps},$ acting on differential forms of degree $p \notin\{n,n+1\}$ are $\cI$-smooth (i.e., polyhomogeneous, see \S\ref{sec:Corners}) on a manifold with corners $\nmidH$ with boundary hypersurfaces capturing the asymptotics as $t\to0$ with $\eps$ bounded, as $t\to 0$ like $\eps^2,$ and as $\eps\to0$ with $t$ bounded.

For differential form degrees $p \in \{n,n+1\},$ the same is true for $\eps^{-2}\Delta_{\eps}$ with boundary hypersurfaces capturing the asymptotics as $t\to0$ with $\eps$ bounded, as $t\to 0$ like $\eps^4,$ as $t\to 0$ like $\eps^2,$ and as $\eps\to0$ with $t$ bounded.\\

ii) For differential forms of degree $p \notin\{n,n+1\},$ the asymptotics of the trace of the heat kernels of the Hodge Laplacians take the form
\begin{equation*}
	\Tr(e^{-t\Delta_{\eps}}) \sim
	\begin{cases}
	t^{-m/2} \sum_{k\geq0} a_k t^k & \Mas t\to 0, \eps>0, \\
	\rho_2^{-(m+1)} \sum_{k\geq 0} A_k \rho_2^k & \Mas \rho_2=\sqrt{t+\eps^2}\to 0,\\
	\Tr(e^{-t(d_{\sH}+d_{\sH}^*)^2}\Pi_{\sE_2^p}) + \cO(\eps) & \Mas \eps\to 0, t>0,
	\end{cases}
\end{equation*}
where $a_k$ and $A_k$ are local (i.e., they are integrals of universal polynomials in the local invariants of the metric), and
$\Pi_{\sE_2^p}$ denotes the orthogonal projection onto $\sE_2^p.$

For differential forms of degree $p \in \{n, n+1\},$ we have
\begin{equation*}
	\Tr(e^{-t\eps^{-2}\Delta_{\eps}}) \sim
	\begin{cases}
	t^{-m/2} \sum_{k\geq0} a_k t^k & \Mas t\to 0, \eps>0, \\
	\rho_4^{-(m+1)} \sum_{k\geq 0} A_k \rho_4^k & \Mas \rho_4=\sqrt{t+\eps^4}\to 0,\\
	\begin{multlined}
	\textstyle
	\Tr(e^{-t(d_{\sH}+d_{\sH}^*)^2}\Pi_{\sE_2^p\setminus \sE_4^p}) \phantom{xxxxxxxxxxx}\\
	\textstyle
	\phantom{xxxxxxxxxxx}+ \rho_2^{-(m+1)} \sum_{k\geq0} B_k \rho_2^k +\cO(\rho_2)
	\end{multlined}
		& \Mas \rho_2=\sqrt{t+\eps^2}\to 0,\\
	\Tr(e^{-t(D_{\sH}+D_{\sH}^*)^2}\Pi_{\sE_4^p}) + \cO(\eps) & \Mas \eps\to 0, t>0,
	\end{cases}
\end{equation*}
where $a_k,$ $A_k,$ and $B_k$ are local, $\Pi_{\sE_2^p\setminus \sE_4^p}$ denotes the orthogonal projection onto $\sE_2^p \setminus \sE_4^p,$
and $\Pi_{\sE_4^p}$ denotes the orthogonal projection onto $\sE_4^p.$\\

iii) If $F$ is a flat bundle on $M$ then the same constructions hold for the heat kernels of the Hodge Laplacians of the de Rham complex with coefficients in $F,$ $\Delta_{\eps}^F.$ In particular, if $F_1$ and $F_2$ are flat bundles of the same rank then 
\begin{equation*}
	\Tr(e^{-t\Delta_{\eps}^{F_1}}) - \Tr(e^{-t\Delta_{\eps}^{F_2}})
\end{equation*}
has the same asymptotics as in (ii) but without the local terms, and similarly for $\eps^{-2}\Delta_{\eps}^{F_i}$ acting on forms in middle degrees. 
\end{theorem}

We then apply this construction to study the $\eta$-invariant and the analytic torsion of contact manifolds.

The $\eta$-invariant of the odd signature operator, $\eta(S),$ was introduced Atiyah, Patodi, and Singer \cite{APS-I} as the boundary contribution to the index formula for the signature operator with appropriate global boundary conditions. In \cite{APS-II} they introduced the $\rho$-invariant which assigns to a flat bundle $F$ the difference of the $\eta$-invariant of the odd signature operator twisted by $F,$ $\eta(S^F),$ with the $\eta$-invariant of the trivial flat bundle of the same rank,
\begin{equation*}
	\rho(F) = \eta(S^F) - \mathrm{rank}(F) \eta(S).
\end{equation*}
They showed that the $\rho$ invariant is a smooth invariant of $M,$ i.e., it is independent of the choice of metric on $M.$
(Note that explicit computations on lens spaces show that it is not a homotopy invariant in general. For manifolds whose fundamental group satisfies the Borel conjecture \cite{Weinberger} or an appropriate version of the Baum-Connes conjecture \cite{Keswani}, the $\rho$-invariant is a homotopy invariant.)

\begin{corollary}
Let $F \lra M$ be a flat bundle, $S_{\eps}^F$ the odd signature operator of $(M,g_{\eps})$ with coefficients in $F,$ and $S_{\sH}^F$ the odd signature operator of Rumin's complex with coefficients in $F.$ As $\eps \to 0,$ the finite part of the difference between $\eta(S_{\eps}^F)$ and $\eta(S_{\sH}^F)$ is the integral of a local invariant of the metric,
\begin{equation}\label{eq:IntEta}
	\eta(S_{\sH}^F) - \FP_{\eps=0} \eta(S_{\eps}^F) = \text{ local }.
\end{equation}
Correspondingly, the $\rho$-invariant of the Rumin complex is equal to the $\rho$-invariant of the de Rham complex,
\begin{equation*}
	\rho(F) = \eta(S_{\sH}^F) - \mathrm{rank}(F) \eta(S_{\sH}).
\end{equation*}
\end{corollary}

The $\eta$-invariant of the Rumin complex was first studied in \cite[\S7]{Rumin:Sub} and subsequently by Biquard, Herzlich, and Rumin \cite{Biquard-Herzlich-Rumin}.
In the latter, they conjecture \cite[\S6]{Biquard-Herzlich-Rumin} that \eqref{eq:IntEta} holds in arbitrary dimension and they study it on three dimensional contact manifolds. If  $M$ is a three dimensional Cauchy-Riemann (CR) Seifert manifold, meaning that the almost complex structure $J$ in the metric \eqref{eq:IntMet} is integrable and $M$ is endowed with a locally free action of $\bbS^1$ that preserves $(\sH, J)$ and is generated by a Reeb field, then they establish \eqref{eq:IntEta} with an explicit formula for the local term, \cite[Theorem 1.4]{Biquard-Herzlich-Rumin}
\begin{equation*}
	\text{$M^3$ CR Seifert} \implies \eta(S_{\sH}) - \FP_{\eps=0} \eta(S_{\eps}) = -\frac1{512}\int_M R^2 \; \theta\wedge d\theta,
\end{equation*}
where $\theta$ is an $\bbS^1$-invariant contact form and $R$ is the curvature of the Tanno-Webster-Tanaka connection.

The analytic torsion of $M$ is a function that uses determinants of the Hodge Laplacians of the de Rham complex to assign to each flat bundle $F \lra M$ and each basis of the cohomology of $M$ with coefficients in $F$ a real number (see \S\ref{sec:Det}).
It was defined by Ray and Singer \cite{Ray-Singer:AT} as an analytic analogue of a combinatorial invariant known as Reidemeister torsion which it was shown to equal by Cheeger \cite{Cheeger:AT} and M\"uller \cite{Muller:AT}, independently, for flat bundles with unitary holonomy, and by M\"uller \cite{Muller:Uni} for flat bundles with unimodular holonomy. Bismut and Zhang \cite{Bismut-Zhang} were able to extend the Cheeger-M\"uller theorem to arbitrary flat bundles.

\begin{corollary}\label{cor:IntAT}
Let $F\lra M$ be a flat bundle with unimodular holonomy, $\{\mu^*\}\subseteq \tH^*(M;F)$ a basis for the cohomology of $M$ with coefficients in $F.$ The difference between the logarithms of the analytic torsion of the de Rham complex of $M$ and the analytic torsion of the Rumin complex of $M$ is given by the integral of a local invariant of the metric,
\begin{equation}\label{eq:IntAT}
	\log \bar{AT}(M,\{\mu^*\}, F) - \log \bar{AT}_{\sH}(M, \{\mu^*\}, F) = \text{ local }.
\end{equation}
In particular, if $F_1$ and $F_2$ are flat bundles of the same rank with isomorphic cohomology (e.g., if they are both acyclic) then the relative analytic torsions of the de Rham and Rumin complexes coincide,
\begin{equation*}
	\log \bar{AT}(M,\{\mu^*\}, F_1) - \log \bar{AT}(M,\{\mu^*\}, F_2)
	= \log \bar{AT}_{\sH}(M,\{\mu^*\}, F_1) - \log \bar{AT}_{\sH}(M,\{\mu^*\}, F_2)
\end{equation*}
\end{corollary}

The analytic torsion of the Rumin complex was introduced and studied by Rumin and Seshadri \cite{Rumin-Seshadri}.
For three dimensional CR Seifert manifolds they are able to establish a stronger version of \eqref{eq:IntAT}, \cite[Theorem 4.2]{Rumin-Seshadri}
\begin{equation*}
	\text{$M^3$ CR Seifert} \implies \log \bar{AT}(M,\{\mu^*\}, F) = \log \bar{AT}_{\sH}(M, \{\mu^*\}, F).
\end{equation*}
On a general $3$-dimensional contact manifold they show that modifying $\log \bar{AT}_{\sH}(M, \{\mu^*\}, F)$ by a local term produces a `CR-torsion,' $\log \bar{AT}_{CR}(M, \{\mu^*\}, F),$ which is independent of the contact form. As, in three dimensions, there are no local invariants of contact metrics that are independent of the choice of contact form 
\cite[proof of Theorem 9.1]{Biquard-Herzlich-Rumin}, it follows from our theorem that
\begin{equation*}
	\text{$M^3$ Contact} \implies \log \bar{AT}(M,\{\mu^*\}, F) = \log \bar{AT}_{CR}(M, \{\mu^*\}, F).
\end{equation*}

In an interesting preprint \cite{Kitaoka}, Kitaoka looks at a modified Rumin complex with the differential $d_{\sH}$ replaced by
\begin{equation*}
	d_{\sK} = |n-p|^{-1/2}d_{\sH} \text{ on forms of degree $p,$ for all }p \notin \{n,n+1\}
\end{equation*}
and studies the corresponding analytic torsion on the odd-dimensional spheres with the standard contact structure and metric, $\log \bar{AT}_{\sK}(\bbS^{2n+1}, \{\mu^*\}, \bbR).$
Using representation theory they are able to prove that
\begin{equation*}
	\log \bar{AT}_{\sK}(\bbS^{2n+1}, \{\mu^*\}, \bbR) = \log \bar{AT}(\bbS^{2n+1}, \{\mu^*\}, \bbR) + \log (n!).
\end{equation*}
A direct computation (see \eqref{eq:KitaokaRS} below) shows that in general Kitaoka's torsion differs from that of Rumin-Seshadri by a local term
\begin{equation}\label{eq:IntKitaokaRS}
	\log \bar{AT}_{\sK}(M,\{\mu^*\}, F) - \log \bar{AT}_{\sH}(M, \{\mu^*\}, F) = \text{ local }.
\end{equation}
$ $

{\bf Further context.}
Pseudodifferential operators on the Heisenberg group, adapted to parabolic dilations, were introduced by Dynin \cite{Dynin}.
The calculus of Heisenberg pseudodifferential operators was consequently developed by Beals-Greiner \cite{Beals-Greiner} and Taylor \cite{Taylor:Noncom} (see also \cite{Melrose:ScatPs} for the relation with the $\theta$-calculus of \cite{Epstein-Melrose-Mendoza}).
A calculus including both the usual pseudodifferential operators and the Heisenberg pseudodifferential operators was developed by Epstein and Melrose \cite{Epstein-Melrose:HHH} and used to study index theory by van Erp \cite{vanErp:Ext} (see also \cite{vanErp:ASI, Baum-vanErp, Epstein:Lect} for more on index theory in this context).
Complex powers and noncommutative residues of operators in the Heisenberg calculus were developed by Ponge \cite{Ponge:Spectral} as well as a construction of associated heat kernels following work of Beals, Greiner, and Stanton \cite{Beals-Greiner-Stanton}. Recently the relations between sub-Riemannian spectral asymptotics and dynamics have been studied by, e.g., Colin de Verdi\`ere, Hillairet, and Tr\'elat \cite{CdV-Hillaret-Trelat}, Fermanian and Fischer \cite{Fermanian-Fischer}, and Savale \cite{Savale}.
Some of these constructions have now been extended to more general filtered manifolds. For example, the groupoid approach to pseudodifferential operators in \cite{vanErp-Yuncken} allows van Erp and Yuncken to handle filtered manifolds, and has been used by Dave and Haller to study BGG sequences generalizing the Rumin complex \cite{Dave-Haller:Graded} and to study heat kernels and their asymptotics \cite{Dave-Haller:Heat}.

The sub-Riemannian limit is closely related to the adiabatic limit introduced by Witten \cite{Witten} in which the metric on the total space of a fiber bundle is blown-up along the fibers. Witten's results were treated and generalized in \cite{Cheeger:Ad, Bismut-Freed:I, Bismut-Freed:II} and since then the adiabatic limit has become a standard tool in geometric analysis. Most relevant to our approach is the adiabatic limit approach to the Leray-Serre spectral sequence by Mazzeo and Melrose \cite{Mazzeo-Melrose:Ad}, and its subsequent extension by Forman \cite{Forman}, as well as the adiabatic limit of analytic torsion studied by Dai and Melrose \cite{Dai-Melrose}. Similar analytic surgery techniques have been used, e.g.,  to study the formation of cylinders \cite{Mazzeo-Melrose:Surgery, Hassell-Mazzeo-Melrose:Surgery, Hassell:AT}, fibered cusps \cite{Albin-Rochon-Sher:FibCusps, Albin-Rochon-Sher:Cusps}, and wedges \cite{Albin-Rochon-Sher:Wedge}. In contrast to these works the degeneration in this project is not happening along a submanifold but rather along a sub-bundle of the tangent bundle.

$ $\\
{\bf Acknowledgements.} 
This material is based upon work supported by the Simons Foundation under grant \#317883 and
the National Science Foundation under grant DMS-1711325 of the first author, 
a graduate research fellowship of the second author, 
and grant DMS-1440140 while the authors were in residence at the Mathematical Sciences Research Institute in Berkeley, California during the Fall 2019 semester.
This project grew out of conversations of the first author with Jing Wang and both authors are grateful for her help. 
The first author is happy to acknowledge helpful conversations with Rafe Mazzeo and Richard Melrose.
The second author is happy to acknowledge helpful conversations with Davide Barilari, Charlie Epstein, and Jeremy Tyson.

\tableofcontents

\paperbody
\section{Rumin's contact complex}\label{sec:Contact Complex}

In this section we review some results of \cite{Rumin:Formes, Rumin:Sub}.

Let $M$ be a contact manifold of dimension $m=2n+1$ with contact distribution $\sH \subseteq TM.$
Let $\mathrm{Ann} (\sH)\subseteq T^*M$ denote the annihilator of $\sH.$
Let $\cI^*$ denote the differential ideal generated by the sections of $\mathrm{Ann} (\sH),$
\begin{equation*}
	\cI^q = \mathrm{span} \{ \theta \wedge \alpha + d\theta \wedge\beta: 
	\alpha \in \Omega^{q-1}M, \beta \in \Omega^{q-2}M, \theta \in \CI(M, \mathrm{Ann}(\sH))\},
\end{equation*}
and let $\cJ^*$ denote the forms that vanish after wedge product with an element of $\cI^*,$
\begin{equation*}
	\cJ^q = \{\omega \in \Omega^qM: \theta\wedge\omega = 0 = d\theta\wedge\omega \Mforall \theta \in \CI(M, \mathrm{Ann}(\sH)) \}.
\end{equation*}

Using techniques from K\"ahler geometry  \cite[\S2]{Rumin:Formes} it is easy to show that
\begin{equation*}
	q\leq n \implies \cJ^q=0, \quad
	q>n \implies \cI^q=\Omega^q M.
\end{equation*}
We define
\begin{equation*}
	\Omega_{\sH}^qM = 
	\begin{cases}
	\Omega^qM/\cI^q & \Mif q\leq n \\
	\cJ^q & \Mif q >n
	\end{cases}
\end{equation*}
The exterior derivative $d$ induces two complexes $(\Omega^q_{\sH}M,d_{\sH})_{q\leq n}$ and $(\Omega^q_{\sH}M,d_{\sH})_{q> n}.$
Rumin 
showed that any $\alpha \in \Omega^n_{\sH}M$ has a unique representative $\beta \in \Omega^nM$ such that 
$d\beta\in \cJ^{n+1}=\Omega_{\sH}^{n+1}M,$
and that setting $D_{\sH}(\alpha)=d\beta$ yields a complex
\begin{equation*}
\xymatrix{
	0 \ar[r] & \Omega^0M \ar[r]^-{d_{\sH}} & \Omega^1_{\sH}M \ar[r]^-{d_{\sH}} & \cdots \ar[r]^-{d_{\sH}} & \Omega^n_{\sH}M 
	\ar `r/8pt[d] `/13pt[l] `^dl[lll]_{D_{\sH}} `^r/8pt[dll] [dlll] \\
	& \Omega^{n+1}_{\sH}M \ar[r]^-{d_{\sH}} & \Omega^{n+2}_{\sH}M \ar[r]^-{d_{\sH}} & \cdots \ar[r]^-{d_{\sH}} & \Omega^m_{\sH}M \ar[r] & 0},
\end{equation*}
now known as Rumin's contact complex.
Moreover, he showed that this complex 
is induced by a complex of sheaves resolving the sheaf of locally constant $\bbR$-valued functions and hence that its cohomology coincides with the de Rham cohomology of $M.$
As pointed out by Rumin-Seshadri \cite{Rumin-Seshadri}, these arguments are local and essentially unchanged by allowing coefficients in a flat vector bundle.\\

For a more geometric description of the Rumin complex, let us assume for simplicity that $\sH$ is coorientable.
We fix a global contact form $\theta,$ let $\cR$ be the Reeb vector field determined by
\begin{equation*}
	\theta(\cR)=1, \quad \cR \hook d\theta =0,
\end{equation*}
and let $\sV$ denote the rank one subbundle of $TM$ spanned by $\cR$ so that $TM=\sH \oplus \sV.$
It is always possible to find an almost complex structure $J$ on $\sH$ such that $g_{\sH}(\cdot, \cdot) = d\theta(\cdot, J(\cdot))$ is a symmetric positive definite bundle metric on $\sH.$
We then fix a Riemannian metric on $M,$
\begin{equation}\label{eq:g1}
	g_1 = g_{\sH} \oplus \theta \otimes \theta.
\end{equation}

The splitting of $TM$ induces a splitting of the cotangent bundle, $T^*M = \sH^* \oplus \sV^*,$ and then of differential forms
\begin{equation}\label{eq:SplittingM}
	\Omega^qM = \Omega^{q}\sH^* \oplus \theta \wedge \Omega^{q-1}\sH^* = \Omega^{0,q}M \oplus \Omega^{1,q}M.
\end{equation}
The exterior derivative $d: \Omega^qM \lra \Omega^{q+1}M$ decomposes into 
\begin{equation*}
	d = d^{1,0} + d^{0,1} + d^{-1,2}, \Mwhere d^{j,k}: \Omega^{p,q}M \lra \Omega^{p+j,q+k}M.
\end{equation*}
The last term reflects the non-integrability of $\sH;$ concretely the fact that $\theta \in \Omega^{1,0}M$ and $d\theta \in \Omega^{0,2}M.$ 
We can also write
\begin{equation}\label{eq:dM}
	d = 
	\pmat{
	d_H & L \\ \cL_{\cR} & -d_H }: 
	\Omega^{q}\sH^* \oplus \theta \wedge \Omega^{q-1}\sH^*
	\lra
	\Omega^{q+1}\sH^* \oplus \theta \wedge \Omega^{q}\sH^*
\end{equation}
where $d_H = d^{0,1},$ $L$ is exterior product with $d\theta,$ and $\cL_{\cR}$ denotes the Lie derivative by $\cR.$\footnote{
Since $d_\eps^2=0$ these operators satisfy the relations $d_H^2=-L\cL_{\cR}$ and $[d_H,L]=[d_H,\cL_{\cR}]=[L,\cL_{\cR}]=0.$}
Correspondingly, we have
\begin{equation}\label{eq:deltaM}
	\delta = 
	\pmat{
	\delta_H & \cL_{\cR}^* \\ L^* & -\delta_H }: 
	\Omega^{q}\sH^* \oplus \theta \wedge \Omega^{q-1}\sH^*
	\lra
	\Omega^{q+1}\sH^* \oplus \theta \wedge \Omega^{q}\sH^*
\end{equation}

We can identify
\begin{equation}\label{eq:RumCmpxDesc}
	\Omega^qM/\cI^q \cong \{ \omega \in \Omega^q\sH^*: L^* \omega=0\}, \quad
	\cJ^q = \{ \omega \in \theta\wedge\Omega^{q-1}\sH^*: L\omega =0\}. 
\end{equation}
Since $L,$ as a map $\Omega^q\sH^* \lra \Omega^{q+2}\sH^*,$ is injective if $q\leq n-1$ and surjective if $q \geq n-1$ \cite[\S2]{Rumin:Formes}
we can express this succinctly as
\begin{equation*}
	\Omega_{\sH}^q M \cong \Omega^qM \cap \Ker \pmat{ 0 & L \\ L^* & 0 }.
\end{equation*}
With this identification, the differential $D_{\sH}:\Omega_{\sH}^nM \lra \Omega_{\sH}^{n+1}M$ is given by \cite[(2)]{Rumin:Sub}
\begin{equation*}
	D_{\sH}\alpha = \theta \wedge (\cL_{\cR} + d_H L^{-1} d_H)\alpha.
\end{equation*}

One can also view the contact complex as arising naturally from the spectral sequence induced by the filtration $\sH \subseteq TM,$ see \cite{Rumin:Sub}.
An analytic approach to spectral sequences for fiber bundles was developed by Mazzeo and Melrose \cite{Mazzeo-Melrose:Ad} and then for arbitrary splittings of the tangent bundle by Forman \cite{Forman}.
The Mazzeo-Melrose approach to the contact complex will be crucial to our construction of the heat kernel, so we will describe it in detail in \S\ref{sec:sRLim}.\\

Because $D_{\sH}$ is a differential operator of order two, the Hodge Laplacians of the Rumin complex are defined to be
\begin{equation}\label{eq:DeltasH}
	\Delta_{\sH}|_{\Omega_{\sH}^pM} =
	\begin{cases}
	d_{\sH}\delta_{\sH} + \delta_{\sH}d_{\sH} & \Mif p \notin\{n,n+1\}\\
	(d_{\sH}\delta_{\sH})^2 + D_{\sH}D_{\sH}^* & \Mif p =n\\
	(\delta_{\sH} d_{\sH})^2 + D_{\sH}^*D_{\sH} & \Mif p=n+1
	\end{cases}
\end{equation}
Rumin \cite[\S3]{Rumin:Formes} showed that these differential operators are not elliptic but they are hypoelliptic.
Recall that this means that if $u$ is a distributional differential form in the Rumin complex, i.e., a distributional section of the corresponding bundle, and it satisfies that $\Delta_{\sH}u$ is smooth, then we can conclude that $u$ was itself smooth.
A discussion of how this hypoellipticity is a general feature of differential form complexes on filtered manifolds can be found in a recent paper of Dave and Haller \cite[Theorem 2]{Dave-Haller:Graded}.

\section{The sub-Riemannian limit} \label{sec:sRLim}

As above, let $M$ be a contact manifold of dimension $m=2n+1$ endowed with a global contact form $\theta.$
We extend the metric $g_1$ from \eqref{eq:g1} to a one-parameter family of Riemannian metrics on $M,$
\begin{equation*}
	g_{\eps} = g_{\sH} \oplus \frac{\theta \otimes \theta}{\eps^2}, \quad \eps >0.
\end{equation*}
We refer to the limit as $\eps \to 0$ as the sub-Riemannian limit.
Gromov \cite{Gromov:CC} has pointed out that the metric geometry of the Riemannian manifolds $(M,g_{\eps})$ converges in the Gromov-Hausdorff topology to the metric geometry of the sub-Riemannian manifold $(M, g_{\sH}).$

The family of dual metrics on $T^*M$ converges to a degenerate metric supported on $\sH^*,$ 
and correspondingly the scalar Laplacians of $g_{\eps}$ converge as $\eps\to0$ to the hypoelliptic Kohn Laplacian $\Delta_{\sH},$
\begin{equation*}
	\Delta_{\sH} = -\sum {H_j}^2,
\end{equation*}
where $\{H_j\}$ is a $g_{\sH}$-orthonormal local basis of $\sH.$
$ $

Let $X = M \times [0,1]_{\eps}$ with projection $\pi_{\eps}: X \lra [0,1]_{\eps}.$
The {\bf sub-Riemannian limit vector fields} are
\begin{equation*}
	\cV_{\sR} = \{ W \in \CI(X;TX): (\pi_{\eps})_*W =0 \Mand
	W|_{\eps=0} \in \CI(M; \sH) \},
\end{equation*}
i.e., they are vector fields on $M$ parametrized by $\eps$ and horizontal at $\eps=0.$ 
As $\cV_{\sR}$ is a finitely generated projective module over $\CI(X),$ we can use the 
Serre-Swan theorem, or proceed directly as in \cite[\S8.2]{Melrose:APS}, 
and find that there is a vector bundle
\begin{equation*}
	{}^{\sR}TX \lra X,
\end{equation*}
together with a bundle map $j:{}^{\sR}TX \lra TX$ such that 
$j_*\CI(X;{}^{\sR}TX) = \cV_{\sR} \subseteq \CI(X;TX).$
Eliding the map $j,$ we say that the space of sections of ${}^{\sR}TX$ is $\cV_{\sR}.$
We will refer to ${}^{\sR}TX$ as the {\bf sub-Riemannian limit tangent bundle}.

At any $p \in M$ there is a Darboux coordinate chart 
$(x_1, \ldots, x_n, y_1, \ldots, y_n, z)$ in which $\theta$ has the form 
\begin{equation*}
	\theta = dz +\frac12 \sum (y_j \; dx_j -x_j \; dy_j).
\end{equation*}
In Darboux coordinates the sections of ${}^{\sR}TX$ are locally spanned by
\begin{equation*}
	\pa_{x_j}, \quad
	\pa_{y_j}, \quad
	\eps \pa_z
\end{equation*}
(but we point out that $\pa_{x_j},$ $\pa_{y_j}$ are not horizontal away from the center of the coordinate chart).
Note that $\eps \pa_z$ does not vanish at $\eps=0$ as a section of ${}^{\sR}TX,$ though it does as a section of $TX.$
Correspondingly the sub-Riemannian family of metrics $g_{\eps}$ defines a non-degenerate bundle metric on the sub-Riemannian tangent bundle.
$ $\\

Let ${}^{\sR}T^*X$ denote the sub-Riemannian cotangent bundle, defined as the dual bundle to ${}^{\sR}TX$, and locally spanned (in a Darboux chart) by 
\begin{equation*}
	dx_j, \quad dy_j, \quad \frac{\theta}{\eps} = \frac{dz + \tfrac{1}{2} \sum (y_j dx_j - x_j dy_j)}{\eps},
\end{equation*}
and define the sub-Riemannian limit differential forms to be sections of its exterior powers 
\begin{equation*}
	{}^{\sR}\Omega^p(X) = \CI(X; \Lambda^p( {}^{\sR}T^*X )) 
	= \CI\left(X; \pi_M^*(\Lambda^p \sH^*)\oplus \tfrac{\theta}{\eps}\wedge \pi_M^*(\Lambda^{p-1} \sH^*)\right),
\end{equation*}
where $\pi_M: X \lra M$ is the natural projection.

The exterior derivative on $M$ (or rather the $\pi_\eps$-vertical exterior derivative on $X$) induces a singular differential operator on ${}^{\sR}\Omega^*(X)$ which, with notation as in \eqref{eq:dM}, has the form
\begin{equation*}
	d_\eps=\begin{pmatrix} d_H & \tfrac{1}{\eps}L \\ \eps \cL_{\cR} & -d_H \end{pmatrix}
\end{equation*}
for $\eps>0.$ The adjoint operator is given by a modified form of \eqref{eq:deltaM},
\begin{equation*}
	\delta_{\eps} = \pmat{\delta_H & \eps\cL^*_{\cR} \\ \tfrac1\eps L^* & -\delta_H}
\end{equation*}
and hence the Hodge Laplacian is given by
\begin{equation*}
	\Delta_{\eps} = \pmat{
	d_H\delta_H + \delta_Hd_H + \eps^2\cL_{\cR}^*\cL_{\cR} + \tfrac1{\eps^{2}}LL^* & \tfrac1\eps[\delta_H, L] + \eps [d_H, \cL_R^*] \\
	\tfrac1\eps[L^*, d_H] + \eps[\cL_R, \delta_H] & d_H\delta_H + \delta_Hd_H + \eps^2\cL_{\cR}^*\cL_{\cR} + \tfrac1{\eps^{2}}LL^*
	}.
\end{equation*}

We will make use of a slightly modified de Rham operator,
\[ d_\eps - \delta_\eps =\begin{pmatrix} d_H-\delta_H & \tfrac{1}{\eps}L - \eps \cL_{\pa_z}^* \\ -\tfrac{1}{\eps}L^* + \eps \cL_{\pa_z} & -(d_H-\delta_H) \end{pmatrix},   \]
which satisfies $(d_\eps-\delta_\eps)^2 = -\Delta_\eps,$ as this will be convenient when we come to study the $\eta$-invariant.

The analysis of these operators is complicated by the singular coefficients as $\eps \to 0,$ which is an effect of the non-integrability of $\sH.$
We will make use of these formulas for $\eps>0$ and understand the behavior of these operators as $\eps \to 0$ by restricting the domain following Mazzeo-Melrose \cite{Mazzeo-Melrose:Ad}.

\subsection{Asymptotically solving Laplace's equation} \label{sec:AsymLap}

Our goal for this section is to understand what we can say about sections $\wt u$ of ${}^{\sR}\Omega^*X$ such that $\Delta_\eps\wt u = \cO(\eps^\ell)$ as $\eps\to 0$ for some $\ell.$ 
Let 
\[ \sE_0^* = {}^{\sR}\Omega^*(X)\rvert_{\eps=0}   \]
and define 
\begin{equation*}\label{eq:EkDef}
	\sE_k^* = \{ u_0 \in \sE_0^*: \exists \wt u\in {}^{\sR}\Omega^*(X) \text{ s.t. } \wt u\rvert_{\eps=0}=u_0 \text{ and } \Delta_\eps(\wt u) = \cO(\eps^{k-2}) \}.
\end{equation*}
We will determine the restrictions placed on the Taylor expansion of elements $\wt u$ of $\sE_k.$

The reason we are interested in these spaces is that they give an analytic realization of the spectral sequence induced by the splitting $TM = \sH \oplus \sV$ \cite{Mazzeo-Melrose:Ad, Forman}.
These satisfy
\begin{equation*}
	\sE^p_0 \supseteq \sE^p_1 \supseteq \ldots
\end{equation*}
and \cite[Theorem 1.3]{Forman} are eventually isomorphic to the singular cohomology of $M$ in that for each $p$ there is an $N \in \bbN$ such that
\begin{equation*}
	\sE^p_N = \sE^p_{\infty} = \bigcap_{j \in \bbN} \sE^p_j \cong \tH^p(M).
\end{equation*}
Indeed, we will see  that $\sE^p_5=\sE^p_{\infty}$ if $p \in \{n,n+1\}$ and that $\sE^p_3 = \sE^p_{\infty}$ if $p \notin\{n,n+1\}.$\\

Let us write 
\[ d_\eps - \delta_\eps = \eps^{-1}a_{-1} + a_0 + \eps a_1  \]
with the $a_i$ independent of $\eps$, and correspondingly
\[ -\Delta_\eps = (d_\eps-\delta_\eps)^2 = \eps^{-2}A_{-2} + \eps^{-1}A_{-1} + A_0 + \eps A_1 + \eps^2 A_2  \]
with the relations
\begin{equation*}
\begin{gathered}
	A_{-2}=a_{-1}^2, \quad 
	A_{-1}=a_{-1}a_0+a_0a_{-1}, \quad 
	A_0 = a_0^2+a_{-1}a_1+a_1a_{-1}, \\
	A_{1}=a_{1}a_0+a_0a_{1}, \quad 
	A_2 = a_1^2
\end{gathered}
\end{equation*}
Note that since $d_\eps-\delta_\eps$ is skew-adjoint, so are the individual $a_i$. \\

As explained in \S\ref{sec:Contact Complex}, the null space of
\begin{equation*}
	a_{-1} = \pmat{ 0 & L \\ -L^* & 0 }
\end{equation*}
is naturally identified with the Rumin complex. The differential of the Rumin complex, outside of middle degree, is correspondingly
\[ d_\sH = \Pi_{\Ker a_{-1}} d_H.   \]
It is also worth pointing out (see \cite[\S4]{Rumin:Sub}) that $a_{-1}$ is a section of the endomorphism bundle of $\Lambda^*({}^{\sR}T^*X)$ whose spectrum is constant. 
In particular it has a generalized inverse
\begin{equation*}
	a_{-1}^{\dagger} \in \CI(X;\End(\Lambda^*({}^{\sR}T^*X))) \quad \Mst \quad
	a_{-1}^\dagger a_{-1} = a_{-1}a_{-1}^\dagger = \text{Id} - \Pi_{\Ker a_{-1}}.
\end{equation*}
$ $

Given $u_0\in {}^{\sR}\Omega^*(X)\rvert_{\eps=0}$, we can smoothly extend it to a form in ${^{\sR}}\Omega^*(X).$
Any such extension $\wt u$ satisfies 
\[ -\Delta_\eps \wt u = a_{-1}^2u_0 \eps^{-2} + \cO(\eps^{-1})   \]
and so
\[ \sE_1^* = \sE_0^*\cap \Ker a_{-1}^2 =  \sE_0^*\cap \Ker a_{-1},     \]
by skew-adjointness. If $u_0\in \sE_1^*$, then any extension $\wt u = u_0 + u_1\eps + \cO(\eps^2)$ satisfies 
\[ -\Delta_\eps (\wt u) = (A_{-1}u_0 + A_{-2}u_1)\eps^{-1} + \cO(1) =  (a_{-1}a_0 u_0 + a_{-1}^2u_1)\eps^{-1} +\cO(1). \]
We can take
\[ u_1 = -a_{-1}^\dagger a_0 u_0,   \]
and with this choice $\Delta_\eps(\wt u)\in \cO(1)$, so 
\[ \sE_2^* = \sE_1^* = \Omega_{\sH}^*M.    \]
$ $

To determine when a form $u_0 \in \sE_2^*$ will be in $\sE_3^*,$ we start by specifying our preferred extension
of a form in $\sE_2^*$ to a form on ${}^{\sR}\Omega^*M.$
Note that if 
$\wt u= u_0+\eps u_1 + \eps^2 u_2 + \cO(\eps^3)$ with 
\begin{equation*}
	u_1 = -a_{-1}^\dagger a_0u_0,
\end{equation*}
then 
\begin{multline*}
	-\Delta_\eps(\wt u) = (A_{-2}u_2 + A_{-1} u_1 + A_0u_0) + \cO(\eps) \\
	= a_{-1}^2 u_2 + (a_{-1}a_0+a_0a_{-1})u_1 + (a_0^2 + a_{-1}a_1 + a_1a_{-1})u_0 + \cO(\eps) \\
	= a_{-1}^2 u_2 + a_{-1}a_0(-a_{-1}^\dagger a_0 u_0) + a_0a_{-1}(- a_{-1}^\dagger a_0u_0) + (a_0^2 +  a_{-1}a_1)u_0 + \cO(\eps) \\
	= a_{-1}^2u_2 - a_{-1}a_0a_{-1}^\dagger a_0 u_0 +  a_{-1}a_1 u_0 + a_0 \Pi_{\Ker  a_{-1}} a_0u_0 + \cO(\eps),
\end{multline*}
so we have a natural choice of $u_2$ by arranging for it to cancel out the part of the right hand side in the image of $a_{-1}$,
\[ u_2 = a_{-1}^\dagger (a_0a_{-1}^\dagger a_0 - a_1 - a_0\Pi_{\Ker a_{-1}} a_0)u_0 . \]
We define an `extension operator' 
\begin{equation*}
\begin{gathered}
	\Phi_2 : \sE_2^* \to {}^{\sR}\Omega^*(X), \\
	\Phi_2(u_0) = u_0 - \eps  a_{-1}^\dagger a_0u_0 + \eps^2  a_{-1}^\dagger (a_0 a_{-1}^\dagger a_0u_0 - a_1u_0 -  a_{-1}^\dagger a_0\Pi_{\Ker  a_{-1}} a_0 u_0) 
\end{gathered}
\end{equation*}
with the property that 
\[ -\Delta_{\eps}(\Phi_2(u_0)) = (\Pi_{\Ker  a_{-1}} a_0)^2 u_0 + \cO(\eps).   \]
This shows that a sufficient condition for $u_0 \in \sE_3^*$ is for $(\Pi_{\Ker  a_{-1}} a_0)^2 u_0=0.$

To see that this is a necessary condition, let $\wt u$ be any extension of $u_0$ such that $\Delta_{\eps}\wt u = \cO(1)$ and write 
\begin{equation*}
	\wt u = \Phi_2(u_0) + \eps \wt v.
\end{equation*}
Here $\wt v \in {}^{\sR}\Omega^*M$ satisfies $\Delta_{\eps}\wt v = \cO(\eps^{-1}),$ witnessing that $v_0 = \wt v|_{\eps=0}$ is an element of $\cE_1^*.$ 
Since $\sE_1^*=\sE_2^*$ we know that $\Delta_{\eps}\Phi_2(v_0) = \cO(1)$ and $\wt v = \Phi_2(v_0) + \eps \wt w,$ with $\wt w \in {}^{\sR}\Omega^*M$ unconstrained.
In this way we have shown that
\begin{equation*}
	\wt u = \Phi_2(u_0) + \eps \Phi_2(v_0) + \eps^2 w_0 + \cO(\eps^3), \Mforsome v_0 \in \sE_2^*, w_0 \in \sE_0^*,
\end{equation*}
and hence $(-\Delta_{\eps} \wt u)|_{\eps=0} =  (\Pi_{\Ker  a_{-1}} a_0)^2 u_0 + a_{-1}^2w_0.$
Since the two terms on the right hand side are orthogonal to each other, we see that the vanishing of $(\Pi_{\Ker  a_{-1}} a_0)^2 u_0$ is also a necessary condition for $u_0 \in \sE_3^*.$

Next we will show that $\sE_3^*=\sE_4^*,$ i.e., that any form $u_0 \in \sE_3^*$ has an extension $\wt u$ such that $\Delta_{\eps}\wt u = \cO(\eps^2).$
If $u_0 \in \sE_3^*,$ we have shown that
\begin{equation}\label{eq:E3}
	a_{-1}u_0=0, \quad \Pi_{\Ker  a_{-1}}a_0 u_0=0,
\end{equation}
thus $\Phi_2(u_0)$ simplifies slightly to
\begin{equation*}
	\Phi_2(u_0) = u_0 + \eps u_1 + \eps^2u_2 = u_0 - \eps  a_{-1}^\dagger a_0u_0 + \eps^2  a_{-1}^\dagger (a_0 a_{-1}^\dagger a_0u_0 - a_1u_0)
\end{equation*}
and for any form $u_3 \in \sE_0^*$ we have that
\begin{multline*}
	-\eps^{-1}\Delta_{\eps}(\Phi_2(u_0)+\eps^3u_3) \\
	= a_{-1}^2(u_3) + (a_{-1}a_0+a_0a_{-1})(u_2)
	+(a_0^2 + a_1a_{-1}+a_{-1}a_1)(u_1)
	+ (a_0a_1+a_1a_0)(u_0) + \cO(\eps) \\
	=a_{-1}^2(u_3)
	+ a_{-1}a_0(u_2)
	+ a_{-1}a_1(u_1) 
	+ a_0 \Pi_{\Im a_{-1}} \Big(a_0a_{-1}^{\dagger}a_0 u_0 - a_1 u_0 
	 \Big) \\
	-a_0^2a_{-1}^{\dagger}a_0 u_0 
	-a_1\Pi_{\Im a_{-1}} a_0 u_0 
	+ a_0a_1u_0 + a_1a_0u_0 + \cO(\eps)\\
	=a_{-1}^2(u_3)
	+ a_{-1}a_0(u_2)
	+ a_{-1}a_1(u_1) 
	+a_0\Pi_{\Ker a_{-1}}(a_1- a_0a_{-1}^{\dagger}a_0)u_0
	+\cO(\eps).
\end{multline*}
We can choose $u_3$ so that it cancels out the part of the right hand side in the image of $a_{-1},$
\begin{equation}\label{eq:u3}
	u_3 = a_{-1}^{\dagger}
	\lrpar{-a_0u_2-a_1u_1
	-a_{-1}^{\dagger}a_0\Pi_{\Ker a_{-1}}(a_1- a_0a_{-1}^{\dagger}a_0)u_0},
\end{equation}
leaving us with
\begin{equation*}
\begin{gathered}
	(\eps^{-1}\Delta(\wt u))|_{\eps=0}
	= \Pi_{\Ker a_{-1}}a_0\Pi_{\Ker a_{-1}}
	\lrpar{ a_1 - a_0a_{-1}^{\dagger}a_0 }u_0.
\end{gathered}
\end{equation*}
It turns out that this vanishes for any $u_0$ satisfying \eqref{eq:E3}.
Indeed, we see from \eqref{eq:RumCmpxDesc} that $a_1 - a_0 a_{-1}^\dagger a_0$ vanishes on forms in $\Ker a_{-1}$ of degree $p\notin\{n,n+1\}$ because it is an off-diagonal operator with respect to \eqref{eq:SplittingM}, while on forms of degree $n$ or $n+1$ this is given by 
\[\Pi_{\Ker  a_{-1}} a_0 \Pi_{\Ker  a_{-1}} (a_1 - a_0 a_{-1}^\dagger a_0 )u_0 = (d_\sH-\delta_\sH)(D_{\sH}-D_{\sH}^*)u_0\] 
and both $d_{\sH}$ and $\delta_{\sH}$ vanish on the images of $D_{\sH},$ $D_{\sH}^*.$
Thus we have found that
\begin{equation*}
	\sE_3^p=\sE_4^p = 
	\begin{cases}
	\Ker_{\Omega_{\sH}^pM}(d_{\sH}\delta_{\sH} + \delta_{\sH}d_{\sH}) & \Mif p \notin\{n,n+1\}\\
	\Ker_{\Omega_{\sH}^nM}(d_{\sH}\delta_{\sH}) &\Mif p = n\\
	\Ker_{\Omega_{\sH}^{n+1}M}(\delta_{\sH}d_{\sH}) & \Mif p=n+1
	\end{cases}
\end{equation*}
For forms of degree $p \notin\{n,n+1\},$ $\sE_3^p$ is isomorphic to the cohomology of the Rumin complex, hence to the de Rham cohomology of $M,$ and so
$\sE_3^p = \sE_{\infty}^p.$\\

Finally, let us analyze when a form in $\sE^p_4$ will be in $\sE^p_5.$
As we did for $\sE_2^*,$ we will start by constructing a preferred extension operator $\Phi_4:\sE_4^*\lra {}^{\sR}\Omega^*M.$
Given $u_0 \in \sE_4^*,$ define $u_1,$ $u_2$ by $\Phi_2(u_0) = u_0 + \eps u_1 + \eps^2 u_2,$ define $u_3$ by \eqref{eq:u3} and note that for any $u_4 \in \sE_0^*,$
\begin{equation*}
	-\eps^{-2}\Delta_{\eps}(u_0 + \eps u_1 + \eps^2 u_2+ \eps^3u_3 + \eps^4 u_4)
	= a_{-1}^2 u_4 + A_{-1}u_3 + A_0 u_2 + A_1 u_{1} + A_2 u_0 + \cO(\eps),
\end{equation*}
so a natural choice for $u_4$ is for it to cancel out the part of the right hand side that is in the image of $a_{-1}.$
So with these choices of $u_1,$ $u_2,$ $u_3,$ we define
\begin{equation*}
\begin{gathered}
	\Phi_4:\sE_4^*\lra {}^{\sR}\Omega^*M, \\
	\Phi_4(u_0) = 
	u_0 + \eps u_1 + \eps^2 u_2+ \eps^3u_3
	- \eps^4 (a_{-1}^{\dagger})^2\Pi_{\Im a_{-1}}(A_{-1}u_3 + A_0 u_2 + A_1 u_{1} + A_2 u_0).
\end{gathered}
\end{equation*}
For any $u_0 \in \sE_4^*$ we have
\begin{multline*}
	-\eps^{-2}\Delta_{\eps}(\Phi_4(u_0))
	= \Pi_{\Ker a_{-1}}(A_{-1}u_3 + A_0 u_2 + A_1 u_{1} + A_2 u_0) + \cO(\eps)\\
	= \Pi_{\Ker a_{-1}} \Big(
	a_0\Pi_{\Im a_{-1}}(-a_0u_2-a_1u_1-a_{-1}^{\dagger}a_0\Pi_{\Ker a_{-1}}(a_1-a_0a_{-1}^{\dagger}a_0)u_0) \\
	+ (a_0^2 + a_1a_{-1})u_2 + (a_0a_1+a_1a_0)u_1 + a_1^2u_0	
	\Big)+\cO(\eps) \\
	= \Pi_{\Ker a_{-1}}a_0\Pi_{\Ker a_{-1}}(a_0 u_2 + a_1u_1) 
	+ (\Pi_{\Ker a_{-1}}(a_1-a_0a_{-1}^{\dagger}a_0))^2u_0 +\cO(\eps)
\end{multline*}
and \eqref{eq:E3} implies that $\Pi_{\Ker a_{-1}}a_0\Pi_{\Ker a_{-1}}(a_0 u_2 + a_1u_1) =0$. Indeed, since $a_1$ preserves $\Im(a_{-1})$, and $u_2\in \Im a_{-1},$ we have 
\[ \Pi_{\Ker a_{-1}} a_0a_1(a_{-1}^\dagger a_0u_0)\equiv 0, \quad \Pi_{\Ker a_{-1}}a_0 u_2=0,  \]
as claimed. So
\begin{equation*}
	-\eps^{-2}\Delta_{\eps}(\Phi_4(u_0)) =  (\Pi_{\Ker a_{-1}}(a_1-a_0a_{-1}^{\dagger}a_0))^2u_0 +\cO(\eps).
\end{equation*}
Thus a sufficient condition for $u_0 \in \sE_5^*$ is  $(\Pi_{\Ker a_{-1}}(a_1-a_0a_{-1}^{\dagger}a_0))^2u_0=0.$

To see that this condition is necessary, we proceed as we did for $\sE_2^*.$
That is, let $u_0 \in \sE_4^*$ and suppose that $\wt u$ is any extension such that $\Delta_{\eps}(\wt u) = \cO(\eps^2).$
We can write $\wt u = \Phi_4(u_0) + \eps \wt v$ with $\Delta_{\eps} (\wt v) = \cO(\eps).$
Thus $v_0 = \wt v|_{\eps=0}$ is in $\sE_3^*$ which we have shown equals $\sE_4^*$ and so $\wt v = \Phi_4(v) + \eps \wt w,$ with $\Delta_{\eps}(\wt w) = \cO(1).$
Thus $w_0 = \wt w|_{\eps=0}$ in $\sE_2^*$ and we can apply our previous analysis of $\sE_2^*$ to conclude that
\begin{multline*}
	\wt u = \Phi_4(u_0) + \eps\Phi_4(v_0) + \eps^2\Phi_2(w_0) + \eps^3\Phi_2(x_0) + \eps^4 y_0 + \cO(\eps^5),  \\
	\Mforsome v_0 \in \sE_4^*, w_0,x_0 \in \sE_2^*, y_0 \in \sE_0^*
\end{multline*}
and hence
\begin{equation*}
	-\eps^{-2}\Delta_{\eps}(\wt u)|_{\eps=0} = (\Pi_{\Ker a_{-1}}(a_1-a_0a_{-1}^{\dagger}a_0))^2u_0
	+ (\Pi_{\Ker a_{-1}} a_0)^2 w_0 + a_{-1}^2 y_0.
\end{equation*}
These three terms are orthogonal as they correspond to separate parts of the Hodge decomposition of Rumin's complex and so the vanishing of 
$(\Pi_{\Ker a_{-1}}(a_1-a_0a_{-1}^{\dagger}a_0))^2u_0$ is also a necessary condition for $u_0 \in \sE_5^*.$
This shows that
\begin{equation*}
	\sE_5^p =
	\begin{cases}
	\Ker_{\Omega_{\sH}^pM}(d_{\sH}\delta_{\sH} + \delta_{\sH}d_{\sH}) & \Mif p \notin\{n,n+1\}\\
	\Ker_{\Omega_{\sH}^nM}(d_{\sH}\delta_{\sH}) \cap \Ker_{\Omega_{\sH}^nM}(D_{\sH}^*D_{\sH}) &\Mif p = n\\
	\Ker_{\Omega_{\sH}^{n+1}M}(\delta_{\sH}d_{\sH}) \cap \Ker_{\Omega_{\sH}^{n+1}M}(D_{\sH}D_{\sH}^*) & \Mif p=n+1
	\end{cases}
\end{equation*}
In all degrees $\sE_5^p$ is isomorphic to the cohomology of the Rumin complex, hence to the de Rham cohomology of $M,$ and so
$\sE_5^p = \sE_{\infty}^p.$\\

In summary, in terms of the associated graded spaces
\begin{equation*}
	\sG_k^p = \sE_k^p \setminus \sE_{k+1}^p,
\end{equation*}
we have established the following decompositions of the space of sub-Riemannian forms,
\begin{equation*}
\begin{aligned}
	{}^{\sR}\Omega^p X|_{\eps= 0}
	&= \underbrace{\sG_0^p}_{(\Omega_{\sH}^pM)^\perp} 
	\oplus \underbrace{\sG_2^p}_{ \Im(\Delta_\sH)} 
	\oplus \underbrace{\sG_\infty^p}_{\Ker \Delta_\sH=\sE_{\infty}^p}, \quad \Mif p \notin\{n,n+1\}\\
	{}^{\sR}\Omega^n X|_{\eps =0} 
	&= \underbrace{\sG_0^n}_{(\Omega_{\sH}^nM)^\perp} 
	\oplus \underbrace{\sG_2^n}_{\Im(d_\sH\delta_\sH)^2} 
	\oplus \underbrace{\sG_4^n}_{\Im(D_{\sH}^*D_{\sH})} 
	\oplus \underbrace{\sG_{\infty}^n}_{\Ker \Delta_{\sH}}, \\
	{}^{\sR}\Omega^{n+1} X|_{\eps =0} 
	&= \underbrace{\sG_0^{n+1}}_{(\Omega_{\sH}^{n+1}M)^\perp} 
	\oplus \underbrace{\sG_2^{n+1}}_{\Im(\delta_\sH d_\sH)^2} 
	\oplus \underbrace{\sG_4^{n+1}}_{\Im(D_{\sH}D_{\sH}^*)} 
	\oplus \underbrace{\sG_{\infty}^{n+1}}_{\Ker \Delta_{\sH}=\sE_{\infty}^{n+1}}.
\end{aligned}
\end{equation*}
It is worth noting that the decompositions above are homogeneous with respect to the natural contact weight and associated filtration, which is defined by assigning weight 1 to forms in $\sH^*$ and weight 2 to forms in $\sV^*$.

We will use the extension operators from above to define two `effective normal operators' of the Hodge Laplacian.
We define the $\sE_2^*$-effective normal operator to be 
\begin{equation*}
\begin{gathered}
	\tN_{\text{eff},2}(\Delta_{\eps}): \sE_2^* \lra \sE_2^* \\
	\tN_{\text{eff},2}(\Delta_{\eps})(u_0) = \Delta_{\eps}(\Phi_2(u_0))\rvert_{\eps=0} = -(\Pi_{\Ker  a_{-1}} a_0)^2 u_0 = (d_\sH\delta_{\sH}+\delta_{\sH}d_{\sH})u_0,
\end{gathered}
\end{equation*}
and we define the $\sE_4^*$-effective normal operator to be
\begin{equation*}
\begin{gathered}
	\tN_{\text{eff},4}(\eps^{-2}\Delta_{\eps}): \sE_4^p \lra \sE_4^p \\
	\tN_{\text{eff},4}(\eps^{-2}\Delta_{\eps})(u_0) = \eps^{-2}\Delta_{\eps}(\Phi_4(u_0))\rvert_{\eps=0} = -(\Pi_{\Ker a_{-1}}(a_1-a_0a_{-1}^{\dagger}a_0))^2u_0 \\
	= \begin{cases}
	D_{\sH}^*D_{\sH}u_0 & \Mif p = n \\
	D_{\sH}D_{\sH}^*u_0 & \Mif p = n+1
	\end{cases}
\end{gathered}
\end{equation*}

These operators will be useful in the construction of the heat kernel.
In similar situations, but for Laplacians without coefficients that are singular in $\eps,$
(e.g., \cite{Albin-Aldana-Rochon:Rel, Albin-Rochon-Sher:FibCusps, Albin-Rochon-Sher:Cusps, Albin-Rochon-Sher:Wedge})
one would define a normal operator for the Laplacian at $\eps=0$
by taking $u_0 \mapsto \Delta_{\eps}(\wt u)|_{\eps=0}$ where $\wt u$ is any extension of $u_0.$
We will be doing the same in our setting save that the choice of extension of $u_0$ to $\wt u$ is 
constrained by the singularity of $\Delta_{\eps}$ as $\eps \to 0,$ and so we have constructed 
explicit extension operators.

\subsection{Asymptotically solving the heat equation} \label{sec:AsymHeat}

In the previous section we showed that we can write
\begin{equation*}
	\Id_{{}^{\sR}\Omega^*X|_{\eps=0}} = \Pi_{\sG_0} \oplus \Pi_{\sG_2} \oplus \Pi_{\sG_4} \oplus \Pi_{\sG_{\infty}},
\end{equation*}
where $\Pi_{\sG_j}$ denotes the orthogonal projection onto $\sG_j.$ (Similarly, we will use $\Pi_{\sE_j}$ to denote the orthogonal projection onto $\sE_j.$)
Heuristically, by ignoring everything except the effective normal operators, we could expect the heat kernel of $\Delta_{\eps}$ to take the form
\begin{equation*}
	e^{-t\eps^{-2}(a_{-1})^2} \Pi_{\sG_0} 
	\oplus 
	e^{-t(d_{\sH}+\delta_{\sH})^2}\Pi_{\sG_2}
	\oplus
	e^{-t\eps^2(D_{\sH} + D_{\sH}^*)^2}\Pi_{\sG_4}
	\oplus 
	\Pi_{\sG_{\infty}}.
\end{equation*}
While the truth is a bit more complicated this exhibits the different regimes necessary to understand the asymptotics of the heat kernel, namely when $t$ and $\eps^2$ go to zero together, when $\eps$ goes to zero with $t$ bounded, and when $t^{-1}$ and $\eps^2$ go to zero together.  For later use, we will explain how to formally solve the heat equation to arbitrarily high order in $\eps$ in each of these regimes.

First suppose we wish to solve $(\pa_t + \Delta_{\eps})W^{(0)}=0$ with initial data $W^{(0)}|_{t=0}=f(\eps)$ where $f(0) \in \sG_0.$
Equivalently, with $s= t/\eps^2,$ we wish to solve $(\pa_s + \eps^2\Delta_{\eps})W^{(0)}=0.$
As $\eps^2\Delta_{\eps}|_{\eps=0} = a_{-1}^2$ we set $W^{(0)}_0(s,\eps)= \exp(-s(a_{-1})^2)f$ and then we have
\begin{equation*}
	\begin{cases}
	(\pa_s + \eps^2\Delta_{\eps})W^{(0)}_0=\cO(\eps), \\
	W^{(0)}_0(0,\eps) = f(\eps).
	\end{cases}
\end{equation*}
To show by induction that we can solve the heat equation to arbitrary order in $\eps,$ suppose that we have found $W^{(0)}_0(s,\eps), \ldots, W^{(0)}_k(s,\eps)$ such that 
\begin{equation*}
	\begin{cases}
	(\pa_s + \eps^2\Delta_{\eps})\Big( \sum_{j=0}^k \eps^j W^{(0)}_j\Big)= \eps^{k+1}R^{(0)}_k(s,\eps), \\
	\sum_{j=0}^k \eps^j W^{(0)}_j(0,\eps) = f(\eps)
	\end{cases}
\end{equation*}
and note that $W^{(0)}_{k+1}(s) = \exp(-sa_{-1}^2) \star R^{(0)}_{k}(s,0),$ where $\star$ denotes convolution, completes the inductive step.

Next suppose we wish to solve $(\pa_t + \Delta_{\eps})W^{(2)}=0$ with initial data $W^{(2)}|_{t=0}=f(\eps)$ where $f(0) \in \sG_2.$
As $(\Delta_{\eps}\circ \Phi_2)|_{\eps=0} = (d_{\sH}+\delta_{\sH})^2 \Pi_{\sG_2}$ we set
\begin{equation*}
	W^{(2)}_0(t,\eps) = \Phi_2( e^{-t(d_{\sH}+\delta_{\sH})^2} f(0))
\end{equation*}
and then we have
\begin{equation*}
	\begin{cases}
	(\pa_t + \Delta_{\eps})W^{(2)}_0=\cO(\eps), \\
	W^{(2)}_0(0,\eps) = \Phi_2(f(0)) = f(\eps) + \cO(\eps).
	\end{cases}
\end{equation*}
If, inductively, we have found
$W^{(2)}_0(t,\eps), \ldots, W^{(2)}_k(t,\eps)$ such that 
\begin{equation*}
	\begin{cases}
	(\pa_t + \Delta_{\eps})\Big( \sum_{j=0}^k \eps^j W^{(2)}_j\Big)= \eps^{k+1}R^{(2)}_k(t,\eps), \\
	\sum_{j=0}^k \eps^j W^{(2)}_j(0,\eps) =  f(\eps) + \eps^{k+1}w^{(2)}_{k+1}(\eps)
	\end{cases}
\end{equation*}
then we can improve the error in the first line by adding 
\begin{equation*}
	-\eps^{k+1}\Big[
	\eps^2\Big((a_{-1}^{\dagger})^2 \Pi_{\sG_0}(R^{(2)}_{k}(t,0))\Big)\\
	+ \Phi_2\Big( e^{-t(d_{\sH}+\delta_{\sH})^2} \Pi_{\sE_2}(R^{(2)}_{k}(t,0)) 
	\Big)
	\Big]
\end{equation*}
and for the error in the initial condition we can use the previous discussion to solve the heat equation to arbitrarily high order with initial data $\Pi_{\sG_0}w^{(2)}_{k+1}(0)$ and then use the inductive hypothesis to solve the heat  equation to order $k+1$ with initial data $\Pi_{\sE_2}w^{(2)}_{k+1}(0).$

Finally suppose we wish to solve $(\pa_t + \Delta_{\eps})W^{(4)}=0$ with initial data $W^{(4)}|_{t=0}=f(\eps)$ where $f(0) \in \sG_4.$
Equivalently, with $s = \eps^{2}t$ we wish to solve $(\pa_s + \eps^{-2}\Delta_{\eps})W^{(4)}=0.$
As $\eps^{-2}(\Delta_{\eps}\circ \Phi_4)|_{\eps=0} = (D_{\sH}+D_{\sH}^*)^2 \Pi_{\sG_4}$ we set
\begin{equation*}
	W^{(4)}_0(s,\eps) = \Phi_4( e^{-s(D_{\sH}+D_{\sH}^*)^2} f(0))
\end{equation*}
and then we have
\begin{equation*}
	\begin{cases}
	(\pa_s + \eps^{-2}\Delta_{\eps})W^{(4)}_0=\cO(\eps), \\
	W^{(4)}_0(0,\eps) = \Phi_4(f(0)) = f(\eps) + \cO(\eps).
	\end{cases}
\end{equation*}
If, inductively, we have found
$W^{(4)}_0(s,\eps), \ldots, W^{(4)}_k(s,\eps)$ such that 
\begin{equation*}
	\begin{cases}
	(\pa_s + \eps^{-2}\Delta_{\eps})\Big( \sum_{j=0}^k \eps^j W^{(4)}_j\Big)= \eps^{k+1}R^{(4)}_k(s,\eps), \\
	\sum_{j=0}^k \eps^j W^{(4)}_j(0,\eps) =  f(\eps) + \eps^{k+1}w^{(4)}_{k+1}(\eps)
	\end{cases}
\end{equation*}
then we can improve the error in the first line by adding 
\begin{multline*}
	-\eps^{k+1}\Big[
	\eps^4\Big((a_{-1}^{\dagger})^2 \Pi_{\sG_0}(R^{(4)}_{k}(s,0))\Big)
	+ \eps^{2}\Phi_2\Big( (d_{\sH}+\delta_{\sH})^{-2} \Pi_{\sG_2}(R^{(4)}_{k}(s,0)) \Big) \\
	+ \Phi_4\Big( e^{-s(D_{\sH}+D_{\sH}^*)^2} \Pi_{\sG_4}(R^{(4)}_{k}(s,0)) \Big)
	\Big],
\end{multline*}
where we use that $(d_{\sH}+\delta_{\sH})^2$ has an inverse on $\sG_2,$
 then for the error in the initial condition we can use the previous discussions to solve the heat equation to arbitrarily high order with initial data $(\Id-\Pi_{\sG_4})w^{(4)}_{k+1}(0),$ and then we can use the inductive hypothesis to solve the heat equation to order $k+1$ with initial data $\Pi_{\sE_4}w^{(4)}_{k+1}(0).$

\section{Interlude: manifolds with corners} \label{sec:Corners}

We briefly review some basic concepts on manifolds with corners and refer the reader to, e.g., \cite[\S2A]{Mazzeo:Edge}, \cite{Melrose:Conormal, Melrose:Corners, Grieser, Dai-Melrose, Melrose:APS, Epstein-Melrose-Mendoza}, for more details.

Recall that a map $[0,\infty)^m \lra [0,\infty)^m$ is smooth if it has a smooth extension to a map between open neighborhoods of $[0,\infty)^m$ in $\bbR^m.$
A smooth $m$-dimensional manifold with corners $W$ is a manifold smoothly modeled on $[0,\infty)^m$ with embedded boundary hypersurfaces. This latter condition is equivalent to the existence, for each boundary hypersurface $H,$ of a smooth function $\rho:W \lra \bbR$ such that
\begin{equation*}
	\rho(W) \subseteq [0,\infty), \quad
	\rho^{-1}(0) = H, \quad
	d\rho \text{ has no zeroes on }H;
\end{equation*}
any such function is known as a boundary defining function for $H.$
A product of boundary defining functions, one per each boundary hypersurface of $W,$ is known as a `total' boundary defining function for $W.$ 

A construction we use repeatedly to obtain new manifolds with corners is (real) blow-up of a `$p$-submanifold' (or `submanifold of product-type').
An embedded submanifold $Y \subseteq W$ is a $p$-submanifold if every point $q \in Y$ has a neighborhood $\cU$ such that 
\begin{equation*}
	W \cap \cU = W' \times W''
\end{equation*}
where $W''$ has no boundary and $Y \cap \cU = W' \times \{q''\}$ for some $q'' \in W''.$ (Thus the interval $\{x=\tfrac12, 0\leq y\leq 1\}$ is a $p$-submanifold of the unit square, while the diagonal $\{0\leq x=y \leq 1\}$ is {\em not} a $p$-submanifold of the unit square.) These are the submanifolds that have `nice' tubular neighborhoods. The blow-up of $W$ along $Y,$ denoted $[W;Y],$ is the manifold with corners obtained by removing $Y$ and replacing it with its inward pointing spherical normal bundle in $W,$
\begin{equation*}
	[W;Y] = W\setminus Y \bigsqcup (N_W^+Y\setminus \{0\})/\sim, \quad
	\Mwhere v \sim \lambda v \Mforall v \in N_W^+Y\setminus \{0\}, \lambda \in \bbR^+.
\end{equation*}
A modification given a subbundle $S$ of the conormal bundle to $Y$ in $W,$ $N^*_WY,$ known as the `parabolic blow-up of $W$ along $Y$ with parabolic directions $S$' and denoted $[W;Y,S],$ is obtained by replacing the radial dilations with anisotropic dilations,
\begin{multline*}
	[W;Y,S] = W\setminus Y \bigsqcup (N_W^+Y\setminus \{0\})/\sim_S, \\
	\Mwhere (v_{S^{\circ}}, v_{S'}) \sim (\lambda v_{S^{\circ}}, \lambda^2 v_{S'}) \Mforall v=(v_{S^{\circ}}, v_{S'}) \in N_W^+Y\setminus \{0\}, \lambda \in \bbR^+,
\end{multline*}
and where we have chosen a complementary sub-bundle $S'$ to the annihilator $S^{\circ}$ of $S.$

A blow-up comes with a blow-down map
\begin{equation*}
	[W;Y] \lra W, \quad [W;Y,S] \lra W,
\end{equation*}
which we usually denote $\beta.$ If $L \subseteq W$ is a submanifold which is equal to the closure of $L\setminus Y,$ then the `interior lift' of $L$ along $\beta$ is defined to be the closure of $\beta^{-1}(L\setminus Y).$

Every manifold with corners can be embedded into a closed manifold (see, e.g., \cite[Theorem 4.2]{Albin-Melrose:Gps}) and a smooth function on a manifold with corners is, by definition, the restriction of a smooth function on a closed manifold. However, it is convenient and often necessary to work within the larger class of {\bf $\cI$-smooth functions}, or functions that are smooth with respect to index sets (also known as polyhomogeneous functions). On a manifold with boundary $W,$ with boundary defining function $\rho,$ an $\cI$-smooth function $f$ with index set $\cE\subseteq \bbC \times \bbN_0$ is a function that is smooth in the interior of $W$ and has an asymptotic expansion
\begin{equation*}
	f \sim \sum_{(s, p) \in \cE} a_{s,p}(y) \rho^{s} (\log \rho)^p \quad \Mas \rho \to 0,
\end{equation*}
where the coefficients, $a_{s,p}(y),$ are smooth functions on $\pa W.$
We denote the set of such functions by
\begin{equation*}
	\sA_{phg}^{\cE}(W).
\end{equation*}
In order for this to make sense, and behave well with respect to change of boundary defining function, we require of $\cE$ that:
\begin{quote}
	i) Any infinite sequence $((s_j,p_j)) \subseteq \cE$ satisfies $\Re s_j \to \infty,$\\
	ii) If $(s,p) \in \cE$ then $(s+k,p') \in \cE$ for all $k \in \bbN_0$ and $0\leq p'\leq p.$
\end{quote}
On a manifold with corners $\cI$-smooth functions have index sets at each boundary hypersurface and joint expansions at corners, see \cite[\S2A]{Mazzeo:Edge} for details.
As $\cI$-smooth functions are $\CI(W)$-modules it is straightforward to define $\cI$-smooth sections of vector bundles.

We say that a function $f$ vanishes to infinite order at a boundary hypersurface $H$ of a manifold with corners $W$ if all of the coefficients in its Taylor expansion at $H$ are identically zero. We denote the set of smooth functions on $W$ that vanish to infinite order at all boundary hypersurfaces of $W$ by $\dCI(W),$ and those that vanish to infinite order at all boundary hypersurfaces except for $H$ by $\dCI_H(W).$

A beautiful and powerful geometric technique of Melrose for understanding the mapping properties of an operator and the composition of two operators is provided by the pull-back and push-forward theorems.
A map $f:W \lra W'$ between manifolds with corners is a $b$-map, i.e., a `boundary map', if the pull-back of any boundary defining function of a boundary hypersurface of $W'$ is the product of boundary defining functions of $W.$ An example is the projection of the square onto one of its sides and a non-example is the map from the unit square onto $[0,2]$ sending $(x,y)$ to $x+y.$
A $b$-map is `simple' if the pull-back of a boundary defining function of $W'$ is a boundary defining function of $W.$
The pull-back of an $\cI$-smooth function by a $b$-map is again an $\cI$-smooth function, \cite[Proposition A.13]{Mazzeo:Edge}.
For a simple $b$-map, $f,$ the index set of $f^*u$ at a boundary hypersurface $H$ of $W'$ is $\bbN_0$ if $f(H) = W'$ and is equal to the index set of $u$ at the boundary hypersurface $f(H)$ otherwise.

The push-forward of a density along a fiber bundle map is the fiberwise integral of that density; in general push-forward is defined as the dual of pull-back.
A $b$-fibration is a $b$-map between manifolds with corners that restricts to a fiber bundle over the interior of each boundary face and does not increase codimension (see \cite[Definition 3.9]{Grieser}).
Push-forward along a $b$-fibration is especially well-behaved for $b$-densities: a density on the interior of $W$ is a $b$-density if its product with a total boundary defining function is a non-degenerate density on $W.$
If $u$ is an $\cI$-smooth density on $W$ and $f:W \lra W'$ is a simple $b$-fibration, then $f_*u$ is an $\cI$-smooth $b$-density on $W'.$ The index set of $f_*u$ at a boundary hypersurface $H$ of $W'$ is the `extended union' of the index sets of $u$ at the boundary hypersurfaces of $W$ that are mapped onto $H$ by $f.$ Here the extended union of two index sets $\cE$ and $\cF$ is 
\begin{equation*}
	\cE \bar\cup \cF
	= \cE \cup \cF \cup \{ (z, p) \in \bbC \times \bbN_0: (z,p-1) \in \cE \cap \cF\}.
\end{equation*}

In this paper we will only require pull-back and push-forward along simple $b$-maps; we refer the reader to the references for the results for non-simple $b$-maps.

\section{The heat kernel construction on model spaces} \label{sec:ModelHeat}

In this section we will describe the construction of the heat kernel of a Laplace-type operator on a closed manifold from \cite[Chapter 7]{Melrose:APS} and the analogous construction of the heat kernel of suitable hypoelliptic operators on `Heisenberg manifolds' following \cite{Beals-Greiner-Stanton, Taylor:Noncom, Ponge:Spectral}.
This will allow us to review the methods involved in situations simpler than the heat kernel of a manifold undergoing the sub-Riemannian limit and will serve as part of that construction.\\

Let us start by considering the structure of the heat kernel for the scalar Laplacian on $Y = \bbR^m$ as a right density,
\begin{equation}\label{eq:HeatKerY}
	\cK_Y(t, \zeta, \zeta') 
	= \frac1{(4\pi t)^{m/2}} \exp\Big( -\frac{|\zeta-\zeta'|^2}{4t} \Big) \; d\zeta'.
\end{equation}
The function $\exp\Big( -\frac{|\zeta-\zeta'|^2}{4t} \Big)$ is smooth everywhere on $\bbR^+_t \times Y \times Y$ except at the submanifold
$\{ t=0, \zeta=\zeta'\}$ where it is not defined and does not have well-defined limits.
We resolve this singularity by blowing-up the diagonal at time zero; that is, we remove this submanifold and we replace it with (a modified version of) its inward pointing spherical normal bundle.
We denote the resulting space by $\wt {HY}$ and the map back to $\bbR^+_t \times Y \times Y,$ the `blow-down map', by
\begin{equation*}
	\beta: \wt{HY} \lra \bbR^+_t \times Y \times Y,
\end{equation*}

In essence this comes down to introducing polar coordinates and then `taking them seriously'. 
\begin{figure}[H]
    \def\svgwidth{\columnwidth}
    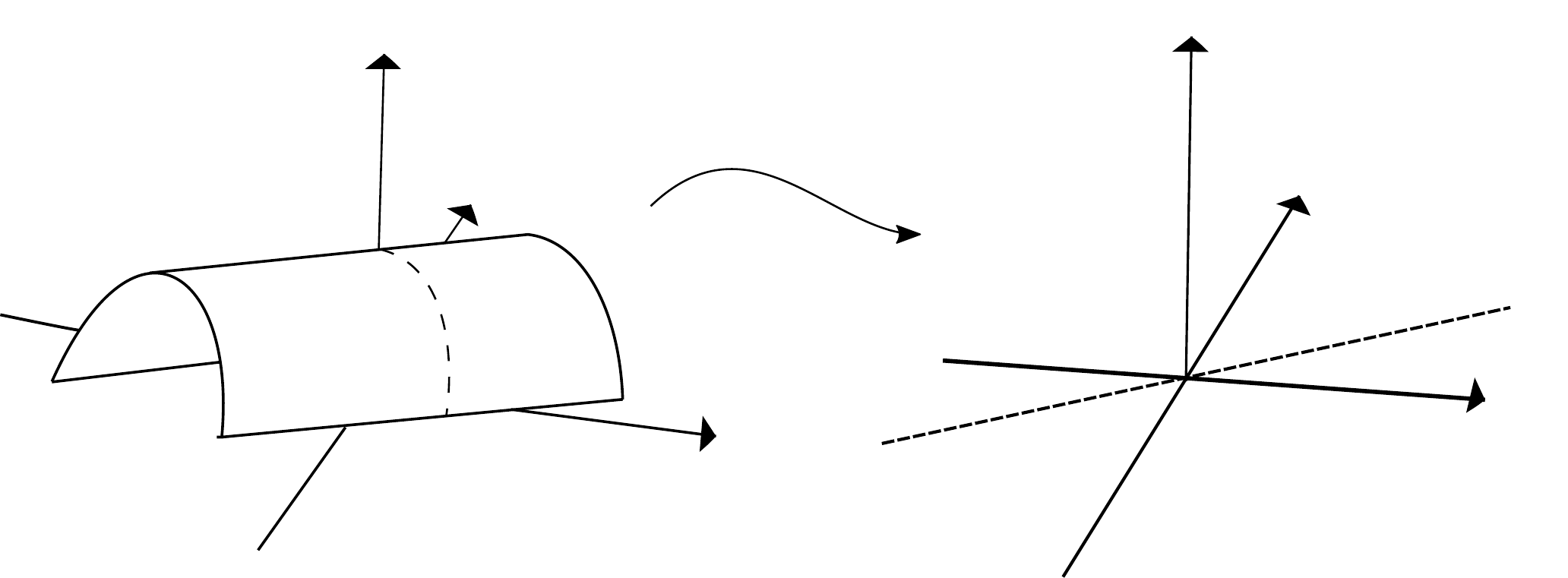

    \caption{The heat space $\wt{HY}$ and its blow-down map $\beta: \wt{HY} \to \bbR_t^+\times Y\times Y$}
    \label{fig:simple heat space}
\end{figure} 
Thus for coordinates on $\wt{HY}$ we can take
\begin{equation*}
	\wt{R}= \oldsqrt[4]{t^2+|\zeta-\zeta'|^4}, \quad 
	\wt\theta = (\wt\theta_t, \wt\theta_{Y}) = \Big(\frac{t}{\wt R^2}, \frac{\zeta-\zeta'}{\wt R}\Big), \quad 
	\zeta'
\end{equation*}
and the map $\beta$ is given by
\begin{equation*}
	\xymatrix @R=1pt{
	\wt{HY} = \bbR^+ \times \{ (\wt \theta_t, \wt \theta_Y): \wt\theta_t^2 + \wt\theta_Y^4=1\} \times \bbR^m \ar[rr]^-{\beta} & &
	\bbR^+ \times Y \times Y \\
	(\wt R, \wt \theta,  \zeta') \ar@{|->}[rr] & &
	(\wt \theta_t \wt R^2,  \zeta' + \wt R\wt \theta_Y, \zeta') }
\end{equation*}

The pull-back the heat kernel of $\Delta_Y$ to $\wt{HY}$ is
\begin{equation*}
	\beta^*\cK_Y(\wt R, \wt \theta, \zeta')
	= \frac1{(4\pi \wt R^2\wt \theta_t)^{m/2}} \exp\Big( -\frac{|\wt \theta_Y|^2}{4\wt \theta_t} \Big) \; d\zeta'
	\in \wt R^{-m} \CI(\wt {HY};\beta_R^*\density_Y),
\end{equation*}
where $\beta_R:\wt {HY} \lra Y$ is the composition of $\beta$ and the projection onto the right factor of $Y$ and $\density_Y$ denotes the density bundle on $Y.$
Its newfound smoothness is due to the pleasant fact that now the numerator and denominator in the exponential are never simultaneously zero.

This construction of the Euclidean heat space, where the $t$-direction is treated with a different weight as the spatial directions, is known as a `parabolic blow-up' and denoted
\begin{equation*}
	\wt{HY} = [\bbR^+_t\times Y^2; \{t=0\}\times \diag_Y, \ang{dt}].
\end{equation*}
An alternate construction that is just as good for understanding the structure of the heat kernel is to first introduce $\tau = \sqrt t$ as a smooth global variable\footnote{To declare that  $\tau = \sqrt t$ is smooth is to change the smooth structure of $\bbR^+ \times Y^2.$ It is the same as carrying out a parabolic blow-up of the boundary hypersurface $\{t=0\},$ i.e., replacing $\bbR^+ \times Y^2$ with $[\bbR^+\times Y^2;\{t=0\},\ang{dt}].$} and then blow-up the diagonal at time zero homogeneously in all directions. We denote the resulting space by
\begin{equation*}
	{HY} = [\bbR^+_\tau \times Y^2; \{\tau=0\} \times \diag_Y]
\end{equation*}
and can obtain coordinates by
\begin{equation*}
	 R = \sqrt{ \tau^2 + |\zeta-\zeta'|^2}, \quad
	 \theta = ( \theta_t,  \theta_{Y}) = \Big(\frac{\tau}{ R}, \frac{\zeta-\zeta'}{ R}\Big), \quad 
	\zeta'.
\end{equation*}
In these coordinates the pull-back of the heat kernel is given by
\begin{equation*}
	\frac1{(4\pi  R^2 \theta_t^2)^{m/2}} \exp\Big( -\frac{| \theta_Y|^2}{4 \theta_t^2} \Big) \; d\zeta'
	\in R^{-m} \CI(HY;\beta_R^*\density_Y).
\end{equation*}
We will make use of both radial blow-ups and parabolic blow-ups in constructing the heat space for the sub-Riemannian limit.
For a discussion of radial and parabolic blow-up we refer back to section \ref{sec:Corners} and the references therein.

\subsection{The heat kernel of a Laplace-type operator} \label{sec:EllHeat}
$ $

Let $M$ be a closed $m$-dimensional manifold and let $\Delta$ denote a Laplace-type operator on sections of a vector bundle $E \lra M.$ 
The heat equation for $\Delta,$ with initial data $f \in \CI(M;E),$ is 
\[ \begin{cases} (\pa_t+\Delta)u=0, \\ u\rvert_{t=0} =f. \end{cases} \]
We denote the linear operator $f\mapsto u$ by $e^{-t\Delta} :\CI(M\times \Rp; E)\to \CI(M;E).$ Its Schwartz kernel 
\begin{equation*}
	\cK_H = \wt\cK_H \cdot \mu_R\in \CI(\bbR^+_t \times M^2; \Hom(E)\otimes \beta_R^*\density_M),
\end{equation*}
where $\mu_R$ is the pull-back of a density to $\bbR^+_t \times M^2$ from the right factor of $M,$
satisfies the equation
\begin{equation*}
	\begin{cases}
	(\pa_t + \Delta_{\zeta})\wt \cK_H(t, \zeta, \zeta') =0\\
	\wt\cK_H(t, \zeta, \zeta')|_{t=0} = \delta_{\diag_M} 
	\end{cases}
\end{equation*}
with $\delta_{\diag_M}$ the delta distribution on the diagonal of $M^2$ at $t=0.$

It is convenient to multiply the first equation by $t,$ 
\begin{equation*}
	(t\pa_t + t\Delta_\zeta)\cK_H(t, \zeta, \zeta')=0,
\end{equation*}
as this does not change the solution and allows us to work with vector fields that are tangent to $\{t=0\}.$
Following \cite[Chapter 7]{Melrose:APS}, we establish the asymptotics of $\cK_H$ and show that it is $\cI$-smooth density on a suitable manifold with corners by constructing the solution to the heat equation within this class of densities. 

We construct the `heat space' of $M$ by first blowing-up $\{t=0\}$ parabolically, i.e., introducing $\tau=\sqrt t$ as a smooth function. We will not reflect this in our notation beyond using $\tau$ instead of $t.$ Secondly we perform the radial blow-up of the diagonal at time zero,
\begin{equation*}
	HM = [\bbR^+_\tau \times M^{ 2}; \{\tau=0,\zeta=\zeta'\}],
\end{equation*}
also represented schematically by Figure \ref{fig:simple heat space}. This space comes with a blow-down map,
\begin{equation*}
	\beta: HM \lra \bbR^+_\tau \times M^{ 2}.
\end{equation*}
We denote the composition of $\beta$ with the projection of $M^2$ to the right factor of $M$ by $\beta_{R}.$
This space has two boundary hypersurfaces,
\begin{equation*}
\begin{gathered}
	\bhs{tf} = \text{`temporal face'} = \bar{\beta^{-1}(\{\tau=0, \zeta\neq\zeta'\})}, \\
	\bhs{\bbE} = \text{`Euclidean face'} = \beta^{-1}(\{\tau=0, \zeta=\zeta'\}).
\end{gathered}
\end{equation*}
The latter can be identified with a fiberwise compactification of the normal bundle to the diagonal, i.e., the tangent bundle of $M,$ and indeed $\beta$ restricted to $\bhs{\bbE}$ is the bundle map.

As $t \to 0$ we expect, by analogy with the Euclidean heat kernel, to have exponential decay at $\bhs{tf}$ and interesting asymptotics only at $\bhs{\bbE},$ so we carry out our construction using the function space $\dCI_{\bbE}(M^2;\Hom(E)\otimes \beta_R^*\density_M)$ discussed in \S\ref{sec:Corners}, consisting of sections that vanish to infinite order at $\bhs{tf}.$
Local coordinates near $\bhs{\bbE}$ are obtained from local coordinates $\zeta$ on $M$ by, e.g., 
\begin{equation*}
	\tau, \quad \wt \omega = \frac{\zeta-\zeta'}\tau, \quad \zeta', 
\end{equation*}
(valid away from $\bhs{tf}$) in which $\tau$ is a boundary defining function for $\bhs{\bbE}.$ 
Since $\beta$ restricts to a diffeomorphism between $HM\setminus\bhs{\bbE}$ and $\bbR^+_{\tau}\times M^2 \setminus (\{\tau=0\} \times \diag_M)$ we can pull-back vector fields and differential operators along $\beta.$ 
The lift of $t\pa_t$ is given, in these coordinates, by
\begin{equation*}
	\beta^*(t\pa_t) = \tfrac12(\tau\pa_{\tau} - \wt \omega \cdot \pa_{\wt \omega}) 
\end{equation*}
and, for any vector field $V = \sum a_j(\zeta)\pa_{\zeta_j},$ the lift of $\tau V$ along $\beta_L$ is 
\begin{equation*}
	\beta_L^*(\tau V) = \sum a_j(\zeta'+\wt \omega\tau)\pa_{\wt\omega_j}.
\end{equation*}
Note that the restriction of $\beta_L^*(\tau V)$ to the fiber over $q \in M$ of $\bhs{\bbE}$ is the constant coefficient vector field obtained from $V$ by freezing coefficients at $q;$ this is none other than the symbol of $V.$

\begin{remark}
There is an equivalent way of obtaining the model operator of a differential operator $L$ at the fiber of $\bhs{\bbE}$ over $q\in M$ (cf. \cite[(2.5)]{Mazzeo:Hodge}).
Choose a chart $\phi$ around $q$ mapping into $T_qM$ with $\phi(q)=0,$ $D\phi(q)=\Id.$
Let $\cD_\tau$ denote the dilation by $\tau$ on $T_qM$ and define $N_q(L)u,$ for say $u$ a Schwartz section of $E_q$ over $T_qM,$ by 
\begin{equation*}
	N_q(L)u = \lim_{\tau \to 0} \cD_\tau^*\phi^*L(\phi^{-1})^*\cD_{1/\tau}^*u.
\end{equation*}
It is easy to see that, just as above, $N_q(\tfrac12\tau\pa_{\tau})u = -\tfrac12\wt \omega\cdot\pa_{\wt \omega}u$ and  $N_q(a(\zeta)\pa_{\zeta_j})u = a(0)\pa_{\wt \omega_j}u.$
Expressing this limit using polar coordinates we see that it is the same as pulling-back along $\beta$ and restricting to $\bhs{\bbE}.$
\end{remark}

It follows that the lift of $\tau^2\Delta$ also restricts to its principal symbol on $\bhs{\bbE},$ which in appropriate coordinates in each fiber is just the Euclidean Laplacian, acting as a fiber-wise differential operator.
This suggests that we take as our first approximation to the heat kernel of $\Delta$ the density we get from 
\eqref{eq:HeatKerY} written in projective coordinates,
\begin{multline}\label{eq:G0E}
	G_0(\tau, \omega, \zeta') = \chi(\tau)\tau^{-m} \frac1{(4\pi)^{m/2}} \exp\Big( -\frac14 |\omega|_{\sigma(\Delta)(\zeta')}^2 \Big)\Id_E \;  \mu_R \\
	\in \tau^{-m}\dCI_{\bbE}(HM; \Hom(E) \otimes \beta_R^*\density_M),
\end{multline}
where $|\cdot|_{\sigma(\Delta)(\zeta')}$ denotes the metric on the tangent space to $M$ at $\zeta'$ defined by the symbol of $\Delta,$  $\Id_E$ denotes the identity on $E,$
$\chi$ is a cut-off function equal to one in a neighborhood of $\bhs{\bbE},$ and $\mu_R$ is the pull-back of a non-degenerate density along $\beta_R.$ 

It follows that $G_0$ solves the heat equation to first order at $\bhs{\bbE},$
\begin{equation*}
	\beta_L^*(t\pa_t+t\Delta)G_0 \in \tau^{-m+1}\dCI_{\bbE}(HM;\Hom(E)\otimes \beta_R^*\density_M)
\end{equation*}
and it satisfies the initial condition $G_0 f \to f$ for all $f,$ i.e., $(\beta_L)_*G_0|_{t=0} = \delta_{\diag_M}.$
The next step is to find $G_j \in \dCI_{\bbE}(HM;\Hom(E)\otimes \beta_R^*\density_M),$ such that for any $N \in \bbN,$
\begin{equation*}
	\beta_L^*(t\pa_t+t\Delta)\sum_{j\leq N} \tau^{-m+j}G_j = \tau^{-m+1+N}R_N
	\in \tau^{-m+1+N}\dCI_{\bbE}(HM;\Hom(E)\otimes \beta_R^*\density_M).
\end{equation*}
Inductively, given $G_0, \ldots, G_N,$ we look for $G_{N+1}$ such that
\begin{multline*}
	\beta_L^*(t\pa_t+t\Delta)\sum_{j\leq N+1} \tau^{-m+j}G_j = \tau^{-m+1+N}R_N + \beta_L^*(t\pa_t+t\Delta)\tau^{-m+N+1}G_{N+1} \\
	\in \tau^{-m+2+N}\dCI_{\bbE}(HM;\Hom(E)\otimes \beta_R^*\density_M),
\end{multline*}
i.e., such that 
\begin{equation*}
	[\tfrac12(-\wt\omega\cdot\pa_{\wt\omega}-m+1+N) + \sigma(\Delta)] G_{N+1}|_{\tau=0} = -R_N|_{\tau=0}.
\end{equation*}
As $R_N|_{\tau=0}$ is a Schwartz function on the fibers of $\bhs{E},$ we can take Fourier transform and solve the resulting equation to get
\begin{equation*}
	\cF(G_{N+1}|_{\tau=0}) = \int_0^1 \exp((r-1)|\xi|^2)\cF(-R_N|_{\tau=0}) \; r^{N+1} \; dr,
\end{equation*}
which shows that $G_{N+1}|_{\tau=0}$ is a Schwartz function on the fibers of $\bhs{E}$ and completes the induction.

Once we have all of the $G_j$ we `asymptotically sum' them; that is we use Borel's Lemma to find 
\begin{equation*}
	G_{\infty} \in \tau^{-m}\CI(HM;\Hom(E)\otimes \beta_R^*\density_M) \text{ such that, for any $k,$ } G_{\infty} - \sum_{j \leq k}G_j = \cO(\tau^{-m+1+k}).
\end{equation*}
It follows that  
\begin{equation*}
	\beta_L^*(t\pa_t+t\Delta)G_\infty=R_{\infty} \in \dCI(HM;\Hom(E)\otimes \beta_R^*\density_M)
\end{equation*}
where $\dCI(HM;\Hom(E)\otimes \cD)$ denotes sections of $\Hom(E)\otimes \beta_R^*\density_M$ that vanish to infinite order at both $\bhs{tf}$ and $\bhs{\bbE}.$

For the final step, we allow $G_{\infty}$ and $R_{\infty}$ to act by convolution as Volterra operators.
Let $\dCI(\bbR^+_t \times M;E)$ denote the sections of $E$ that vanish to infinite order at $t=0$ (which is equivalent to vanishing to infinite order in $\tau$), and define
\begin{equation*}
\begin{gathered}
	G_{\infty}*: \dCI(\bbR^+_t \times M;E) \lra \dCI(\bbR^+_t \times M;E) \\
	(G_{\infty}*u)(t,\zeta) = \int_0^t \int_M G_{\infty}(t-s,\zeta, \zeta')  u(s,\zeta') \; ds 
\end{gathered}
\end{equation*}
This operator satisfies
\begin{equation*}
	(\pa_t + \Delta)(G_{\infty}*) = \Id + t^{-1}R_{\infty}*
\end{equation*}
and we can invert the right hand side (as a convolution operator) with a convergent Neumann series \cite[Proposition 7.17]{Melrose:APS}
\begin{equation*}
	(\Id + t^{-1}R_{\infty}*)^{-1}= \sum (-1)^k(t^{-1}R_{\infty}*)^k = \Id + S_{\infty}, \quad
	S_{\infty} \in \dCI(HM;\Hom(E)\otimes \beta_R^*\density).
\end{equation*}
It follows that the heat kernel of $\Delta$ (pulled back to $HM$) is equal to
\begin{equation*}
	\cK_{e^{-t\Delta}} = G_{\infty} + G_{\infty}*S_{\infty} \in \tau^{-m}\dCI_{\bbE}(HM;\Hom(E)\otimes \beta_R^*\density_M)
\end{equation*}
and, since $\cK_{e^{-t\Delta}} - G_{\infty} = G_{\infty}*S_{\infty} \in \dCI(HM;\Hom(E)\otimes\beta_R^*\density_M),$ we see that the 
Taylor expansion of the heat kernel at $\bhs{\bbE}$ is equal to the Taylor expansion of $G_{\infty},$ i.e., to the
one we constructed step by step above.

\subsection{The heat kernel of a contact Laplacian} \label{sec:CtHeat}
$ $

For $m=2n+1$ the Heisenberg group is the multiplicative subgroup of $\mathrm{GL}_m(\bbR)$ given by 
\begin{equation*}
	\bbH^{2n+1} = \lrbrac{
	\pmat{ 1 & x & z \\ 0 & \Id_n & y \\ 0 & 0 & 1 }: x,y \in \bbR^n, \; z \in \bbR },
\end{equation*}
and is naturally diffeomorphic to $\bbR^m.$ An important feature of $\bbH^m$ is that the dilations $(x,y,z)\mapsto(ax, by, cz)$ that are group homomorphisms are precisely those with $c=ab;$ the dilation with $a=b=c^{1/2} =\lambda>0$ is known as the `Heisenberg scaling' by $\lambda.$

In the construction of the heat kernel on a closed manifold we have used that the principal symbol of a differential operator is, at each point $p \in M,$ a translation-invariant differential operator on $T_pM\cong \bbR^m.$
On a contact manifold $M$ we can identify the tangent bundle with a bundle of Heisenberg groups (the osculating group of \cite{Folland-Stein}).
The `Heisenberg symbol' of a differential operator is, at each point $p \in M,$ a left-invariant differential operator on $T_pM \cong \bbH^m.$ 
A self-adjoint operator is said to be `Rockland' if its Heisenberg symbol is invertible and this is known to be the case for the Hodge Laplacians of the Rumin complex, $\Delta_{\sH}$ \cite[\S3]{Rumin:Formes}, \cite[Example 4.21]{Dave-Haller:Graded}. The heat kernel of a Rockland operator has been studied in, e.g., \cite{Beals-Greiner-Stanton, Dave-Haller:Heat}, \cite[\S5]{Ponge:Spectral}.
For future use we recast the construction of these heat kernels in parallel to the construction in \S\ref{sec:EllHeat}.

Let $P \in \Diff^{\ell}(M;E)$ be a non-negative self-adjoint Rockland operator of order $\ell$ (recall that $\Delta_{\sH}$ is of order four on $n$-forms and $n+1$-forms, and order two otherwise).
We construct the Heisenberg heat space of $M$ from $\bbR^+_t \times M^2$ by first introducing $\tau = t^{1/\ell}$ as a smooth function (i.e., by performing a quasi-homogeneous blow-up of $\{t=0\}$).
We will not include this blow-up in our notation below beyond using $\tau$ instead of $t.$
Secondly we blow-up the diagonal at time zero, but parabolically with respect to $\Ann(\sH),$
\begin{equation}\label{eq:CtHeatSpace}
	H_{\bbH}M = [\bbR^+_{\tau}\times M^2; \{\tau=0, \zeta=\zeta'\}, \Ann(\sH)].
\end{equation}
As before this space comes with a blow-down map
\begin{equation*}
	\beta:H_{\bbH}M \lra \bbR^+_{\tau}\times M^2
\end{equation*}
and has two boundary hypersurfaces,
\begin{equation*}
\begin{gathered}
	\bhs{tf} = \text{`temporal face'} = \bar{\beta^{-1}(\{\tau=0, \zeta\neq\zeta'\})},\\
	\bhs{\bbH} = \text{`Heisenberg face'} = \beta^{-1}(\{\tau=0, \zeta=\zeta'\}).
\end{gathered}
\end{equation*}
The latter is a fiberwise compactification of the tangent bundle of $M$ and $\beta$ restricts to the bundle map.

Note that the pull-back of the density bundle along the blow-down map satisfies
\begin{equation*}
	 \beta^*\density(\bbR_\tau^+\times M^2) = \rho_{\bbH}^{m+1} \;  \density(H_{\bbH}M),
\end{equation*}
explaining why the asymptotic expansion of the trace of the heat kernel of a contact Laplacian will begin at $-(m+1)$ instead of $-m.$ 

As above, we will work with the space $\dCI_{\bbH}(M^2;\Hom(E)\otimes \beta_R^*\density_M)$ of densities that vanish to infinite order at $\bhs{tf}$ but have non-trivial asymptotics at $\bhs{\bbH}.$ 
Local coordinates near $\bhs{\bbH}$ can be obtained from a Darboux chart $(x_j, y_j, z) = (\zeta, z)$ by, e.g.,
\begin{equation}\label{eq:CtProjCoords}
	\tau, \quad 
	\omega_\zeta = \frac{\zeta-\zeta'}\tau, \quad
	\omega_z = \frac{z-z'}{\tau^2}, \quad 
	\zeta', \quad 
	z',
\end{equation}
valid away from $\bhs{tf},$ in which $\tau$ is a boundary defining function for $\bhs{\bbH}.$ 

In these coordinates, the lift of $t\pa_t$ along $\beta$ is 
\begin{equation*}
	\beta^*(t\pa_t) 
	= \tfrac{1}{\ell}(\tau\pa_\tau - \omega_{\zeta} \cdot \pa_{\omega_{\zeta}} - 2\omega_z\pa_{\omega_z}) 
	= \tfrac{1}{\ell}(\tau\pa_\tau- \cR_{\bbH}),
\end{equation*}
and we recognize that $\cR_{\bbH}$ is the infinitesimal generator of dilations adapted to the Heisenberg group;
for a vector field $V = \sum a_j(\zeta,z) \pa_{\zeta_j} + a_0(\zeta,z) \pa_z$ the lift of $\tau V$ along $\beta$ is 
\begin{equation*}
	\beta^*(\tau V) = \sum a_j(\zeta'+\tau\omega_{\zeta}, z'+\tau^2\omega_z)\pa_{\omega_{\zeta_j}} 
	+  \tau^{-1} a_0(\zeta'+\tau\omega_{\zeta}, z'+\tau^2\omega_z)\pa_{\omega_{z}}.
\end{equation*}
Significantly, at the center of the Darboux chart, the lift of
$\pa_{x_j}-\tfrac12y_j\pa_z$ is $\pa_{\omega_{x_j}}-\tfrac12\omega_{y_j}\pa_{\omega_z}$ and the lift of 
$\pa_{y_j}+\tfrac12x_j\pa_z$ is $\pa_{\omega_{y_j}}+\tfrac12\omega_{x_j}\pa_{\omega_z};$
i.e., the lift of the horizontal vector fields on $M$ are the corresponding vector fields on $T_qM.$

\begin{remark}
Equivalently we can find the model operator of a differential operator $L$ at the fiber of $\bhs{\bbH}$ over $q \in M$ by choosing, e.g., a Darboux chart $\phi$ around $q$ mapping into $T_qM$ with $\phi(q)=0,$ $D\phi(q)=\Id.$
Let $\cD_\lambda$ denote the anisotropic dilation map
\begin{equation*}
	T_qM = \sH_q \oplus \sV_q \ni (\omega_\zeta, \omega_z) 
	\xmapsto{\phantom{x}\cD_{\lambda}\phantom{x}} (\lambda \omega_{\zeta}, \lambda^2\omega_z) \in \sH_q \oplus \sV_q = T_qM
\end{equation*}
and define $N_q(L)u,$ for say $u$ a Schwarz section of $E_q$ over $T_qM,$ by 
\begin{equation*}
	N_q(L)u = \lim_{\lambda\to0} \cD_\lambda^*\phi^*L(\phi^{-1})^*\cD_{1/\lambda}^*u.
\end{equation*}
It is easy to see that this gives the same operators as above, and coincides with other definitions of the Heisenberg symbol, cf. \cite{Julg-VanErp}, \cite[\S2.1]{Ponge:Spectral}.
\end{remark}

It follows that the lift of $\tau^\ell P$ restricts to each fiber of $\bhs{\bbH}$ to its Heisenberg symbol $\sigma^{\bbH}(P),$ e.g., the corresponding Hodge Laplacian on the Heisenberg group. Also we recognize the restriction of the lift of $t\pa_t$ to $\bhs{\bbH}$ as $-\ell^{-1}$ times the infinitesimal generator of the Heisenberg scaling.
It is shown in, e.g., \cite[Lemma 4]{Dave-Haller:Heat} that the heat kernel of $\sigma^{\bbH}(P)(q),$ $k_{\sigma^{\bbH}(P)(q)}(\tau, \omega_{\zeta}, \omega_z) \; \mu,$ is a Schwartz section of the density bundle on $T_qM,$ homogeneous with respect to the anisotropic dilation
\begin{equation*}
	k_{\sigma^{\bbH}(P)(q)}(\lambda \tau, \lambda \omega_{\zeta}, \lambda^2\omega_z) 
	= \lambda^{-(m+1)}k_{\sigma^{\bbH}(P)(q)}( \tau,  \omega_{\zeta}, \omega_z) 
\end{equation*}
whose integral over $T_qM$ is equal to the identity on $E_q.$
This suggests a first approximation to the heat kernel of $P$ analogous to \eqref{eq:G0E},
\begin{multline*}
	G_0(\tau, \omega_{\zeta}, \omega_z, \zeta', z')
	= \chi(\tau) \tau^{-(m+1)} k_{\sigma^{\bbH}(P)(\zeta')}(1, \omega_{\zeta}, \omega_z) \; \mu_R \\
	\in \tau^{-(m+1)}\dCI_{\bbH}(H_{\bbH}M; \Hom(E) \otimes \beta_R^*\density_M),
\end{multline*}
with $\chi$ a cut-off function equal to one in a neighborhood of $\bhs{\bbH}$, and $\mu_R$ is a density on $M$ pulled-back along $\beta_R.$
It follows that $G_0$ solves the heat equation to first order at $\bhs{\bbH},$
\begin{equation*}
	\beta_L^*(t\pa_t+t\Delta)G_0 \in \tau^{-(m+1)}\dCI_{\bbH}(H_{\bbH}M;\Hom(E)\otimes \beta_R^*\density_M)
\end{equation*}
and it satisfies the initial condition $\lim_{\tau\to0}G_0 f = f$ for all $f,$ i.e., $(\beta_L)_*G_0|_{t=0} = \delta_{\diag_M}.$
Since convolution on the Heisenberg group preserves Schwartz functions, we may construct $G_j$ as in \S\ref{sec:EllHeat} and asymptotically sum them to find
\begin{equation*}
\begin{gathered}
	G_{\infty} \in \tau^{-(m+1)}\dCI_{\bbH}(H_{\bbH}M; \Hom(E)\otimes \beta_R^*\density_M), \\
	\beta_L^*(t\pa_t + tP)G_{\infty} = R_{\infty} \in \dCI(H_{\bbH}M;\Hom(E)\otimes \beta_R^*\density_M).
\end{gathered}
\end{equation*}
Continuing as in \S\ref{sec:EllHeat}, we view this as a Volterra operator acting by convolution and invert it by a Neumann series to find
\begin{equation*}
	\cK_{e^{-tP}} \in \tau^{-(m+1)}\dCI_{\bbH}(H_{\bbH}M;\Hom(E)\otimes \beta_R^*\density_M).
\end{equation*}
$ $

Finally, we will have need of the Schwartz kernel of $Ae^{-tP}$ where $P$ is a non-negative self-adjoint Rockland operator of order $\ell$ and $A$ is a pseudodifferential operator of order zero in the Heisenberg calculus. The structure of these operators, when $P$ is elliptic, has been studied in several contexts as, e.g., a way of obtaining the Wodzicki-Guillemin residue of $A$ (see for example \cite{Loya:Res} for a nice overview, \cite{Lesch:Traces, Gil-Loya, Grubb-Seeley} for the use of $Ae^{-tP},$ and \cite{Ponge:NonCR, Ponge:NonCt, Ponge:NonH} for contact manifolds).
Fortunately, as we will explain, the case we will encounter below is simpler than the general case.

The Schwartz kernel of $Ae^{-tP}$ is again $\cI$-smooth on the heat space $H_{\bbH}M$ but the coefficients of the expansion at $\bhs{\bbH}$ are not necessarily local. Correspondingly it is convenient to understand $Ae^{-tP}$ as an integral transform of the resolvent of $P$ or of the complex powers of $P.$ Indeed, the Mellin transform of $Ae^{-tP}$ is equal to $AP^{-s},$ and, for $P$ a positive self-adjoint Rockland differential operator, the structure of $P^{-s}$ is detailed in \cite[\S5.3]{Ponge:Spectral}; namely,if $P$ has differential order $\ell,$ then $P^{-s}$ (and hence $AP^{-s}$) is a Heisenberg pseudodifferential operator of order $-s\ell$ and together these operators form a holomorphic family. By taking inverse Mellin transform, it follows that $Ae^{-tP}$ is an $\cI$-smooth function on $H_{\bbH}M,$ vanishing to infinite order at $\bhs{tf}$ and with asymptotic expansion at $\bhs{\bbH}$ of the form
\begin{equation}\label{eq:ExpAeP}
	\cK_{Ae^{-tP}} \sim \tau^{-(m+1)}\sum_{j\geq 0} a_j \tau^j
	+ 
	\sum_{k\geq 0} (b_k + \wt b_k \log \tau)\tau^{\ell k}, \Mas \tau \to 0,
\end{equation}
where we recall that $\tau = t^{1/\ell}.$ Recall from \S\ref{sec:Corners} that this is expressed as
\begin{equation*}
	Ae^{-tP} \in \cA_{phg}^{\cJ}(H_{\bbH}M;\Hom(E)\otimes \beta_R^*\density_M) \Mwith
	\cJ(\bhs{tf}) = \emptyset, \quad \cJ(\bhs{\bbH}) = -(m+1) \bar\cup 0,
\end{equation*}
where we use $-(m+1)$ and $0$ to stand for $-(m+1)+\bbN_0$ and $\bbN_0,$ respectively.
Interestingly, while $a_j$ and $\wt b_k$ are local (i.e., they each depend on finitely many terms in the asymptotic expansion into homogeneous terms of the symbols of $A$ and $P$), the terms $b_k$ need not be local. For example if $A$ is trace-class then the trace of $A$ will be equal to the integral of $a_{m+1}+b_0.$

It follows that the trace of $Ae^{-tP}$ has an expansion of the same form as \eqref{eq:ExpAeP} (the coefficient of the first log term is the noncommutative residue of $A$).
In our construction below, $A$ will be the projection $\Pi_{\sE_4}$ from \S\ref{sec:AsymLap} and $P$ will be $\Delta_{\sH}$ from \eqref{eq:DeltasH}. This allows us to use an observation of Branson-Gover \cite[\S 3]{Branson}, \cite{Branson-Gover} to rewrite the trace as follows (e.g., on forms of degree $n$)
\begin{multline}\label{eq:BransonTrick}
	\Tr(e^{-t\Delta_{\sH}}\rest{\sE_4^n}) 
	= \Tr_{\Omega^n_{\sH}M}(e^{-tD_{\sH}^*D_{\sH}}) \\
	= \Tr_{\Omega^n_{\sH}M}(e^{-t\Delta_{\sH}}) - \Tr_{\Omega^n_{\sH}M}(e^{-t(d_{\sH}\delta_{\sH})^2}) 
	= \Tr_{\Omega^n_{\sH}M}(e^{-t\Delta_{\sH}}) - \Tr_{\Omega^{n-1}_{\sH}M}(e^{-t(\delta_{\sH}d_{\sH})^2}) \\
	= \Tr_{\Omega^n_{\sH}M}(e^{-t\Delta_{\sH}}) 
	- \Tr_{\Omega^{n-1}_{\sH}M}(e^{-\Delta_{\sH}^2})
	+ \Tr_{\Omega^{n-2}_{\sH}M}(e^{-(d_{\sH}\delta_{\sH})^2}) \\
	= \ldots
	= \Tr_{\Omega^n_{\sH}M}(e^{-t\Delta_{\sH}}) 
	+ \sum_{k=1}^n (-1)^k \Tr_{\Omega^{n-k}_{\sH}M}(e^{-t\Delta_{\sH}^2}),
\end{multline}
which implies that the short-time asymptotic expansion of the trace does not have any $\log \tau$ terms or terms that are global in the manifold.
The analogous result holds for forms of degree $n+1.$

\section{The heat kernel of a sub-Riemannian limit} \label{sec:HeatSub}

In \cite[Theorems 3.5, 3.6]{Rumin:Sub}, Rumin  established convergence of the spectrum of the Hodge Laplacian of metrics undergoing a sub-Riemannian limit. Specifically he showed that, for some $\lambda \in \bbC,$ (in fact all $\lambda \in \bbC\setminus \bbR$)
\begin{equation*}
	\begin{cases}
	\displaystyle
	(\Delta_{\eps}-\lambda)^{-1} \to (\Delta_{\sH}|_{\sE_2^p}-\lambda)^{-1} & \Mon {}^{\sR}\Omega^pX, \quad p \notin\{n,n+1\} \\
	\\
	\displaystyle
	(\eps^{-2}\Delta_{\eps}-\lambda)^{-1} \to (\Delta_{\sH}|_{\sE_4^p}-\lambda)^{-1} & \Mon {}^{\sR}\Omega^pX, \quad p \in\{n,n+1\},
	\end{cases}
\end{equation*}
where $\sE_*^p$ are the spaces defined in \eqref{eq:EkDef}.
For positive time the behavior of the heat kernel of $\Delta_{\eps}$ as $\eps \to 0$ is entirely analogous, \cite[\S7]{Rumin:Sub}.
However the behavior of the heat kernel as both time and $\eps$ go to zero is a bit more intricate.
We will understand this degeneration by constructing a manifold with corners, the sub-Riemannian limit heat space, on which the heat kernel 
is $\cI$-smooth for a smooth index set $\cI.$\\

As sketched in \S\ref{sec:AsymHeat}, the heat kernel will have three interesting regimes as $\eps\to0,$ one where $\eps^2$ and $t$ go to zero at the same rate, one where $t$ stays bounded, and one where $\eps^2$ and $t^{-1}$ go to zero at the same rate. The latter only shows up when the form degree is $n$ or $n+1,$ so we will construct different heat spaces for forms in middle degrees, $\midH,$ and forms that are outside of middle degrees, $\nmidH.$

\subsection{The heat kernel outside of middle degree} 
$ $

In this section and the next we construct the heat kernels of the Hodge Laplacians $\Delta_{\eps}$ undergoing a sub-Riemannian limit. The construction is  parallel to those in \S\ref{sec:ModelHeat}; we construct an appropriate heat space and find a $\cI$-smooth density that solves the induced model problems at each boundary hypersurface, then improve this to an $\cI$-smooth density that solves the heat equation to infinite order at each boundary hypersurface, and finally use a Volterra series to solve away the remaining error.
In this section we work with ${}^{\sR}\Omega^pX$ for a fixed $p \notin \{n,n+1\}.$

To construct $\nmidH$ we start with $\bbR^+_t \times M^2 \times [0,1]_{\eps}$ and we first blow-up $\{t=0\}$ parabolically; i.e., we introduce $\tau = \sqrt t$ as our global time variable. We will not include this blow-up in the notation beyond using $\tau$ instead of $t.$
Next, to capture the asymptotics of the heat kernel as $\tau$ and $\eps$ both go to zero, we blow-up up the submanifold 
\begin{equation*}
	\{\tau=0\} \times \diag_M \times \{\eps=0\} \subseteq \bbR^+_\tau \times M^2 \times [0,1]_{\eps},
\end{equation*}
parabolically  in the directions of $\Ann(\sH),$ as in \eqref{eq:CtHeatSpace}. We denote the resulting boundary hypersurface by $\bhs{d,1}.$
Finally, to capture the asymptotics of the heat kernel as $\tau \to 0$ for positive $\eps,$ we blow-up  the (interior lift of the) diagonal of $M^2$ at time zero for all $\eps,$ and denote the resulting boundary hypersurface by $\bhs{d,0}.$

\begin{figure}[ht]
    \def\svgwidth{\columnwidth}
    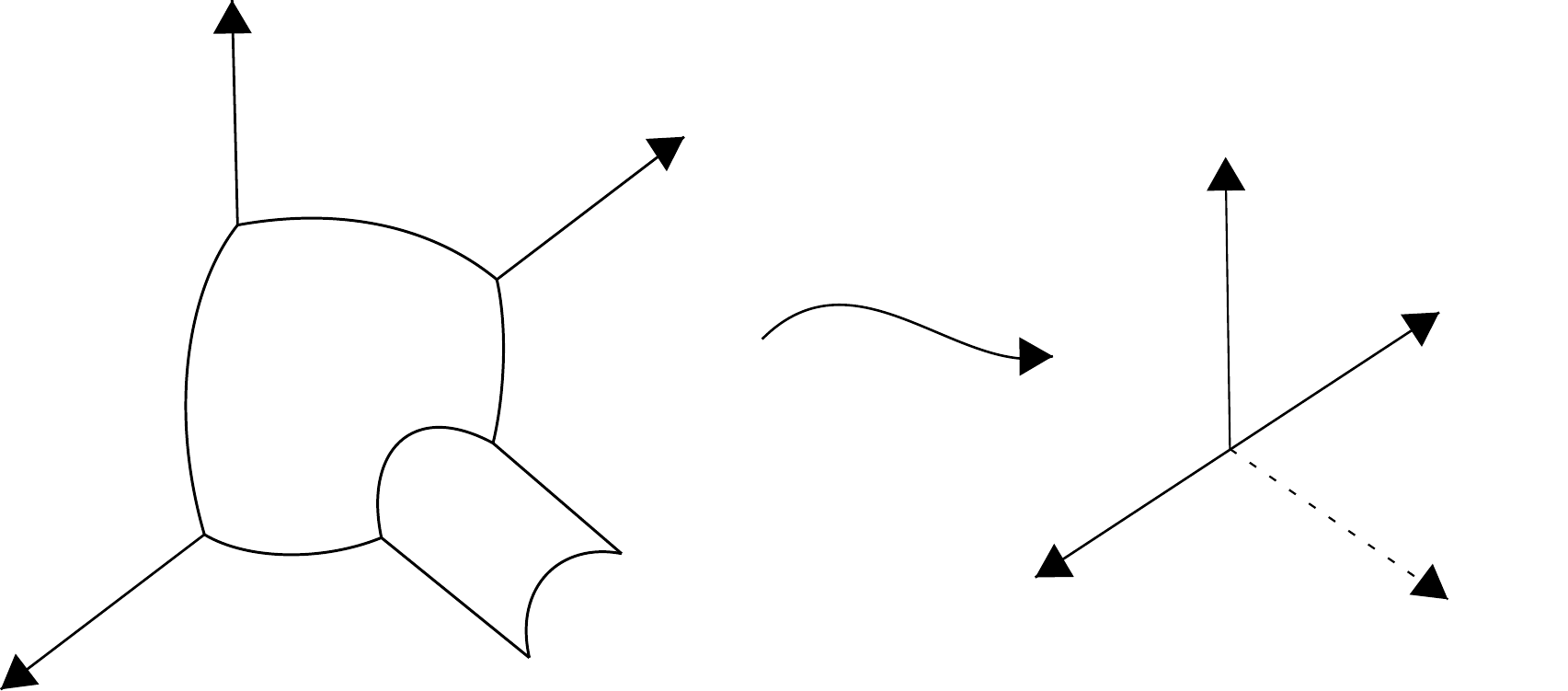

    \caption{The heat space $H_{\sR}^{\nmid} X$ and its blow-down map}
    \label{heatspace}
\end{figure}

Thus altogether we have
\begin{equation*}
\begin{gathered}
\begin{multlined}
	\nmidH
	= [ \bbR^+_\tau \times M^2 \times [0,1]_{\eps};
	\{\tau=0\} \times \diag_M \times \{\eps=0\}, \Ann(\sH); \\
	\{\tau=0\} \times \diag_M \times [0,1]_{\eps}], 
\end{multlined} 
\end{gathered}
\end{equation*}
together with its blow-down map
\begin{equation*}
	\beta: \nmidH \lra \bbR^+_\tau \times M^2 \times [0,1]_{\eps},
\end{equation*}
see Figure \ref{heatspace}.

Ignoring, as we will, the boundary hypersurface $\{\eps=1\},$ $\nmidH$ has four boundary hypersurfaces
\begin{equation*}
\begin{aligned}
	\bhs{0,1}(\nmidH) 
	&= \text{`$\sE_2$ face'}
	= \bar{\beta^{-1}(\{\eps=0\} \setminus \{\tau=0\} \times \diag_M)},  \\
	\bhs{d,1}(\nmidH) 
	&= \text{`$\sE_0$ face'}
	= \beta^{-1}(\{\tau=0\} \times \diag_M \times \{\eps=0\}), \\
	\bhs{1,0}(\nmidH) 
	&= \text{`temporal face'}  = \bar{\beta^{-1}(\{\tau=0\}\times (M^2\setminus \diag_M) \times [0,1]_{\eps})}, \\
	\bhs{d,0}(\nmidH) 
	&= \text{`Euclidean face'}
	= \bar{\beta^{-1}(\{\tau=0\} \times \diag_M \times (0,1]_{\eps})}.
\end{aligned}
\end{equation*}
where in the notation $\bhs{i,j}$ the first subindex is different from zero when the face lies over $\tau=0$ and the second subindex is different from zero when the face lies over $\eps=0.$
The restriction of the blow-down map to the boundary faces produced by blow-up induces fiber bundle structures,
\begin{equation*}
	\beta:\bhs{d,0}(\nmidH) \lra M \times [0,1]_{\eps}, \quad
	\beta:\bhs{d,1}(\nmidH)  \lra M.
\end{equation*}
The first of these is the fiberwise compactification of the sub-Riemannian limit tangent bundle ${}^{\sR}TX;$
in the second case the fiber over $q \in M$ is a compactification of $T_qM \oplus \bbR^+$ consistent with the Heisenberg scaling in $T_qM.$

We will denote a boundary defining function for, e.g., $\bhs{d,0}(\nmidH)$ by $\rho_{d,0},$ and similarly for the other boundary hypersurfaces.\\

Having constructed the heat space, we now focus on solving the model heat equations at each boundary hypersurface and constructing a parametrix for the heat kernel.
We will use $\beta_R^*\density_M$ to denote the density bundle of $M$ pulled-back to $\nmidH$ along the blow-down map followed by the projection onto the right factor of $M.$

The face $\bhs{d,0}(\nmidH)$ can be treated just as in \S\ref{sec:EllHeat}.
Indeed, this face fibers over $X = M \times [0,1]_{\eps}$ and the fiber over each point is a compactification of the corresponding fiber of the sub-Riemannian limit tangent bundle ${}^{\sR}TX.$
The model operator of $t\Delta_{\eps}$ is its principal symbol as a constant coefficient operator on the fiber of ${}^{\sR}TX.$ For horizontal vector fields, this is the usual symbol, while for vertical vector fields we divide out by $\eps$ and then fix coefficients (more correctly, in $\eps \pa_z$ as a section of ${}^{\sR}TX,$ the $\eps$ is not a coefficient but an inseparable part of the section).

Analogously to \S\ref{sec:EllHeat}, let $G_0$ be the section of $\Hom({}^{\sR}\Omega^pX)\otimes \beta_R^*\density_M$ which, in coordinates valid for $\eps>0,$
\begin{equation*}
	\tau, \quad
	\wt\omega_\zeta' = \frac{\zeta'-\zeta}{\tau}, \quad
	\wt\omega_z' = \frac{z'-z}{\tau}, \quad
	\zeta, \quad
	z, \quad
	\eps
\end{equation*}
is given by 
\begin{equation*}
	\frac\chi{(4\pi)^{m/2}}  \exp\Big( -\frac14 |\wt \omega'|_{g_{\eps}(\zeta,z,\eps)}^2 \Big) \Id_{{}^{\sR}\Omega^*(X)} \; d\wt\omega_{\zeta}'\;d\wt\omega_z',
\end{equation*}
where $|\cdot|_{g_{\eps}(\zeta',z',\eps)}^2$ denotes the metric $g_{\eps}$ on the corresponding fiber of ${}^{\sR}TX,$
$\chi$ is a smooth cut-off function equal to one in a neighborhood of $\bhs{d,0}(\nmidH).$
Coordinates valid near the intersection $\bhs{d,0}(\nmidH)\cap \bhs{d,1}(\nmidH)$ include
\begin{equation*}
	\sigma = \frac{\tau}{\eps}, \quad
	\bar\omega_\zeta' = \frac{\zeta'-\zeta}{\tau}, \quad
	\bar\omega_z' = \frac{z'-z}{\eps\tau}, \quad
	\zeta, \quad
	z, \quad
	\eps
\end{equation*}
in which $G_0$ takes the form
\begin{equation*}
\begin{gathered}
	 \frac\chi{(4\pi)^{m/2}}  \exp\Big( -\frac14 |\bar \omega'|_{g_{\eps}(\zeta,z,1)}^2 \Big) \Id_{{}^{\sR}\Omega^p(X)} \; d\wt\omega_{\zeta}'\;d\wt\omega_z'\\
	=  \eps^{-(m+1)}\sigma^{-m} \frac\chi{(4\pi)^{m/2}}  \exp\Big( -\frac14 |\bar \omega'|_{g_{\eps}(\zeta,z,1)}^2 \Big) \Id_{{}^{\sR}\Omega^*(X)} \; \bar\mu_R.
\end{gathered}
\end{equation*}
Thus we have 
\begin{equation*}
\begin{gathered}
	G_0 \in \rho_{d,1}^{-(m+1)}\rho_{d,0}^{-m}\rho_{1,0}^{\infty}\CI(\nmidH;\Hom({}^{\sR}\Omega^pX)\otimes \beta_R^*\density_M), \\
	\beta_L^*(t\pa_t + t\Delta_{\eps})G_{0} \in \rho_{d,1}^{-(m+1)}\rho_{d,0}^{-m+1}\rho_{0,1}^{\infty}\CI(\nmidH;\Hom({}^{\sR}\Omega^pX)\otimes \beta_R^*\density_M)
\end{gathered}
\end{equation*}
and we can proceed to remove the subsequent errors at $\bhs{d,0}(\nmidH)$ just as in \S\ref{sec:EllHeat} to find,
\begin{equation}\label{eq:Ginfty}
\begin{gathered}
	G_{\infty} \in \rho_{d,1}^{-(m+1)}\rho_{d,0}^{-m}\rho_{1,0}^{\infty}\CI(\nmidH;\Hom({}^{\sR}\Omega^pX)\otimes \beta_R^*\density_M),\\
	\beta_L^*(t\pa_t + t\Delta_{\eps})G_{\infty} 
		\in \rho_{d,1}^{-(m+1)}(\rho_{d,0}\rho_{0,1})^{\infty}\CI(\nmidH;\Hom({}^{\sR}\Omega^pX)\otimes \beta_R^*\density_M).
\end{gathered}
\end{equation}
This gives a parametrix for the heat equation with error vanishing to infinite order at the `Euclidean face', $\bhs{d,0}(\nmidH).$ \\

Next let us consider the boundary hypersurface $\bhs{d,1}(\nmidH),$ which fibers over $M$ and on which we obtain a model problem on each of these fibers.
This model problem takes a different aspect depending on the choice of coordinates.
Consider first coordinates projective with respect to $\tau,$ as in \eqref{eq:CtProjCoords},
\begin{equation*}
	\tau, \quad 
	\omega_\zeta = \frac{\zeta-\zeta'}\tau, \quad
	\omega_z = \frac{z-z'}{\tau^2}, \quad 
	\zeta', \quad 
	z', \quad
	\alpha = \frac\eps\tau,
\end{equation*}
valid away from $\bhs{1,0}(\nmidH)\cup \bhs{d,0}(\nmidH),$ in which $\tau$ is a boundary defining function for $\bhs{d,1}(\nmidH).$
In these coordinates, the fiber of the interior of $\bhs{d,1}(\nmidH)$ above $(\zeta',z') \in M$ is $\bbH^{m}_{(\omega_{\zeta},\omega_z)} \times \bbR^+_{\alpha},$
and the model operator of\footnote{In this section, to avoid confusion, we will sometimes denote the Laplacian on $M$ by $\Delta^M$ and the Laplacian on the Heisenberg group by $\Delta^{\bbH}.$}  $t\Delta_{\eps}^M$ is $\Delta_{\alpha}^{\bbH},$ i.e., the sub-Riemannian limit of the Hodge Laplacian on the Heisenberg group parametrized by $\alpha.$
So at this face, in these coordinates, we still have to deal with a sub-Riemannian limit, albeit on a simpler space.

Let us consider instead coordinates near $\bhs{d,1}(\nmidH)$ that are projective with respect to $\eps$ such as
\begin{equation*}
	\sigma = \frac{\tau}{\eps}, \quad
	\theta_{\zeta} = \frac{\zeta-\zeta'}{\eps}, \quad
	\theta_z = \frac{z-z'}{\eps^2}, \quad
	\zeta', \quad
	z', \quad
	\eps,
\end{equation*}
valid away from $\bhs{0,1}(\nmidH),$ in which $\eps$ is a boundary defining function for $\bhs{d,1}(\nmidH).$
In these coordinates, the fiber of the interior of $\bhs{d,1}(\nmidH)$ above $(\zeta',z') \in M$ is $\bbH^{m}_{(\theta_{\zeta},\theta_z)} \times \bbR^+_{\sigma},$
and the model operator of $t(\pa_t+\Delta_{\eps}^M)$ is $\tfrac12\sigma\pa_\sigma + \sigma^2\Delta_1^{\bbH},$ i.e., the heat equation for the Hodge Laplacian on the Heisenberg group.
The solution of the model heat problem on each fiber is thus, with notation as in \S\ref{sec:CtHeat}, 
\begin{equation*}
	k_{\Delta^{\bbH}}(\sigma, \theta_{\zeta}, \theta_z) \; \mu_R.
\end{equation*}
In terms of the previous set of coordinates, this becomes
\begin{equation*}
	k_{\Delta^{\bbH}}(\alpha^{-1}, \alpha^{-1}\omega_{\zeta}, \alpha^{-2}\omega_{z}) \; \mu_R 
	= 
	\alpha^{m+1}k_{\Delta^{\bbH}}(1, \omega_{\zeta}, \omega_{z}) \; \mu_R
\end{equation*}
(note that $\alpha^{m+1}\mu_R$ is homogeneous of degree zero with respect to Heisenberg dilation).

Returning to the parametrix construction, we have so far constructed $G_{\infty}$ satisfying \eqref{eq:Ginfty}. Let
\begin{multline*}
	R_0 = \rho_{d,1}^{m+1}\beta_L^*(t\pa_t + t\Delta_{\eps})G_{\infty}|_{\bhs{d,1}(\nmidH)} \\
	\in (\rho_{d,0}\rho_{0,1})^{\infty}\CI(\bhs{d,1}(\nmidH); \Hom({}^{\sR}\Omega^pX)\otimes \beta_R^*\density_M).
\end{multline*}
In coordinates projective with respect to $\eps,$ $R_0$ vanishes to infinite order as $\sigma \to 0$ and hence we may find $\wt L_0$ in the same space
(by convolving $R_0$ with the heat kernel of the Hodge Laplacian on the Heisenberg group) such that $(\tfrac12\sigma\pa_{\sigma}+\sigma^2\Delta_1^{\bbH})(\wt L_0) = R_0.$ Thus if we extend $\rho_{d,1}^{-(m+1)}\wt L_0$ to 
\begin{equation*}
	L_0 \in  \rho_{d,1}^{-(m+1)}(\rho_{d,0}\rho_{0,1})^{\infty}\CI(\nmidH;\Hom({}^{\sR}\Omega^pX)\otimes \beta_R^*\density_M)
\end{equation*}
then we can add it to $G_{\infty}$ and solve the heat equation to one order better at $\bhs{d,1}(\nmidH),$ i.e.,
\begin{equation*}
	\beta_L^*(t\pa_t + t\Delta_{\eps})(G_{\infty} + L_0)
		\in \rho_{d,1}^{-(m+1)+1}(\rho_{d,0}\rho_{0,1})^{\infty}\CI(\nmidH;\Hom({}^{\sR}\Omega^pX)\otimes \beta_R^*\density_M).
\end{equation*}
Continuing in this way we can solve away the error of the parametrix at this face and find 
\begin{equation*}
\begin{gathered}
	L_\infty \in \rho_{d,1}^{-(m+1)}(\rho_{d,0}\rho_{0,1})^{\infty}\CI(\nmidH;\Hom({}^{\sR}\Omega^pX)\otimes \beta_R^*\density_M), \\
	R_{\infty} = \beta_L^*(t\pa_t + t\Delta_{\eps})(G_{\infty} + L_\infty)
		\in (\rho_{d,1}\rho_{d,0}\rho_{0,1})^{\infty}\CI(\nmidH;\Hom({}^{\sR}\Omega^pX)\otimes \beta_R^*\density_M).
\end{gathered}
\end{equation*}
This gives a parametrix for the heat equation with error vanishing to infinite order at the `Euclidean face', $\bhs{d,0}(\nmidH)$ and the `$\sE_0$ face', $\bhs{d,1}(\nmidH).$ \\

We next consider the situation at the face $\bhs{0,1}(\nmidH)$ which we can identify with the Heisenberg heat space of $M$ from \eqref{eq:CtHeatSpace}.
The initial condition for the model problem at this face is that as time goes to zero, the heat kernel must match the solution of the model problem already constructed at
$\bhs{d,1}(\nmidH).$ In projective coordinates with respect to $\eps,$ we found that the model operator of $t\Delta_{\eps}$ at $\bhs{d,1}(\nmidH)$ is the sub-Riemannian limit of the Hodge Laplacian on the Heisenberg group. Hence the initial problem for the heat equation at $\bhs{0,1}(\nmidH)$ is the projection onto the null space of $a_{-1},$ i.e., the projection onto $\sE_2$\footnote{If we were to use projective coordinates with respect to $\tau,$ this is a statement about the large-time limit of the heat kernel of the Hodge Laplacian on the Heisenberg group and was already noticed by Rumin \cite[Theorem 7.14]{Rumin:Sub}.}. The upshot is that for the model heat equation we may work in $\sE_2$ and the solution is then $e^{-t(d_{\sH}+\delta_{\sH})^2}.$

Returning to the parametrix construction, since the solutions to the model equations on $\bhs{d,1}(\nmidH)$ are kernels valued in $\sE_2$ at $\bhs{d,1}(\nmidH)\cap \bhs{0,1}(\nmidH),$ we may assume that $R|_{\bhs{0,1}(\nmidH)}$ is valued in $\sE_2.$ 
We may then proceed as explained in \S\ref{sec:AsymHeat} to solve the heat equation asymptotically at $\bhs{0,1}(\nmidH),$ i.e., to find
\begin{equation*}
\begin{gathered}
	P_{\infty} \in (\rho_{d,1}\rho_{d,0}\rho_{0,1})^{\infty}\CI(\nmidH;\Hom({}^{\sR}\Omega^pX)\otimes \beta_R^*\density_M),\\
\begin{multlined}
	\beta_L^*(t\pa_t + t\Delta_{\eps})(G_{\infty} + L_\infty+ P_{\infty})
		\in (\rho_{0,1}\rho_{d,1}\rho_{d,0}\rho_{0,1})^{\infty}\CI(\nmidH;\Hom({}^{\sR}\Omega^pX)\otimes \beta_R^*\density_M)\\
		=\dCI(\nmidH;\Hom({}^{\sR}\Omega^pX)\otimes \beta_R^*\density_M).
\end{multlined}
\end{gathered}
\end{equation*}

Finally, we can remove the remaining error by considering $Q = G_{\infty} + L_\infty+ P_{\infty}$ and $S = -\beta_L^*(t\pa_t + t\Delta_{\eps})(Q)$ as Volterra operators acting by convolution. The section $Q$ satisfies
\begin{equation*}
	\beta_L^*(\pa_t + \Delta_{\eps})(Q\star) = \Id - \tfrac1tS\star
\end{equation*}
and we can invert the right hand side using the convergent Neumann series
\begin{equation*}
	(\Id-\tfrac1tS\star)^{-1} = \Id +\sum_{j\geq1} (\tfrac1tS)^{\star j} = \Id + S_{\infty}\star,
	\quad
	S_{\infty} \in \dCI(\nmidH;\Hom({}^{\sR}\Omega^pX)\otimes \beta_R^*\density_M).
\end{equation*}
It follows that the sub-Riemannian limit heat kernel, for differential forms outside of middle degrees is given by
\begin{equation*}
	\cK_{e^{-t\Delta_{\eps}}} = Q + Q \star S_{\infty} 
	\in \rho_{d,1}^{-(m+1)}\rho_{d,0}^{-m}\rho_{1,0}^{\infty} \CI(\nmidH;\Hom({}^{\sR}\Omega^pX)\otimes \beta_R^*\density_M),
\end{equation*}
and satisfies 
\begin{equation*}
\begin{gathered}
	\beta_L^*(t\pa_t + t\Delta_{\eps})\cK_{e^{-t\Delta_{\eps}}} = 0, \\
	\cK_{e^{-t\Delta_{\eps}}}\rvert_{\bhs{0,1}(\nmidH)} = \cK_{e^{-t\Delta^M_{\sH}}\Pi_{\sE_2}}, \quad
	\rho_{d,1}^{(m+1)}\cK_{e^{-t\Delta_{\eps}}}\rvert_{\bhs{d,1}(\nmidH)} = k_{e^{-t\Delta^{\bbH}_{\sH}}}, \\
	\rho_{d,0}^{m}\cK_{e^{-t\Delta_{\eps}}}\rvert_{\bhs{d,0}(\nmidH)} = 
	\frac {\eps^{-(m+1)}}{(4\pi)^{m/2}}  \exp\Big( -\frac14 |\wt \omega|_{g_{\eps}(\zeta',z',1)}^2 \Big) \Id_{{}^{\sR}\Omega^p(X)} \; \bar\mu_R,
\end{gathered}
\end{equation*}
as required.
This proves part (i) of Theorem \ref{thm:IntroThm} for differential forms in degree $p \notin\{n,n+1\}.$

\subsection{The heat kernel in middle degrees} \label{sec:HeatMid}
$ $

In this section we first discuss how to modify the construction of the heat kernel of $\Delta_{\eps}$ in the previous section for forms ${}^{\sR}\Omega^pX,$ with $p \in \{n,n+1\}.$ Then we will discuss the how to construct the heat kernel of $\eps^{-2}\Delta_{\eps}$ and why this is necessary for a uniform description of the long-time behavior.\\

For the Hodge Laplacian $\Delta_{\eps}$ in middle degrees, the construction in the previous section of a parametrix for $t\pa_t + t\Delta_{\eps}$ works without change at the faces $\bhs{d,0}(\nmidH)$ and $\bhs{d,1}(\nmidH).$ The model heat equation at $\bhs{0,1}(\nmidH)$ has initial condition $\Pi_{\sE_2}$ and so we may work in $\sE_2$ and the model heat equation to solve is then 
$(t\pa_t + t(d_{\sH}+\delta_{\sH})^2)\Pi_{\sE_2},$ whose solution is 
\begin{equation*}
	e^{-t(d_{\sH}+\delta_{\sH})^2\Pi_{\sE_2}}
	=\Pi_{\sG_2}e^{-t(d_{\sH}+\delta_{\sH})^2}\Pi_{\sG_2} + \Pi_{\sE_4}
	=\begin{cases}
	\Pi_{\sG_2}e^{-t(d_{\sH}\delta_{\sH})}\Pi_{\sG_2} + \Pi_{\sE_4} & \Mif q=n \\
	\Pi_{\sG_2}e^{-t(\delta_{\sH}d_{\sH})}\Pi_{\sG_2} + \Pi_{\sE_4} & \Mif q=n+1
	\end{cases}
\end{equation*}
since $d_{\sH}|_{\Omega_{\sH}^n}=0=\delta_{\sH}|_{\Omega_{\sH}^{n+1}}.$
To understand the structure of this solution, recall that Rumin showed that $\Delta_{\sH,n} = (d_{\sH}\delta_{\sH})^2 + D_{\sH}^*D_{\sH}$ is Rockland and hence it has a generalized inverse in the Heisenberg calculus, which we denote $\Delta_{\sH,n}^{\dagger},$ such that
\begin{equation*}
	\Delta_{\sH,n}\Delta_{\sH,n}^{\dagger} = \Delta_{\sH,n}^{\dagger}\Delta_{\sH,n} = \Id - \Pi_{\Ker(\Delta_{\sH,n})}.
\end{equation*}
The decomposition
\begin{equation*}
	\Omega^n_{\sH}M 
	= \Ker(\Delta_{\sH,n}) \oplus \Im(\Delta_{\sH,n})
	= \Ker(\Delta_{\sH,n}) \oplus \Im((d_{\sH}\delta_{\sH})^2) \oplus \Im(D_{\sH}^*D_{\sH})
\end{equation*}
is preserved by $\Delta_{\sH,n}^{\dagger}$ and hence 
\begin{equation*}
	\Pi_{\sG_2} = (d_{\sH}\delta_{\sH})^2\Delta_{\sH,n}^{\dagger}, \quad
	\Pi_{\sE_4} = \Id_{\Omega^n_{\sH}M} - (d_{\sH}\delta_{\sH})^2\Delta_{\sH,n}^{\dagger}
\end{equation*}
are both Heisenberg pseudodifferential operators of order $0$ (they add up to the identity on $\sE_2^n = \Omega^n_{\sH}M$).
Moreover, the operator $(d_{\sH}\delta_{\sH} + D_{\sH}^**)$ is Rockland, as it squares to $\Delta_{\sH,n},$ and we can write
\begin{equation*}
	e^{-t(d_{\sH}+\delta_{\sH})^2\Pi_{\sE_2}}
	= \Pi_{\sG_2}e^{-t(d_{\sH}\delta_{\sH} + D_{\sH}^**)}\Pi_{\sG_2} + \Pi_{\sE_4} \text{ on } \Omega_{\sH}^nM
\end{equation*}
and analogously for $q=n+1.$
Thus we can see that the first term in this equality is smooth on $\bhs{0,1}(\nmidH)$ but the second term is not.

\begin{figure}[ht]
    \def\svgwidth{\columnwidth}
    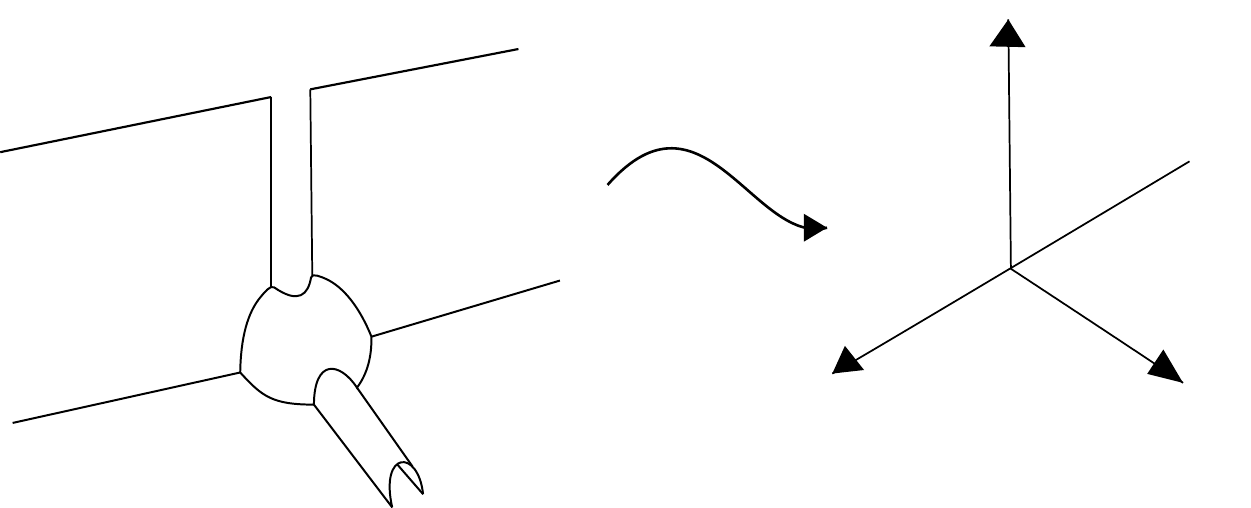

    \caption{The heat space $H_{\sR}^{+} X$ and its blow-down map.}
    \label{plus heatspace}
\end{figure}

As we know that on $\sE_4,$ $\Delta_{\eps}\circ \Phi_4 = \eps^2\Delta_{\sH} \Pi_{\sE_4},$ we see that we should blow-up the diagonal at $\eps=0$ to capture the asymptotics of this term.
Let
\begin{equation*}
	\pmidH
	=[\nmidH; \bbR^+_{\tau} \times \diag_M \times \{\eps=0\}, \mathrm{Ann}(\sH)].
\end{equation*}
Together with the lifts of the boundary hypersurfaces of $\nmidH,$ 
\begin{equation*}
	\bhs{0,1}(\pmidH), \quad 
	\bhs{d,1}(\pmidH), \quad 
	\bhs{1,0}(\pmidH), \quad
	\bhs{d,0}(\pmidH), \quad
\end{equation*}
there is the boundary hypersurface produced by the new blow-up
\begin{equation*}
	\bhs{d,2}(\pmidH) = \text{`local $\sE_2$ face'} = \bar{\beta^{-1}(\bbR^+_{\tau} \times \diag_M \times \{\eps=0\})},
\end{equation*}
see Figure \ref{plus heatspace}.

The boundary hypersurface $\bhs{d,2}(\pmidH)$ fibers over $M$ with fiber over $p \in M$ given by $\bbR^+_{\tau}$ times a parabolic compactification of $T_pM.$
Identifying $T_pM$ with the Heisenberg group $\bbH^m,$ the lift of $\eps^2\Delta_{\sH} \Pi_{\sE_4}$ is $\Delta_{\sH}^{\bbH} \Pi_{\sE_4^{\bbH}}$ and so the model heat equation is $\tfrac12\tau\pa_\tau + \tau^2 \Delta_{\sH}^{\bbH}\Pi_{\sE_4^{\bbH}}.$ The solution to this equation is
\begin{equation*}
	e^{-\tau^2\Delta_{\sH}^{\bbH}\Pi_{\sE_4^{\bbH}}} = \Pi_{\sE_4^{\bbH}}\; e^{-\tau^2\Delta_{\sH}^{\bbH}}\; \Pi_{\sE_4^{\bbH}}.
\end{equation*}
We can solve the heat equation $\tfrac12\tau\pa_{\tau}+ \tau^2\eps^2\Delta_{\sH}\Pi_{\sE_4}$ to infinite order at $\bhs{d,2}(\pmidH)$ by using the expansion \eqref{eq:ExpAeP}
and the equation $\tfrac12\tau\pa_{\tau}+\tau^2\Delta_{\sH}\Pi_{\sG_2}$ at $\bhs{0,1}(\pmidH)$ by using the construction described in $\S\ref{sec:AsymHeat}$ and the previous section.
Combining this with the constructions of $G_{\infty}$ and $L_{\infty}$ from the previous section, we obtain
\begin{equation*}
\begin{gathered}
	W_{\infty} \in \sA_{phg}^{\cJ}(\pmidH;\Hom({}^{\sR}\Omega^pX)\otimes \beta_R^*\density_M), \\
	\Mwith 
	\cJ(\bhs{d,1}(\pmidH)) = 
	\cJ(\bhs{d,0}(\pmidH)) = 
	\cJ(\bhs{1,0}(\pmidH)) = \emptyset, \\
	\cJ(\bhs{0,1}(\pmidH)) = 0, \quad
	\cJ(\bhs{d,2}(\pmidH)) = -(m+1) \bar\cup 0
\end{gathered}
\end{equation*}
such that the heat equation has remainder
\begin{equation*}
	 \beta_L^*(t\pa_t + t\Delta_{\eps})(G_{\infty} + L_{\infty} + W_{\infty}) 
		\in \CI(\pmidH;\Hom({}^{\sR}\Omega^pX)\otimes \beta_R^*\density_M).
\end{equation*}
We may use a Neumann series of Volterra operators as in the previous section to improve this to the heat kernel itself.
In this way we have shown that, for degrees $p \in \{n, n+1\},$
\begin{equation*}
\begin{gathered}
	\cK_{e^{-t\Delta_{\eps}}} 
	\in \sA_{phg}^{\cW}(\pmidH;\Hom({}^{\sR}\Omega^pX)\otimes \beta_R^*\density_M),\\
	\Mwith
	\cW(\bhs{1,0}(\pmidH)) = \emptyset, \quad
	\cW(\bhs{d,0}(\pmidH)) = -m, \quad
	\cW(\bhs{0,1}(\pmidH)) = 0, \\
	\cW(\bhs{d,1}(\pmidH)) = -(m+1), \quad
	\cW(\bhs{d,2}(\pmidH)) = -(m+1) \bar\cup 0,
\end{gathered}
\end{equation*}
satisfies 
\begin{equation*}
\begin{gathered}
	\beta_L^*(t\pa_t + t\Delta_{\eps})\cK_{e^{-t\Delta_{\eps}}} = 0, \\
	\cK_{e^{-t\Delta_{\eps}}}\rvert_{\bhs{0,1}(\pmidH)} 
		= \cK_{e^{-t\Delta^M_{\sH}}\Pi_{\sG_2}}, \\
	\rho_{d,1}^{(m+1)}\cK_{e^{-t\Delta_{\eps}}}\rvert_{\bhs{d,1}(\pmidH)} = k_{e^{-t\Delta^{\bbH}_{\sH}}}, \quad
	\rho_{d,2}^{(m+1)}\cK_{e^{-t\Delta_{\eps}}}\rvert_{\bhs{d,2}(\pmidH)} = k_{e^{-t\Delta^{\bbH}_{\sH}}\Pi_{\sE_4^{\bbH}}}, \\
	\rho_{d,0}^{m}\cK_{e^{-t\Delta_{\eps}}}\rvert_{\bhs{d,0}(\pmidH)} = 
	\frac {\eps^{-(m+1)}}{(4\pi)^{m/2}}  \exp\Big( -\frac14 |\wt \omega|_{g_{\eps}(\zeta',z',1)}^2 \Big) \Id_{{}^{\sR}\Omega^p(X)} \; \bar\mu_R.
\end{gathered}
\end{equation*}

It follows that, as $\eps \to 0$ and $t \to \infty,$ the heat kernel of $\Delta_{\eps}$ converges to the projection onto $\sE_4.$
Since this is infinite dimensional this is unsatisfactory for understanding the limit of, for example, the zeta function of the Hodge Laplacian.
The way around this, as already studied by Rumin \cite{Rumin:Sub}, is to construct the heat kernel of $\eps^{-2}\Delta_{\eps}$ instead.
Notice that the operator $t(\pa_t + \eps^{-2}\Delta_{\eps}),$ expressed in the rescaled time variable $T = \frac{t}{\eps^2}$ is $T(\pa_T + \Delta_{\eps})$
so solving the heat equation for $\eps^{-2}\Delta_{\eps}$ is the same as solving the heat equation for $\Delta_{\eps}$ but with a rescaled time.
This means that we have already done most of the work of solving this equation.\\

The heat kernel of $\eps^{-2}\Delta_\eps$ will be $\cI$-smooth on a different heat space, $\midH,$ obtained from $\pmidH$ by rescaling the time variable.
We again start with $\bbR^+_t\times M^2 \times [0,1]_{\eps}$ and replace $t$ with $\tau = \sqrt t.$
Then we blow-up the submanifold $\{\tau=\eps=0\}$ and denote the resulting boundary hypersurface by $\bhs{1,1}(\midH).$
Let $R = \sqrt{\eps^2 + \tau^2}$ and $\Theta = \arctan(\frac\eps\tau)$ denote polar coordinates on the resulting space.

Next we blow-up the submanifold $\{R=0, \Theta = \tfrac\pi2\}\times\diag_M$ parabolically in the directions of $\text{Ann}(\sH)$ and denote the resulting boundary hypersurface by $\bhs{d,1}(\midH).$ Coordinates valid near the interior of this face include
\begin{equation*}
	\sigma = \frac \tau{\eps^2}, \quad
	\omega_{\zeta} = \frac{\zeta-\zeta'}\eps, \quad
	\omega_{z} = \frac{z-z'}{\eps^2}, \quad
	\zeta', \quad
	z', \quad
	\eps,
\end{equation*}
in which $\eps$ is a boundary defining function for $\bhs{d,1}(\midH).$

Thirdly, we blow-up the interior lift of the submanifold $\{R=0\} \times \diag_M,$ parabolically in the directions of $\text{Ann}(\sH)$ and denote the resulting hypersurface by $\bhs{d,2}(\midH).$ Coordinates valid near the interior of this face include
\begin{equation*}
	s = \frac \tau{\eps}, \quad
	\omega_{\zeta} = \frac{\zeta-\zeta'}\eps, \quad
	\omega_{z} = \frac{z-z'}{\eps^2}, \quad
	\zeta', \quad
	z', \quad
	\eps,
\end{equation*}
in which $\eps$ is a boundary defining function for $\bhs{d,2}(\midH).$

Finally we blow-up the interior lift of the submanifold $\{\tau=0\}\times \diag_M \times [0,1]_{\eps}$ and denote the resulting boundary hypersurface by $\bhs{d,0}(\midH).$

\begin{figure}[ht]
    \def\svgwidth{\columnwidth}
    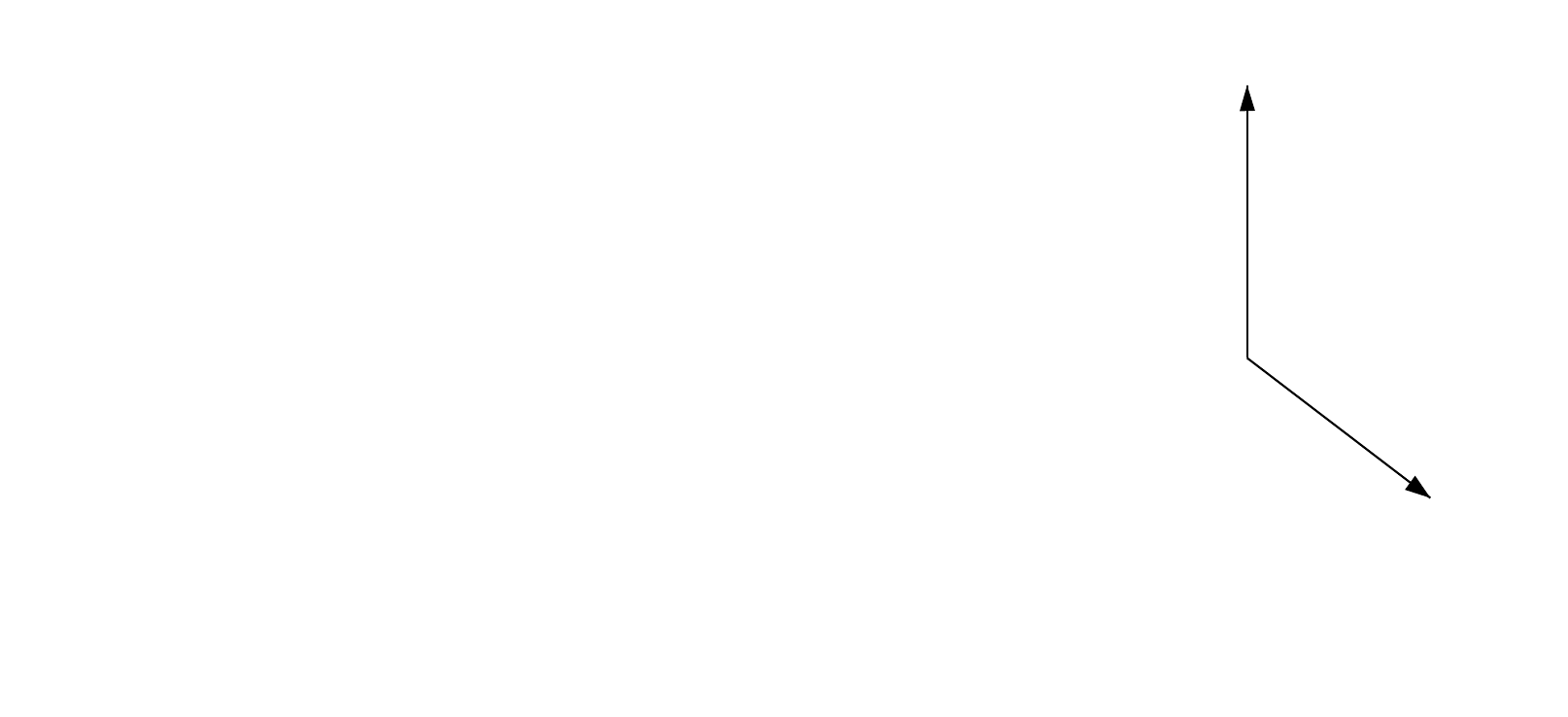

    \caption{The heat space $H_{\sR}^{\mid} X$ and its blow-down map}
    \label{mid heatspace}
\end{figure}

Thus altogether we have
\begin{multline*}
	\midH 
	= \Big[
	[\Rp_\tau \times M^2\times [0,1]_\eps; 
	\{\tau=0=\eps\}];
	\{R=0, \Theta = \tfrac\pi2\}\times\diag_M,\text{Ann}(\sH)\oplus\ang{dt}; \\
	\{R=0\}\times \diag_M, \text{Ann}(\sH);
	\{\tau=0\}\times \diag_M\times [0,1]_\eps
	\Big],
\end{multline*}
with blow-down map 
\begin{equation*}
	\beta: \midH \lra \bbR^+_{\tau} \times M^2 \times [0,1]_{\eps},
\end{equation*}
see Figure \ref{mid heatspace}.

There are six boundary hypersurfaces of $\midH$ in $\{\eps<1\},$\footnote{Note that we are using the same notation $\bhs{1,1}(\midH)$ for the boundary hypersurface of $\midH$ and for the boundary hypersurface produced by the blow-up of $\{\tau=\eps=0\}$ in $\bbR^+_\tau \times M^2 \times [0,1]_{\eps},$ and similarly for the other faces. We trust that this abuse of notation will not lead to confusion.}
\begin{equation*}
\begin{aligned}
	\bhs{0,1}(\midH) 
	&= \text{`$\sE_4$ face'}
	= \bar{\beta^{-1}(\{\eps=0\} \setminus \{\tau=0\} \times \diag_M)},  \\
	\bhs{1,1}(\midH) 
	&= \text{`$\sE_2$ face'} 
	= \bar{\beta^{-1}(\{\tau=0\} \times M^2\times \{\eps=0\} \setminus ( \diag_M\times \{\eps=0\} ) )}, \\
	\bhs{d,2}(\midH) 
	&= \text{`local $\sE_2$ face'}
	= \overline{ \beta^{-1}(\{R=0, \Theta=\tfrac\pi2\}\times \diag_M ) }, \\
	\bhs{d,1}(\midH) 
	&= \text{`$\sE_0$ face'} 
	= \overline{ \beta^{-1} (\{R=0\}\times\diag_M ) }, \\
	\bhs{1,0}(\midH) 
	&= \text{`temporal face'}  
	= \bar{\beta^{-1}(\{\tau=0\}\times (M^2\setminus \diag_M) \times [0,1]_{\eps})}, \\
	\bhs{d,0}(\midH) 
	&= \text{`Euclidean face'}
	= \bar{\beta^{-1}(\{\tau=0\} \times \diag_M \times (0,1]_{\eps})}.
\end{aligned}
\end{equation*}

The parametrices we constructed for $t(\pa_t + \Delta_{\eps})$ at all of the boundary hypersurfaces of $\pmidH$ carry over to $\midH$ by rescaling the time variable and yield a parametrix
\begin{multline*}
	V
	\in \sA_{phg}^{\cW}(\midH;\Hom({}^{\sR}\Omega^pX)\otimes \beta_R^*\density_M),\\
	\Mwith
	\cW(\bhs{1,0}(\midH)) = \emptyset, \quad
	\cW(\bhs{d,0}(\midH)) = -m, \\
	\cW(\bhs{0,1}(\midH)) = 
	\cW(\bhs{1,1}(\midH)) = 
	0, \\
	\cW(\bhs{d,1}(\midH)) = -(m+1), \quad
	\cW(\bhs{d,2}(\midH)) = -(m+1) \bar\cup 0, 
\end{multline*}
satisfying 
$\beta_L^*(t\pa_t + t\eps^{-2}\Delta_{\eps})V \in \dCI_{\bhs{0,1}}(\midH;\Hom({}^{\sR}\Omega^pX)\otimes \beta_R^*\density_M). $
Moreover we may assume that its restriction to $\bhs{0,1}(\midH)$ is valued in $\sE_4$ and proceed as explained in \S\ref{sec:AsymHeat} to solve the heat equation asymptotically at this face and find
\begin{equation*}
\begin{gathered}
	V' \in \dCI_{\bhs{0,1}}(\midH;\Hom({}^{\sR}\Omega^pX)\otimes \beta_R^*\density_M), \\
	\beta_L^*(t\pa_t + t\Delta_{\eps})(V + V') \in \dCI(\midH;\Hom({}^{\sR}\Omega^pX)\otimes \beta_R^*\density_M).
\end{gathered}
\end{equation*}
Again we can now use a Neumann series of Volterra operators to improve this to the heat kernel itself.
Thus we have shown that, for degrees $p \in \{n, n+1\},$
\begin{multline}\label{eq:HeatKerMid}
	\cK_{e^{-t\eps^{-2}\Delta_{\eps}}} 
	\in \sA_{phg}^{\cW}(\midH;\Hom({}^{\sR}\Omega^pX)\otimes \beta_R^*\density_M),\\
	\Mwith
	\cW(\bhs{1,0}(\midH)) = \emptyset, \quad
	\cW(\bhs{d,0}(\midH)) = -m, \\
	\cW(\bhs{0,1}(\midH)) = 
	\cW(\bhs{1,1}(\midH)) = 
	0, \\
	\cW(\bhs{d,1}(\midH)) = -(m+1), \quad
	\cW(\bhs{d,2}(\midH)) = -(m+1) \bar\cup 0, 
\end{multline}
satisfies 
\begin{equation*}
\begin{gathered}
	\beta_L^*(t\pa_t + t\eps^{-2}\Delta_{\eps})\cK_{e^{-t\eps^{-2}\Delta_{\eps}}} = 0, \\
	\cK_{e^{-t\Delta_{\eps}}}\rvert_{\bhs{0,1}(\midH)} 
		= \cK_{e^{-t\Delta^M_{\sH}}\Pi_{\sG_4}}, \quad
	\cK_{e^{-t\Delta_{\eps}}}\rvert_{\bhs{1,1}(\midH)} 
		= \cK_{e^{-t\Delta^M_{\sH}}\Pi_{\sG_2}}, \\
	\rho_{d,1}^{(m+1)}\cK_{e^{-t\Delta_{\eps}}}\rvert_{\bhs{d,1}(\midH)} = k_{e^{-t\Delta^{\bbH}_{\sH}}}, \quad
	\rho_{d,2}^{(m+1)}\cK_{e^{-t\Delta_{\eps}}}\rvert_{\bhs{d,2}(\midH)} = k_{e^{-t\Delta^{\bbH}_{\sH}}\Pi_{\sE_4^{\bbH}}}, \\
	\rho_{d,0}^{m}\cK_{e^{-t\Delta_{\eps}}}\rvert_{\bhs{d,0}(\midH)} = 
	\frac {\eps^{-(m+1)}}{(4\pi)^{m/2}}  \exp\Big( -\frac14 |\wt \omega|_{g_{\eps}(\zeta',z',1)}^2 \Big) \Id_{{}^{\sR}\Omega^p(X)} \; \bar\mu_R.
\end{gathered}
\end{equation*}
This finishes the proof of part (i) of Theorem \ref{thm:IntroThm}.

\subsection{Trace of the heat kernel}\label{Heat Trace Asymptotics} $ $

Having found a precise description of the structure of the heat kernel we now deduce the consequences for its trace.
Recall that, by Mercer's theorem \cite{Brislawn} \cite[Proposition 4.55]{Melrose:APS}, 
\begin{equation*}
	\Tr(e^{-t\Delta_{\eps}}) = \int_M \tr\lrpar{\beta_*\cK_{e^{-t\Delta_{\eps}}}\rest{\diag_M}}
\end{equation*}
where $\tr$ denotes the pointwise trace of $\Hom({}^{\sR}\Omega^*X)|_{\diag}$ and $\beta_*$ denotes the push-forward along the blow-down map from the heat space to $\bbR^+_t\times M^2 \times [0,1]_{\eps}.$ 
Equivalently, instead of pushing-forward the heat kernel and then restricting to the diagonal, we can directly restrict the heat kernel to the interior lift of the diagonal.

Let us start by considering differential forms of degree $p \in \{n, n+1\},$ as the other form degrees will be simpler.
\begin{figure}[ht]
    \def\svgwidth{\columnwidth}
    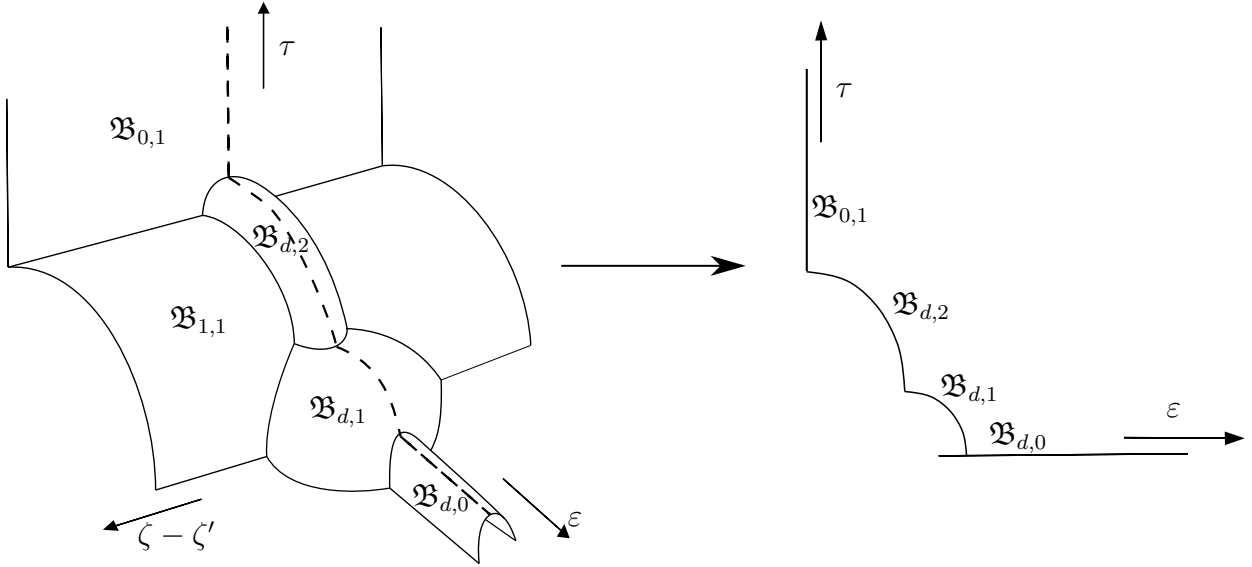

    \caption{The heat space $H_{\sR}^{\mid} X$ and corresponding diagonal $\Delta_{H_\sH^{\mid}X}$}
    \label{mid heatspace and diagonal}
\end{figure} 
The interior lift of the diagonal to $H^{\mid}_{\sR}X$ is shown in Figure \ref{mid heatspace and diagonal}, and can be identified with
\begin{equation*}
	\diag_{\midH}
	= M \times \Big[
	[\Rp_\tau \times [0,1]_\eps; 
	\{\tau=0=\eps\}];
	\{R=0, \Theta = \tfrac\pi2\}
	\Big],
\end{equation*}
with $R = \sqrt{\eps^2 + \tau^2}$ and $\Theta = \arctan(\frac\eps\tau)$ polar coordinates valid after the first blow-up.
We denote the boundary hypersurfaces of $\sT\sE^{\mid}$ with the same symbols we used for the boundary hypersurfaces of $\midH,$ thus
\begin{equation*}
\begin{aligned}
	\bhs{d,0}(\sT\sE^{\mid}) &= \bar{\beta^{-1}(\tau =0, \eps>0)}, \\
	\bhs{d,1}(\sT\sE^{\mid}) &= \text{ interior lift of $\{R=0\}$ after first blow-up},\\
	\bhs{d,2}(\sT\sE^{\mid}) &= \text{ interior lift of $\{R=0, \Theta = \tfrac\pi2\}$ after first blow-up},\\
	\bhs{0,1}(\sT\sE^{\mid}) &= \bar{\beta^{-1}(\eps =0, \tau>0)}.
\end{aligned}
\end{equation*}

Let us write this as $M \times \sT\sE^{\mid}$ and denote the natural projection off of $M$ by 
\begin{equation*}
	p_{\sT\sE^{\mid}}: \diag_{\midH}  \lra \sT\sE^\mid.
\end{equation*}
Mercer's theorem is then that
\begin{equation*}
	\Tr(e^{-t/\eps^2 \Delta_{\eps}}) = (p_{\sT\sE})_*\cK_{e^{-t/\eps^2\Delta_{\eps}}}\rest{\diag_{\midH}},
\end{equation*}
and this allows us to read off the asymptotics of the trace.

The analysis of the sub-Riemannian limit Hodge Laplacian on forms outside of middle degree is similar but simpler, with the interior lift of the diagonal given by
\begin{equation*}
	\diag_{\nmidH} = M \times [\bbR^+_{\tau} \times [0,1]_{\eps}; \{\tau=\eps=0\}].
\end{equation*}
Denote $[\bbR^+_{\tau} \times [0,1]_{\eps}; \{\tau=\eps=0\}]$ by $\sT\sE^{\nmid}$ and denote its boundary hypersurfaces by
\begin{equation*}
\begin{aligned}
	\bhs{d,0}(\sT\sE^{\nmid}) &= \bar{\beta^{-1}(\tau =0, \eps>0)}, \\
	\bhs{d,1}(\sT\sE^{\nmid}) &= \beta^{-1}(\tau=\eps=0),\\
	\bhs{0,1}(\sT\sE^{\nmid}) &= \bar{\beta^{-1}(\eps =0, \tau>0)}.
\end{aligned}
\end{equation*}
Finally recall from the discussion at the end of \S\ref{sec:CtHeat} that, although the heat kernel for differential form degrees $p \in \{n,n+1\}$ has log-terms in its asymptotic expansion, its trace does not.

\begin{theorem}\label{thm:AsympTrace}
Let $\Delta_{\eps}$ be the Hodge Laplacian of the sub-Riemannian limit metrics $g_{\eps}.$ 

For $p \notin\{n,n+1\},$ we have
\begin{equation*}
	\Tr(e^{-t\Delta_{\eps}}) 
	\in \rho_{d,1}^{-(m+1)}\rho_{d,0}^{-m} \CI(\sT\sE^{\nmid}), 
\end{equation*}
with a local expansion at  $\bhs{d,1}(\sT\sE^{\nmid})$ and $\bhs{d,0}(\sT\sE^{\nmid})$ (i.e., 
the coefficients of the expansion at $\bhs{d,1}(\sT\sE^{\nmid})$ and $\bhs{d,0}(\sT\sE^{\nmid})$ are given by universal polynomials in the corresponding symbol of $\Delta_{\eps}$) and the leading term at $\bhs{0,1}$ given by
\begin{equation*}
	\Tr(e^{-t\Delta_{\eps}^{\nmid}}) \rest{\bhs{0,1}} = \Tr(e^{-\tau\Delta_{\sH}}|_{\sE_2^*}).
\end{equation*}

For $p \in\{n,n+1\},$ we have
\begin{equation*}
	\Tr(e^{-t/\eps^2\Delta_{\eps}}) 
	\in (\rho_{d,1}\rho_{d,2})^{-(m+1)}\rho_{d,0}^{-m} \CI(\sT\sE^{\mid}), 
\end{equation*}
with a local expansion  at $\bhs{d,1}(\sT\sE^{\mid})$ and $\bhs{d,0}(\sT\sE^{\mid}),$ with expansion at $\bhs{d,2}(\sT\sE^{\mid})$
given by the sum of $\Tr(e^{-\tau\Delta_{\sH}}|_{\sG_2^*}) + \cO(\rho_{d,2})$ and a local expansion, and with leading term at $\bhs{0,1}$ given by
\begin{equation*}
	\Tr(e^{-t\Delta_{\eps}^{\nmid}}) \rest{\bhs{0,1}} = \Tr(e^{-\tau\Delta_{\sH}}|_{\sE_4^*}).
\end{equation*}
\end{theorem}

This establishes part (ii) of Theorem \ref{thm:IntroThm}. The local coefficients in the asymptotic expansion of the trace of the heat kernel at boundary hypersurfaces over $\{\eps=0\}$ are  integrals of universal polynomials in the curvature and torsion of the Tanno connection as we now discuss. 

\subsection{The heat invariants and Tanno's connection}\label{sec:Conn}
$ $

It is well-known \cite{Minakshisundaram-Pleijel} that the coefficients in the short-time asymptotic expansion of the trace of  heat kernel of the Hodge Laplacian are integrals of universal polynomials in the curvature of the metric and its covariant derivatives.
This is an easy consequence of the construction of the heat kernel in \S\ref{sec:EllHeat} since the terms in the Taylor expansion of the heat kernel at $\bhs{d,1}$ are obtained from the Taylor expansion of the symbol of the Lapacian and this is the Riemannian metric. The heat invariants arising from $\bhs{d,1}$ for the Hodge Laplacian undergoing a sub-Riemannian limit have a similar description as integrals of universal polynomials. In this case the universal variables are the coefficients of different powers of $\eps$ in the expansion of the Levi-Civita connection $\nabla^{\eps}$ as $\eps\to0$ and, as we know explain, these are the Tanno connection and its torsion.

Let 
\begin{equation*}
	W_0=\eps \cR, \quad
	W_1=X_1, \quad
	\ldots, \quad
	W_n=X_n, \quad
	W_{n+1}=Y_1, \quad
	\ldots, \quad
	W_{2n}=Y_n
\end{equation*}
be a Darboux frame for $(TM, g_{\eps})$ so that the only non-zero brackets are $[X_i,Y_i] = \cR.$
Define $\alpha(s,t)$ and $\wt s$ by
\begin{equation*}
	[W_s,W_t] = \alpha(s,t) \cR, \quad
	\wt s =
	\begin{cases}
	0 & \Mif s=0\\
	s+n & \Mif 0<s\leq n \\
	s-n & \Mif s>n
	\end{cases}
\end{equation*}
so that $J(W_i) = \alpha(i, \wt i) W_{\wt i}$ for all $i.$
Define the Christoffel symbols of this frame by
\begin{equation*}
	\nabla^{\eps}_{W_i}W_j = \Gamma_{ij}^k(\eps) W_k.
\end{equation*}
Then an easy computation using the Koszul formula shows that the non-zero Christoffel symbols are given by
\begin{equation*}
\begin{gathered}
	\Gamma_{ij}^k(\eps) = \Gamma_{ij}^k(1), \quad \Gamma_{ij}^0(\eps) = \tfrac{1}{2}(-\eps\cR g_{ij} +\tfrac{1}{\eps}\alpha(i,j)), \\
	\Gamma_{i0}^k(\eps) = \tfrac{1}{2}(g^{k\ell}(\eps\cR g_{\ell i}) - \tfrac{1}{\eps}\alpha(i,\wt i)g^{k\wt i} ) = \Gamma_{0i}^k(\eps),
\end{gathered}
\end{equation*}
where $i,j,k,\ell>0.$
It follows that as $\eps\to 0$ the Levi-Civita connection acting on differential forms splits into
\begin{equation*}
	\nabla^{\eps} = \eps \nabla^{1,0} + \nabla^{0,1} + \eps^{-1} \nabla^{-1,2}
\end{equation*}
(where the notation corresponds to the splitting of the exterior derivative in \S\ref{sec:Contact Complex}) in which
$\nabla^{-1,2}$ has the contribution of $J$ to the Christoffel symbols, and $\nabla^{1,0}$ has the derivative of the metric with respect to the Reeb vector field.
Thus $\nabla^{0,1}$ is the connection on $TM$ obtained from the Levi-Civita connection by removing the vertical contributions; 
this is known as the Tanno connection, see \cite[(3.1)]{Tanno}.
The parts of the connection forms that depend on $\eps$ as $\eps \to 0$ then make up the torsion forms of the Tanno connection.

\section{The $\eta$ invariant}

In this section we consider the $\eta$ invariant of the signature operator for sub-Riemannian limit metrics $g_{\eps}.$
This has been considered in three dimensions by Biquard-Herzlich-Rumin \cite[Theorem 1.4]{Biquard-Herzlich-Rumin} and in general dimension by Rumin \cite[\S7]{Rumin:Sub}.

The Hodge star induces a natural involution on the complexified differential forms on $M,$
\begin{equation*}
	\sI:\Omega_{\bbC}^*M \lra \Omega_{\bbC}^*M, \quad \sI^2=\Id.
\end{equation*}
The odd signature operator is
\begin{equation*}
	S = -i(d\sI+\sI d) = -i\sI (d-\delta) = -i(d-\delta)\sI,
\end{equation*}
its square is the Hodge Laplacian, 
and the $\eta$-invariant of an odd dimensional manifold is 
\begin{equation*}
	\eta(S) = \int_0^{\infty} \Tr(\sqrt t S e^{-t\Delta})\; \frac{dt}t.
\end{equation*}
Since $\sI$ maps forms of degree $p$ to forms of degree $m-p=2n+1-p,$
\begin{equation*}
	Se^{-t\Delta}: \Omega^p_{\bbC}M \lra \Omega^{2n+2-p}_{\bbC}M \oplus \Omega^{2n-p}_{\bbC}(M),
\end{equation*}
and so the only forms that can contribute to the trace are those of the middle degrees $p\in \{n,n+1\},$ i.e., we have
\begin{equation*}
	\eta(S) = \int_0^{\infty} \Tr(\Pi^{\mid}\sqrt t S e^{-t\Delta}\Pi^{\mid})\; \frac{dt}t,
\end{equation*}
with $\Pi^{\mid}$ the projection onto middle degree forms.

\begin{theorem}
Let $g_{\eps}$ be a sub-Riemannian limit family of metrics on a contact manfiold and let $S_{\eps}$ denote its odd signature operator.
The Schwartz kernel of the operator $\Pi^{\emph{mid}}\frac{\sqrt t}{\eps} S_{\eps} e^{-\frac t{\eps^2}\Delta_{\eps}}\Pi^{\emph{mid}}$ (pulled-back to $H_{\sR}^{\emph{mid}}X$) is an $\cI$-smooth right density
\begin{equation*}
	\beta^*\left(\Pi^{\emph{mid}}\frac{\sqrt t}{\eps} S_{\eps} e^{-\frac t{\eps^2}\Delta_{\eps}}\Pi^{\emph{mid}}\right) 
	\in \sA_{\phg}^{\cE_{\eta}}(H_{\sR}^{\emph{mid}}X; \Hom(\Lambda^*({}^{\sR}T^*X)) \beta_R^*\Lambda)
\end{equation*}
with the same index sets as the middle degree heat kernel \eqref{eq:HeatKerMid}.

The $\eta$-invariant of $S_{\eps}$ is an $\cI$-smooth function on $[0,1]_{\eps},$ 
\begin{multline*}
	\eta(S_{\eps})
	= \int_0^{\infty} \Tr(\Pi^{\emph{mid}}\sqrt t S_{\eps} e^{-t\Delta_{\eps}}\Pi^{\emph{mid}})\; \frac{dt}t
	=
	\int_0^{\infty} \Tr\left(\Pi^{\emph{mid}}\frac{\sqrt u}{\eps} S_{\eps} e^{-\frac u{\eps^2}\Delta_{\eps}}\Pi^{\emph{mid}}\right) \; \frac{du}u \\
	\in \sA_{\phg}^{E_{\eta}}([0,1]) 
	\Mwith E_{\eta}(\{\eps=0\}) = (-m+\bbN_0) \bar\cup (-m+1+\bbN_0)
\end{multline*}
and we have
\begin{multline*}
	\FP_{\eps=0}
	\eta(S_{\eps})
	= \int_0^{\infty} \Tr(\Pi^{\emph{mid}}\sqrt t (D_{\sH}*+*D_{\sH}) e^{-t\Delta_{\sH}\Pi_{\sE_4}}\Pi^{\emph{mid}})\; \frac{dt}t \\
	+ \int_0^{\infty} \Tr(\Pi^{\emph{mid}}\sqrt t (d_{\sH}*+*d_{\sH}) e^{-t\Delta_{\sH}\Pi_{\sE_2\setminus\sE_4}}\Pi^{\emph{mid}})\; \frac{dt}t 
	+ \text{local}
	\\
	= \eta_{\emph{contact}}
	+ \text{local}
\end{multline*}
where the final term is an integral of a universal polynomial in the torsion and curvature of the Tanno connection and their covariant derivatives.
\end{theorem}

\begin{proof}
We can write the differential operator $\beta_L^*(\Pi^{\mid}\frac{\sqrt t}{\eps} S_{\eps}\Pi^{\mid})$ using vector fields on $\midH$ that are tangent to every boundary hypersurface (except $\bhs{0,1},$ and $\bhs{1,1}$), and with coefficients that are smooth on $\midH.$ Applying such an operator to $e^{-t\Delta_{\eps}^{\mid}}$ yields a function that is $\cI$-smooth at all boundary hypersurfaces but $\bhs{0,1}$ and $\bhs{1,1}$ with the same index sets.

The same argument applies to $\beta_L^*(\eps^2\Pi^{\mid}\frac{\sqrt t}{\eps} S_{\eps}\Pi^{\mid})$ at $\bhs{0,1},$ as it is made up of vector fields tangent to that face, but not to $\beta_L^*(\Pi^{\mid}\frac{\sqrt t}{\eps} S_{\eps}\Pi^{\mid})$ since it has coefficients singular at this face.
However, we constructed the Taylor expansion of $e^{-t\Delta_{\eps}^{\mid}}$ at $\bhs{0,1}$ so that it was valued in $\sE_4,$ and hence it follows from \S \ref{sec:sRLim} that $\beta^*(\Pi^{\mid}\frac{\sqrt t}{\eps} S_{\eps}e^{-t\Delta_{\eps}^{\mid}})$ is $\cI$-smooth at $\bhs{d,2}$ with the same (smooth) index set as $e^{-t\Delta_{\eps}^{\mid}},$ and leading term
\begin{equation*}
\beta_L^*(\Pi^{\mid}\sqrt t (D_{\sH}*+*D_{\sH})) \cK_{e^{-t\Delta_{\sH}}}.
\end{equation*}
Similarly, we have a contribution at $\bhs{d,2}$ given by
\[\Tr(\Pi^{\emph{mid}}\sqrt t (d_{\sH}*+*d_{\sH}) e^{-t\Delta_{\sH}\Pi_{\sG_2}}\Pi^{\emph{mid}}) \]
however we can observe that this will not contribute to the renormalized $\eta$-invariant since this trace is identically zero. Since $d_\sH$ vanishes on $\Omega_\sH^nM$, this signature operator $(*d_\sH+d_\sH *)$ maps $\Omega_\sH^\mid M$ to $\Omega_\sH^\nmid M$.

The rest of the theorem follows by applying Melrose's push-forward theorem as above.
\end{proof}

\section{Determinant of the Hodge Laplacian} \label{sec:Det}

Ray-Singer \cite{Ray-Singer:AT} defined the determinant of a Laplacian by generalizing the relation, valid for any finite set of positive numbers $\{\mu_i\},$
\begin{equation*}
	\log \prod \mu_i 
	= -\pa_{s}|_{s=0}\sum \mu^{-s} 
	= -\pa_{s}|_{s=0}\lrpar{\frac1{\Gamma(s)}\int_0^{\infty} t^s \sum e^{-t\mu_i} \; \frac{dt}t}
	= -\pa_s|_{s=0}\zeta_{\{\mu_i\}}(s),
\end{equation*}
in a way that we now briefly review.

Let $e^{-tA}$ be the heat kernel of a positive operator, such as an elliptic Laplace-type operator or a contact Laplacian, whose trace satisfies
\begin{equation*}
	\Tr(e^{-tA}) \sim t^{-n/\ell}\sum_{k\geq 0} a_k(A) t^{k/\ell} \quad \Mas t \to 0, \quad
	\Tr(e^{-tA}-\Pi_{\Ker A}) \to 0 \text{ exponentially } \Mas t \to \infty,
\end{equation*}
for some $n,\ell \in \bbN,$ and with $\dim\Ker A<\infty.$
The zeta function of $A$ is defined, for $\Re(s)>n/\ell,$ by
\begin{equation*}
	\zeta_A(s) = \frac1{\Gamma(s)} \int_0^{\infty} t^s \Tr(e^{-tA} - \Pi_{\Ker A}) \; \frac{dt}t.
\end{equation*}
This is a holomorphic function of $s$ on this half-plane and the short-time asymptotic expansion of $\Tr(e^{-tA})$ induces a meromorphic continuation
(of $\Gamma(s)\zeta_A(s)$ and hence of $\zeta_A(s)$) with potential poles at $s \in \{(n-k)/\ell: k \in \bbN_0\}.$
We define (cf. \cite[\S3]{Hassell:AT})
\begin{equation*}
	\sideset{^R}{_0^{\infty}}\int \Tr(e^{-tA}) \; \frac{dt}t = \FP_{s=0} \int_0^{1} t^s \Tr(e^{-tA}) \; \frac{dt}t 
	+\FP_{s=0} \int_1^{\infty} t^s \Tr(e^{-tA}) \; \frac{dt}t 
\end{equation*}
and then note that $\zeta_A(s),$ near $s=0,$ is given by
\begin{equation*}
	\Big(s + \gamma s^2 + \cO(s^3)\Big)
	\lrpar{ \frac{a_{n}-\dim\Ker A}s + \sideset{^R}{_0^{\infty}}\int \Tr(e^{-tA}) \; \frac{dt}t + \cO(s) },
\end{equation*}
where $\gamma$ is the Euler-Mascheroni constant.
Thus (the meromorphic continuation of) $\zeta_A(s)$ is regular at $s=0,$ its derivative is equal to
\begin{equation*}
	\zeta_A'(0) 
	= \gamma (a_n - \dim\ker A) + \sideset{^R}{_0^{\infty}}\int \Tr(e^{-tA}) \; \frac{dt}t,
\end{equation*}
and we define this to be $-\log\det A.$\\

Let us consider the Hodge Laplacian for sub-Riemannian limit metrics.
First outside of middle degrees, $p \notin \{n,n+1\},$ we have that
\begin{equation*}
	\Tr(e^{-t\Delta_{\eps}}-\Pi_{\Ker \Delta_{\eps}}) \to 0 \text{ exponentially } \Mas t \to \infty,
\end{equation*}
at a rate independent of $\eps$ by \cite[Theorem 7.1]{Rumin:Sub}.  
Secondly, in this case $a_n=0$ since $\ell=2,$ $n=m,$ and $a_{k}(\Delta_{\eps})=0$ for $k$ odd.
Finally the push-forward theorem for renormalized integrals \cite[page 128]{Hassell-Mazzeo-Melrose:Surgery}, \cite[Lemma 11.1]{Albin-Rochon-Sher:FibCusps} allows us to conclude from Theorem \ref{thm:AsympTrace} that
\begin{equation*}
	-\FP_{\eps=0}\log\det \Delta_{\eps}
	= - \gamma b_p + \sideset{^R}{_0^{\infty}}\int \Tr(e^{-t\Delta_{\sH}}) \; \frac{dt}t
	+ \sideset{^R}{_0^{\infty}}\int A_{m+1} \; \frac{d\sigma}\sigma,
\end{equation*}
where $b_p$ is the $p$-th Betti number of $M,$ $\sigma$ is a rescaled time-variable on $\bhs{d,1}(\nmidH)$ (e.g., $\sigma =\frac{\sqrt{t}}{\eps}$) and $A_{m+1}$ is, with notation from Theorem \ref{thm:IntroThm}, the constant term in the expansion of the trace at $\bhs{d,1}(\nmidH).$
Comparing with $\log\det\Delta_{\sH}$ we have, for forms of degree $p \notin\{n,n+1\},$
\begin{equation*}
	\FP_{\eps=0}\log\det \Delta_{\eps,p}
	= \log\det\Delta_{\sH,p} + \gamma a_{m+1}(\Delta_{\sH,p})
	- \sideset{^R}{_0^{\infty}}\int A_{m+1} \; \frac{d\sigma}\sigma.
\end{equation*}

Next let us consider middle degree forms, $p \in \{n,n+1\},$ starting with $p=n.$
Arguing as above we find that
\begin{align*}
	-\FP_{\eps=0}\log\det \eps^{-2}\Delta_{\eps,n} &
	= - \gamma b_n 
	+ \sideset{^R}{_0^{\infty}}\int \Tr(e^{-t(D_{\sH}+D_{\sH}^*)^2}\rest{\sE_4^n}) \; \frac{dt}t \\
	&+ \sideset{^R}{_0^\infty}\int \Tr(e^{-t(d_{\sH}+\delta_{\sH})^2}\rest{\sG_2^n})\; \frac{dt}{t} 
	+ \sideset{^R}{_0^{\infty}}\int B_{m+1,n} \; \frac{d\sigma}\sigma
	+ \sideset{^R}{_0^{\infty}}\int A_{m+1,n} \; \frac{d\sigma'}{\sigma'},
\end{align*}
where $\sigma,$ and $\sigma'$ are rescaled time variables (e.g., $\frac{\sqrt t}{\eps}$ and $\frac{\sqrt t}{\eps^2},$ respectively) and, with notation from Theorem \ref{thm:IntroThm}, $A_{m+1,n}$ and $B_{m+1,n}$ are the constant terms in the local part of the expansion of the trace at $\bhs{d,1}(\midH)$ and $\bhs{d,2}(\midH),$ respectively.

Recall \eqref{eq:BransonTrick}, that we have 
\begin{equation*}
	\Tr(e^{-t\Delta_{\sH}}\rest{\sE_4^n}) 
	= \Tr_{\Omega^n_{\sH}M}(e^{-t\Delta_{\sH}}) 
	+ \sum_{k=1}^n (-1)^k \Tr_{\Omega^{n-k}_{\sH}M}(e^{-t\Delta_{\sH}^2}),
\end{equation*}
and similarly,
\begin{multline*}
	\Tr(e^{-t(d_{\sH}+\delta_{\sH})^2}\rest{\sG_2^n})=\Tr_{\Omega_\sH^nM} (e^{-td_\sH\delta_\sH})=\Tr_{\Omega_\sH^{n-1}M} (e^{-t\delta_\sH d_\sH}) \\
	= \Tr_{\Omega_\sH^{n-1}M} (e^{-t\Delta_\sH}) - \Tr_{\Omega_\sH^{n-2}M} (e^{-td_\sH\delta_\sH}) = \ldots = \sum_{k=1}^{n} (-1)^{k-1} \Tr_{\Omega_\sH^{n-k}M} (e^{-t\Delta_\sH})
\end{multline*}
so that
\begin{multline*}
	- \gamma b_n 
	+ \sideset{^R}{_0^{\infty}}\int \Tr(e^{-t(D_{\sH}+D_{\sH}^*)^2}\rest{\sE_4^n}) \; \frac{dt}t 
	+ \sideset{^R}{_0^\infty}\int \Tr(e^{-t(d_{\sH}+\delta_{\sH})^2}\rest{\sG_2^n})\; \frac{dt}{t} \\
	=  \sideset{^R}{_0^{\infty}}\int \Tr_{\Omega^n_{\sH}M}(e^{-t\Delta_{\sH}}) \; \frac{dt}t 
	+ \sum_{k=1}^n (-1)^k \sideset{^R}{_0^{\infty}}\int \Tr_{\Omega^{n-k}_{\sH}M}(e^{-t\Delta_{\sH}^2}) - \Tr_{\Omega_\sH^{n-k}M} (e^{-t\Delta_\sH}) \; \frac{dt}t, \\
	= \zeta_{\Delta_{\sH,n}}'(0) - \gamma(a_{m+1}(\Delta_{\sH,n})) \\
	+ \sum_{k=1}^n (-1)^k \lrpar{
	\zeta_{\Delta_{\sH,n-k}^2}'(0) - \zeta_{\Delta_{\sH,n-k}}'(0) - \gamma(a_{m+1}(\Delta_{\sH,n-k}^2) - a_{m+1}(\Delta_{\sH,n-k}^2)}
\end{multline*}
and since $\zeta_{A^2}(s)=\zeta_A(2s)$ we have
\begin{multline*}
	\FP_{\eps=0} \log\det \Delta_{\eps,n} 
	= \log\det \Delta_{\sH,n} + \gamma a_{m+1}(\Delta_{\sH,n}) \\
	- \sideset{^R}{_0^{\infty}}\int B_{m+1,n} \; \frac{d\sigma}\sigma
	- \sideset{^R}{_0^{\infty}}\int A_{m+1,n} \; \frac{d\sigma'}{\sigma'}, \\
	+ \sum_{k=1}^n (-1)^k \log\det \Delta_{\sH,n-k} + \gamma a_{m+1}(\Delta_{\sH,n-k})
\end{multline*}

The corresponding analysis for $p = n+1$ yields
\begin{multline*}
	\FP_{\eps=0}\log\det \eps^{-2}\Delta_{\eps,n+1}
	= \log\det \Delta_{\sH,n+1} + \gamma(a_{m+1}(\Delta_{\sH,n+1})) \\
	- \sideset{^R}{_0^{\infty}}\int B_{m+1,n+1} \; \frac{d\sigma}\sigma
	- \sideset{^R}{_0^{\infty}}\int A_{m+1,n+1} \; \frac{d\sigma'}{\sigma'}, \\
	+ \sum_{k=1}^n (-1)^k \log\det \Delta_{\sH,n+k+1} + \gamma a_{m+1}(\Delta_{\sH,n+k+1})
\end{multline*}
$ $

Finally, let us discuss what this means for analytic torsion. The constructions above are essentially unchanged by allowing the Laplacians to act differential forms twisted by a flat bundle $F \lra M.$ Let $g_F$ be a bundle metric on $F.$
Our convention, following \cite{Ray-Singer:AT}, is to set
\begin{equation*}
	\log \mathrm{AT}(M,g_{\eps},F,g_F) = \sum_{p=0}^m (-1)^p \; p\; \zeta_{\Delta_{\eps,p}}'(0).
\end{equation*}
This is independent of the metric $g_{\eps},$ and in particular independent of $\eps,$ if the bundle $F$ is acyclic\footnote{A flat bundle is said to be acyclic if $\tH^j(M;F)=0$ for all $j.$} 
and its holonomy is orthogonal \cite{Cheeger:AT, Muller:AT} or unimodular \cite{Muller:Uni}.
For flat bundles that are not acyclic we can remove the dependence on the metric by assigning to each basis $\{\mu_j^q\}$ of $\Ker \Delta_{g_{\eps},q}$ the number
\begin{equation*}
	\bar{\log \mathrm{AT}}(M,\{\mu_j^q\}, F) = \log\mathrm{AT}(M,g_{\eps},F,g_F) - \log \Big( \prod_{q=0}^{m} [\mu^q|\omega^q]^{(-1)^q} \Big),
\end{equation*}
where $\{\omega^q\}$ denotes an orthonormal basis of harmonic forms and $[\mu^q|\omega^q] = |\det W^q|$ with $W^q$ the change-of-basis matrix satisfying
\begin{equation*}
	\mu_i^q = (W^q)_i^j \omega_j^q.
\end{equation*}
Note that, by e.g., \cite[Theorem 7.1]{Rumin:Sub}, an orthonormal basis of harmonic forms for $\eps>0$ converges to an orthonormal basis of harmonic forms for the Rumin complex, so we only need to understand the asymptotics of
$\log \mathrm{AT}_{\eps}.$

We start by noting that, for any $\eps>0$ and $\Re s>m/2,$
\begin{multline*}
	\zeta_{\eps^{-2}\Delta_{\eps}}(s)
	= \frac1{\Gamma(s)} \int_0^{\infty} t^s \Tr(e^{-t\eps^{-2}\Delta_{\eps}}- \Pi_{\Ker \Delta_{\eps}}) \; \frac{dt}t\\
	\xrightarrow[\frac{du}u=\frac{dt}t]{\phantom{xx}u=\eps^{-2}t\phantom{xx}}
	\frac{\eps^{2s}}{\Gamma(s)} \int_0^{\infty} u^s \Tr(e^{-u\Delta_{\eps}}- \Pi_{\Ker \Delta_{\eps}}) \; \frac{du}u
	= \eps^{2s}\zeta_{\Delta_{\eps}}(s),
\end{multline*}
and hence we have
\begin{equation*}
	\zeta_{\eps^{-2}\Delta_{\eps}}'(0) = 2(\log \eps) \zeta_{\Delta_{\eps}}(0) + \zeta_{\Delta_{\eps}}'(0), 
	\Mand
	\FP_{\eps=0} \log \det \eps^{-2}\Delta_{\eps} = \FP_{\eps=0} \log \det \Delta_{\eps}.
\end{equation*}

Thus the finite part of analytic torsion is given by
\begin{multline*}
	\FP_{\eps=0} \log\mathrm{AT}(M, g_{\eps}, F, g_F)
	= \frac12\sum_{p=0}^{2n+1} (-1)^p p \FP_{\eps=0} \zeta_{\Delta_{\eps,p}}'(0) \\
	=  \frac{1}{2} \sum_{p\neq n,n+1}  (-1)^p p \FP_{\eps=0} \zeta_{\Delta_{\eps,p}}'(0) 
	+\frac{1}{2} \sum_{p=n}^{n+1} (-1)^p p \FP_{\eps=0}\zeta_{\eps^{-2}\Delta_{\eps,p}}'(0)  \\
	= \frac12\sum_{p\neq n,n+1} (-1)^p p \Big(
	\zeta_{\Delta_{\sH,p}}'(0)-\gamma a_{m+1}(\Delta_{\sH,p}) 
	+ \sideset{^R}{_0^{\infty}}\int A_{m+1,p} \; \frac{d\sigma}\sigma \Big) \\
	+ \frac12\sum_{p=n}^{n+1}(-1)^p p \Big(
	\zeta_{\Delta_{\sH,p}}'(0)-\gamma a_{m+1}(\Delta_{\sH,p}) 
	+ \sideset{^R}{_0^{\infty}}\int B_{m+1,p} \; \frac{d\sigma}\sigma
	+ \sideset{^R}{_0^{\infty}}\int A_{m+1,p} \; \frac{d\sigma'}{\sigma'} \Big)    \\ 
	+ (-1)^n \frac{n}{2}\sum_{k=1}^n (-1)^k \zeta_{\Delta_{\sH,n-k}}'(0) - \gamma a_{m+1}(\Delta_{\sH,n-k}) \\
	+(-1)^{n+1}\frac{n+1}{2} \sum_{k=1}^n (-1)^k \zeta_{\Delta_{\sH,n+k+1}}'(0) - \gamma a_{m+1}(\Delta_{\sH,n+k+1})\\
	=\frac12\sum_{p\neq n,n+1} (-1)^p p\sideset{^R}{_0^{\infty}}\int A_{m+1,p} \; \frac{d\sigma}\sigma 
	+ \frac12\sum_{p=n}^{n+1} (-1)^p p\Big(\sideset{^R}{_0^{\infty}}\int B_{m+1,p} \; \frac{d\sigma}\sigma
	+ \sideset{^R}{_0^{\infty}}\int A_{m+1,p} \; \frac{d\sigma'}{\sigma'} \Big)  \\
	+ \frac12\sum_{p=0}^{n-1} (-1)^p \Big( (p + n)\zeta_{\Delta_{\sH,p}}'(0) - \gamma((p+ n)a_{m+1}(\Delta_{\sH,p}) )\Big)\\
	+ \frac12\sum_{p=n}^{n+1} (-1)^p p\Big(\zeta_{\Delta_{\sH,p}}'(0) - \gamma ((p) a_{m+1}(\Delta_{\sH,p}) ) \Big)\\
	+  \frac12\sum_{p=n+2}^{2n+1} (-1)^p \Big((p - (n+1))\zeta_{\Delta_{\sH,p}}'(0) - \gamma((p- (n+1))a_{m+1}(\Delta_{\sH,p}))\Big) \\
	=\frac12\sum_{p\neq n,n+1} (-1)^p p\sideset{^R}{_0^{\infty}}\int A_{m+1,p} \; \frac{d\sigma}\sigma 
	+ \frac12\sum_{p=n}^{n+1} (-1)^p p\Big(\sideset{^R}{_0^{\infty}}\int B_{m+1,p} \; \frac{d\sigma}\sigma
	+ \sideset{^R}{_0^{\infty}}\int A_{m+1,p} \; \frac{d\sigma'}{\sigma'} \Big)  \\
	+ \log\mathrm{AT}_{\sH}(M, g_{\sH}, F, g_F).
\end{multline*}
We make one final remark: the torsion, $\log \bar{AT}_\sH(M,\{\mu^*\},F)$, associated to the weight function $\wt w(p)$ arising above, and the weight function $w(p)$ originally given in \cite{Rumin-Seshadri},
\[ \wt w(p) = \begin{cases} p+n & p< n\\ p & p\in \{n,n+1\} \\ p -(n+1) & p>n+1  \end{cases}, \quad w(p) = \begin{cases} p & p\leq n \\ p+1 &  n+1 \leq p \end{cases}\]
coincide whenever the  complex twisted by $F$ satisfies Poincar\'e duality, i.e. arises from a unitary representation. This was the only case considered in \cite{Rumin-Seshadri}, however in this section we are working with arbitrary flat bundles and thus take this to extend their original definition.
This establishes the proof of Corollary \ref{cor:IntAT}.

\begin{remark}
When comparing with the definition of Rumin-Seshadri, note that their convention is to use $\Delta_{\sH,p}^2$ for $p \notin\{n,n+1\}$ and that their definition of analytic torsion is the multiplicative inverse of ours.
\end{remark}

Finally we establish the relation between Kitaoka's torsion and that of Rumin-Seshadri and justify \eqref{eq:IntKitaokaRS}. 

On a $2n+1$-dimensional contact manifold, Kitaoka modifies the definition of the Rumin complex $(\sK^*,d_\sK)$ via a new definition of the differential: $d_\sK^p= a_pd_{\sH}^p= \tfrac{1}{\sqrt{|n-p|}}d_{\sH}^p$ for $p\neq \{n,n+1\}$. We can relate sub-Laplacians $\Delta_{\sH,p}$ of the original complex to those of the modified complex as follows
\begin{equation*}
\begin{gathered}
	\Delta_{\sK,p} = (\delta_\sK^{p+1} d_\sK^p)^2+(d_\sK^{p-1} \delta_\sK^{p})^2 = \tfrac{1}{(n-p)^2}(\delta_\sH d_\sH)^2 + \tfrac{1}{(n-p+1)^2}(d_\sH\delta_\sH)^2 , \Mfor p \notin \{n, n+1\}, \\
	\Delta_{\sK,n}= \Delta_{\sH,n} \Mand 
	\Delta_{\sK,n+1}= \Delta_{\sH,n+1}.
\end{gathered}\end{equation*}
Now, since $d_\sH$ and $\delta_\sH$ commute with $\Delta_\sH$, we can use the equivalence between the non-zero eigenvalues of $\delta_\sH^{p+1} d_\sH^{p}$ and the non-zero eigenvalues of $d_\sH^p \delta_\sH^{p+1}$ to conclude that
\[ \zeta (\delta_\sH^{p+1} d_\sH^{p})(s) = \zeta (d_\sH^p \delta_\sH^{p+1})(s)   \]
where $\zeta(A)(s)=\sum\limits_{0\neq \lambda_\ell\in \sigma(A)}\lambda_\ell^{-s}$ (with eigenvalues repeated with multiplicity)\footnote{This differs from the definition for $\zeta(A)(s)$ of \cite{Rumin-Seshadri}, but only by a cohomological term, which is $s$-independent. Thus it does not contribute to $\pa_s\cK_\sH$.}. In particular, since 
\[ \zeta(\Delta_{\sK,p})(s) = \zeta (\tfrac{1}{(n-p)^2}(\delta_\sH^{p+1} d_\sH^{p})^2)(s) + \zeta (\tfrac{1}{(n-p+1)^2} (d_\sH^{p-1} \delta_\sH^{p})^2)(s)    \]
we have
\begin{align*}
\zeta_{AT,\sK}(s) &= \sum_{p=0}^n (-1)^{p+1}(n+1-p)\zeta(\Delta_{\sK,p})(s)  \\
& =  -n \zeta (\tfrac{1}{n^2}\Delta_\sH^0)(s) \\
& +  (n-1)\zeta (\tfrac{1}{(n-1)^2} (\delta_\sH^{2}d_\sH^1)^2)(s) + (n-1)\zeta (\tfrac{1}{n^2}\Delta_\sH^0)(s) \\
& -  (n-2)\zeta (\tfrac{1}{(n-2)^2} (\delta_\sH^{3}d_\sH^2)^2)(s) - (n-1)\zeta (\tfrac{1}{(n-1)^2} (\delta_\sH^{2}d_\sH^1)^2)(s) \\
& + \ldots  \\
&+ (-1)^{n+1} \zeta (D_{\sH}^* D_{\sH})(s) + (-1)^{n+1} \zeta (d_\sH^{n-1}\delta_\sH^n)(s),
\end{align*}
so after grouping together the different terms with the same half-Laplacian,
\begin{equation*}
\begin{gathered}
	\zeta_{AT,\sK}(s)= -\zeta (\tfrac{1}{n^2}(\delta_\sH^1d_\sH^0)^2)(s) + \zeta (\tfrac{1}{(n-1)^2}(\delta_\sH^2d_\sH^1)^2)(s) -\ldots + (-1)^{n+1} \zeta (D_{\sH}^* D_{\sH})(s)  \\
	= -n^{2s}\zeta ((\delta_\sH^1d_\sH^0)^2)(s) + (n-1)^{2s}\zeta ((\delta_\sH^2d_\sH^1)^2)(s) -\ldots + (-1)^{n+1} \zeta (D_{\sH}^* D_{\sH})(s) \\
	= \sum_{p=0}^{n-1} (-1)^{p+1} (n-p)^{2s}\zeta (\delta_\sH^{p+1}d_\sH^p)(s) +(-1)^{n+1}\zeta (D_{\sH}^* D_{\sH})(s)  \\
	= \sum_{p=0}^{n-1} (-1)^{p+1} (n-p)^{2s}\zeta (d_\sH^p\delta_\sH^{p+1})(s) +(-1)^{n+1}\zeta (D_{\sH}^* D_{\sH})(s).
\end{gathered}\end{equation*}
If we now use Branson's observation \eqref{eq:BransonTrick} in the form
\[ \zeta (d_{\sH}^p\delta_{\sH}^{p+1})(s) = \sum_{j=0}^p (-1)^j \zeta (\Delta_\sH^{p-j})(s)  \]
we can rewrite this in terms of the Rumin Hodge Laplacians as
\begin{align*}  
\zeta_{AT,\sK}(s) &=\sum_{p=0}^{n-1} (-1)^{p+1} (n-p)^{2s}  \sum_{j=0}^p (-1)^j \zeta (\Delta_\sH^{p-j})(s) + (-1)^{n+1} \zeta (D_{\sH}^* D_{\sH})(s)  \\
 & = \sum_{p=0}^{n-1} (-1)^{p+1}\left( \sum_{j=1}^{n-p} j^{2s} \right)  \zeta (\Delta_\sH^p)(s) + (-1)^{n+1} \zeta (D_{\sH}^* D_{\sH})(s).
\end{align*}

For comparison, Rumin-Seshadri's torsion function is equal to
\begin{align*}
\zeta_{AT,\sH}(s) &= -\zeta ((\delta_\sH^1d_\sH^0)^2)(s) + \zeta ((\delta_\sH^2d_\sH^1)^2)(s) -\ldots + (-1)^{n+1} \zeta (D_{\sH}^* D_{\sH})(s)  \\
& =\sum_{p=0}^{n-1}(-1)^{p+1}(n-p)\zeta (d_\sH^p \delta_\sH^{p+1})(s) +(-1)^{n+1}\zeta (D_{\sH}^* D_{\sH})(s) \\
&= \sum_{p=0}^{n-1} (-1)^{p+1}\left( \sum_{j=1}^{n-p} j \right)  \zeta (\Delta_\sH^k)(s) + (-1)^{n+1} \zeta (D_{\sH}^* D_{\sH})(s) .
\end{align*}
Taking derivatives at $s=0$ we find that
\begin{multline}\label{eq:KitaokaRS}
	\log AT_\sK(M,g_\sH) = \log AT_{\sH}(M,g_\sH) + 2\sum_{p=0}^{n-1} (-1)^{p+1} \log[(n-p)!]\zeta(\Delta_\sH^p)(0) \\
	 = \log AT_{\sH}(M,g_\sH) + \text{local}
\end{multline}
which shows that these two definitions of analytic torsion of the Rumin complex differ by a local term.

\end{document}